\documentclass[draft]{article}
\usepackage[normalem]{ulem}
\usepackage{amsmath,amsfonts,amsthm,amssymb,amscd,cancel,bbm,mathabx}
 \usepackage{enumitem}
\usepackage{verbatim}
\usepackage[mathscr]{euscript}
\renewcommand{\Re}{\mathop{\rm Re}\nolimits}
\renewcommand{\Im}{\mathop{\rm Im}\nolimits}
\newcommand{\abs}[1]{|#1|}

\newcommand{\at}[1]{\big\vert\sb{#1}}

\theoremstyle{plain} \newtheorem{theorem}{Theorem}[section]
\newtheorem{lemma}[theorem]{Lemma}
\newtheorem{proposition}[theorem]{Proposition}
\newtheorem{corollary}[theorem]{Corollary} \theoremstyle{definition}
\newtheorem{definition}[theorem]{Definition} \theoremstyle{remark}
\newtheorem{remark}[theorem]{Remark}

\newtheorem{claim}[theorem]{Claim}
 
\newcommand{\R}{\mathbb{R}}

 \newcommand{\C}{\mathbb{C}}
\newcommand{\Z}{\mathbb{Z}}
 \newcommand{\U}{\mathcal{U}}
\newcommand{\N}{\mathbb{N}}

\def\({\left(}
\def\){\right)}
\def\<{\left\langle}
\def\>{\right\rangle}

 \def \im{\mbox{$\imath\hskip -4pt\imath\,$}}

\newcommand{\resto}{{\mathcal R}}
\newcommand{\s}{\mathfrak{ s}}

\def\uno{{\kern+.3em {\rm 1} \kern -.22em {\rm l}}}

\numberwithin{equation}{section}

\begin{document}

 \title{ On asymptotic stability of ground states of some systems of nonlinear Schr\"odinger equations}

\author{
{\sc Andrew Comech}
\\
{\it\small Texas A\&M University, College Station, Texas 77843, USA}
{\small and}
{\it\small IITP, Moscow 127051, Russia}
\\~\\
{\sc Scipio Cuccagna}
\\
{\it\small
University of Trieste, Trieste 34127, Italy, scuccagna@units.it}
}

 \maketitle

\begin{abstract} We extend to a specific class of systems of nonlinear Schr\"odinger equations (NLS)   the theory
of asymptotic stability of ground states already proved for the scalar NLS. Here the key point is the
choice of an adequate system of modulation coordinates and the novelty, compared to the scalar NLS, is
the fact that the group of symmetries of the system is non-commutative.
\end{abstract}

\section{Introduction}

In this paper we consider the system of coupled
nonlinear Schr\"odinger equations,
\begin{equation}\label{Eq:NLS}
\left\{\begin{array}{l}
\mathbbm{i}\sigma_3 \dot u + \Delta u -\beta (|u |^2) u=0,
\\[1ex]
u(0,x)=u_0(x) \in \C^2,\qquad x\in \R^3,
\end{array}
\right.
\end{equation}
where $ \mathbbm{i}$ is the imaginary unit and
the Pauli matrices are given by
\begin{equation}
\label{pauli-matrices}
\begin{aligned} &
\sigma_1=\begin{pmatrix} 0 &
1 \\
1 & 0
 \end{pmatrix} \,,
\qquad
\sigma_2= \begin{pmatrix} 0 &
-\mathbbm{i} \\
 \mathbbm{i} & 0
 \end{pmatrix} \,,
\qquad
\sigma_3=\begin{pmatrix} 1 &
0 \\
0 & -1
 \end{pmatrix}. \end{aligned}
\end{equation}
We assume that $\beta $ satisfies the following two hypotheses, which guarantee the local well-posedness of \eqref{Eq:NLS} in $H^1(\R^3, \C^2)$:
\begin{itemize}
\item[(H1)]
$\beta (0)=0$, $\beta\in C^\infty(\R,\R)$;
\item[(H2)] there exists
$\alpha \in(1,5)$ such that for every $k\in\N_{0}$
there is a fixed $C_k$ with $$\left|
\frac{d^k}{dv^k}\beta(v^2)\right|\le C_k \abs{v}^{\alpha-k-1}
\qquad\text{for}\ \ v\in\R,\quad \abs{v}\ge 1.
$$ \end{itemize}

\noindent We recall that under further hypotheses, there exist ground state solutions of the scalar NLS
\begin{equation}\label{Eq:scalNLS}
 \mathbbm{i} \dot u + \Delta u -\beta (|u |^2) u=0\,, \qquad
u(t,x)\vert_{t=0}=u_0(x) \in \C,
\qquad x\in \R^3
\end{equation}
 in $H^1 (\R^3, \C)$  which are of  the form $e^{\mathbbm{i}\omega t}\phi(x)$ with $\omega >0$ and $\phi(x)>0$.
Here we assume:
\begin{itemize} \item[(H3)]
there is an open interval $\mathcal{O}\subseteq (0, \infty)$
and a $C^\infty $-family $\mathcal{O}\ni\omega\mapsto \phi_{\omega}\in
\mathop{\cap}\sb{n\in\N}\Sigma_{n}(\R^3, \C),$
with $ \Sigma_{n}(\R^3, \C) $ defined in \eqref{eq:sigma}, such that $\phi_{\omega}$ is a positive radial solution of \begin{equation}
 \label{eq:B}
 -\Delta u + \omega u+\beta(|u|^2)u=0\qquad\text{for $x\in \R^3$;}
\end{equation}

\item[(H4)] we have $\frac d {d\omega } \| \phi_{\omega }\|^2_{L^2}>0$
for $\omega\in\mathcal{O}$;

\item[(H5)] for $L_+:=-\Delta +\omega +\beta (\phi_\omega^2)+2\beta
 '(\phi_\omega^2) \phi_\omega^2$ with the domain $H^2 (\R^3,\C)$,
  $L_+$ has one negative eigenvalue
 and
 $\ker L_+ =\text{Span}\{ \partial_{x_j}\phi_{\omega}:\; j=1,2,3\}$.
\end{itemize}
The above hypotheses guarantee that the ground states are orbitally
stable solutions of the scalar NLS \eqref{Eq:scalNLS}, see \cite{GSS2,W2}.  In  \cite{MR2805462,Cu3} it has been proved that, under some additional hypotheses, these ground states are asymptotically stable,
in a sense that will be clarified later.  This paper shows that some solitary waves of
\eqref{Eq:NLS}   are    asymptotically stable. To state the result,  we denote by $\mathbf{C}:\;\C^n\to\C^n$   the operator if complex conjugation in $\C^n$. We consider the group
$\mathbf{SU}(2)$, given  by
\begin{equation} \label{eq:su2}
 \mathbf{SU}(2) =\left \{ \begin{bmatrix}
a&b\\- \mathbf{C} b& \mathbf{C} a
\end{bmatrix} :\; (a,b) \in \C^2 \text{ such that } \abs{a}^2+\abs{b}^2=1 \right \}.
\end{equation}
We consider the group
\begin{equation}\label{group G}
 \mathbf{G}= \R^3 \times \mathbb{T}\times \mathbf{SU}(2).
\end{equation}
There is a natural representation of $\mathbf{G}$  on $H^1
 (\R^3, \C^2)$, with $\vartheta \in  \mathbb{T}$  acting  on $u_0$   like $e^{\mathbbm{i}\vartheta} u_0$, $x_0\in \R^3$ acting like a translation operator, and with an element
 of $ \mathbf{SU}(2)$ acting on $u_0$  by transforming it into
$
       (a +b \sigma_2\mathbf{C})u_0.$
  System \eqref{Eq:NLS} admits solitary waves of the form
  \begin{equation}\label{eq:psi}
\psi _{\omega ,v} (t) =
 e^{{\mathbbm{i}} t \left (\omega + \frac{v^2}{4} \right) } e^{\frac{\mathbbm{i}} 2v\cdot (x -tv) } \phi_\omega (x-tv)    \overrightarrow{e}_1 \text{ ,  with the column vector $\overrightarrow{e}_1= {^{t}(1,0)} $.}
 \end{equation}
We will show later that, along with mass, which we will denote by $\Pi _4$,    linear momenta, which we will denote by $\left .\Pi _i \right | _{i=1}^{3}$, and energy,  system \eqref{Eq:NLS}
 admits 3 further invariants related to $ \mathbf{SU}(2)$ which we will denote by $\left .\Pi _i \right | _{i=5}^{7}$. By $\Pi$ we will denote the vector $\left .\Pi _i \right | _{i=1}^{7}$. We will see later, that acting with
 $\mathbf{G}$ on $\psi _{\omega ,v}$ we can generalize the solitary waves. We will have solitary waves $\Phi _p$ characterized by  $\Pi (\Phi _p)=p$. We will prove the following theorem.

 \begin{theorem}\label{theorem-1.1}
Assume (H1)--(H5) stated above, (H6)--(H8) stated in Section~\ref{sec:speccoo},
and (H9) stated in Sect.~\ref{sec:system}. Pick $\omega ^1\in \mathcal{O} $.  Then there exist $\epsilon_0 = \epsilon_0 (\omega ^1)>0$
and $C= C (\omega ^1)>0$
such that
if $u $ solves \eqref{Eq:scalNLS}
with $u\vert_{t=0}=u_0$ and if we set
\begin{equation}\label{def:epsilon} \epsilon:=\inf_{ g\in \mathbf{G} }\|u_0- T(g)\psi _{\omega ^1,0}(0)
\|_{H^1(\R^3, \C^2)}<\epsilon_0,
\end{equation}
% and moreover \begin{eqnarray}\label{boxed} \Pi (u_0) =\Pi (\psi _{\omega ^0,0}  ),\end{eqnarray}
 then
there exist a solitary wave $\psi _{\omega ^+,v^+}$,  a function $ g\in C^1(\R_+, \mathbf{G})$ and an element $h_+ \in H^1(\R^3, \C^2) $ with
$\| h_+ \|_{H^1(\R^3, \C^2) }+|\omega _+ -\omega ^1 |+|v^+| \le C \epsilon $,
such that
\begin{equation}\label{eq:scattering} \lim_{t\to +\infty} \|u(t)- T(g(t)) \psi _{\omega ^+,v^+} (t) -
 e^{ - {\mathbbm{i}} \sigma_3 \Delta t} h_+ \|_{H^1(\R^3, \C^2) }=0. \end{equation}
\end{theorem}

 This result is a  transposition to a system, of the result proved for scalar equations in   \cite{MR2805462,Cu3,CM1}, see also \cite{MR3053771}.
We are not aware of previous similar results   for systems of  PDE's. Papers such as \cite{boussaid,PelinovskyStefanov,MR2924465}
treat special cases of the Dirac Equation, but don't touch the novel issues treated here.
For the orbital stability of systems of NLS we refer to Grillakis et  al. \cite{GSS2}, see also
\cite{Bhattarai} and therein.

The proof of Theorem~\ref{theorem-1.1} goes along the lines of the proof for the scalar NLS.
If we look at the analogous classical problem of the asymptotic stability of the equilibrium $0$ for a system
$\dot r =Ar +g (r)$, where $g (r)=o(r)$ at $r=0$ and with a matrix $A$, of key importance is the location
of the spectrum $\sigma (A)$. Stability requires that if $\varsigma \in \sigma (A)$ then $\Re \varsigma\le 0$.
Isolated eigenvalues on the imaginary axis correspond to central directions whose contribution to stability
or instability can be ascertained only analyzing the nonlinear system, and not just the linearization $\dot r =Ar $.
This classical framework is also used for Theorem \ref{theorem-1.1}. First of all,
an appropriate expansion of $u$ at the ground states
(see Lemma~\ref{lem:modulation} below) gives us the variable $r$.
The analogue of $A$ is given by \eqref{eq:lincom}. In our case the spectrum is all contained in the imaginary axis,
but the continuous spectrum  plays the same role of the stable spectrum of $A$, thanks to dispersion and along the lines
described in pp. 36--37 of Strauss's introduction to nonlinear wave equations \cite{strauss}. The discrete spectrum of \eqref{eq:lincom} plays the role of central directions.
The nonlinear mechanism acting on the corresponding discrete modes and responsible for the stabilization
indicated in \eqref{eq:scattering}  has been termed \textit{Nonlinear Fermi Golden Rule}
in \cite{sigal} and was explored initially in \cite{MR1334139,SW3}. A detailed description,
by means of some elementary examples, is in \cite[Introduction]{CM4},  see also \cite{W3}.
The same mechanisms described in \cite{CM4} and used in \cite{MR3053771,MR2843104,MR1334139,MR2805462, Cu3,CM1,SW3} and in a number
of other papers referenced therein,
 are applied here
to prove Theorem~\ref{theorem-1.1}. The novel difficulty occurs   with the choice of modulation.
Here the   the idea is to use the representation \eqref{eq:rep1}.
The rest of the paper is not very different from \cite{MR2805462,Cu0,Cu3,CM1}. In the course of the proof there
are some difficulties related to the fact that the Lie algebra of $\mathbf{G}$ is not commutative,
and correspondingly, the Poisson brackets $\{ \Pi_j, \Pi_l \}$ are not identically zero like in the earlier papers. This   is solved quite naturally by exploiting conservation laws and
considering the reduced manifold, see \cite[Ch. 6]{Olver}. Thanks to an appropriate uniformity
with respect to the conserved quantities
of the coordinate changes, we obtain the desired result.

\section{Notation and preliminaries}\label{sec:notation}

We start with  some notation. For $\varsigma\in \C^n$ we consider the Japanese Bracket $\langle \varsigma \rangle =\sqrt{1+|\varsigma |^2}$.

\noindent Given two Banach spaces $\mathbb{X}$ and $\mathbb{Y}$
 let $B(\mathbb{X},\mathbb{Y})$ be the Banach space of bounded linear transformations from $\mathbb{X}$ to $\mathbb{Y}$.

\noindent
Let
$m,\,k,\,s\in\R$.
Given a Banach space $\mathbb{E}$ and functions $\R^3\to \mathbb{E}$,
we denote by
$\Sigma_m(\R^3,\mathbb{E}) $
and $H^{ k,s}(\mathbb{R}^3,\mathbb{E})$
the Banach spaces with the norms
\begin{align}\label{eq:sigma}
&
\| u \|_{\Sigma_m}^2:= \|
\big\langle
\sqrt{-\Delta +|x|^2}\,
\big\rangle^{m} u \|_{L^2(\R^3,\mathbb{E})}^2,
\\[1ex]
&
\label{eq:weigt}
 \| f\|_{H^{k,s}(\R^3,\mathbb{E})}:=\|
\langle x \rangle^s
\big\langle \sqrt{-\Delta }
\,\big\rangle^{k} f
 \|_{L^2(\R^3,\mathbb{E})},
\end{align}
where we will use mostly $\mathbb{E} =\C^2$.
We also consider
\begin{align}\label{Schwartz}
  \text{the space of Schwartz functions }
\mathcal{S}(\R^3, \mathbb{E}) := \cap_{m\in \R }\Sigma_{m}(\R^3, \mathbb{E} ) \\
\text{ and the space of tempered distributions }   \mathcal{S}'(\R^3, \mathbb{E}) := \cup_{m\in \R }\Sigma_{m}(\R^3, \mathbb{E} ) .
\label{tempered}
\end{align}
We denote by $ {^t } {{v}}$ the transpose of $ {v}\in\C^n$,
so that the hermitian conjugate of $v\in\C^n$ is given by ${^t }({\mathbf{C}v})$.
For $ u,\,v \in \C^n$ we set
 $|v|^2= {^t(\mathbf{C}}v) v.$
We denote the hermitian form
in $L^2(\R^3,\C^2)$ by
\begin{equation} \label{eq:hermitian}\begin{aligned}& \langle u,v\rangle =\Re \int
_{\mathbb{R}^3}{^t(\mathbf{C}u (x))}\,{v}(x)\,dx,
\qquad
u,\,v\in L^2(\R^3,\C^2),
%% \text{ for $u,v:\mathbb{R}^3\to \mathbb{C}^2$ }
 \end{aligned}\end{equation}
and we consider the symplectic form  \begin{equation} \label{eq:Omega} \Omega (X,Y) := \langle {\mathbbm{i}} \sigma_3 X, Y \rangle,
\qquad
X,\,Y\in L^2 (\R^3, \C^2).
\end{equation}

\begin{definition}
Given a differentiable function $F$,
its  Hamiltonian vector field  with respect to a {\it strong}
symplectic form $\Omega $ is the field $X_F $ such that $ \Omega (X_F,Y)= dF (Y)$ for any
tangent vector $Y $, with $dF$ the Fr\'echet derivative.
  For $F,G $ differentiable functions
their Poisson bracket is
$
 \{ F,G \} := dF (X_G)
$
if $G$ is scalar valued and $F$ is either scalar or has values in a
Banach space $\mathbb{E}$.
\end{definition}
Notice that since $X\to  \langle {\mathbbm{i}} \sigma_3 X, \  \rangle$  defines an isomorphism of $L^2(\R^3,\C^2)$, or of $H^1(\R^3,\C^2)$, into itself,
our symplectic form \eqref{eq:Omega} is strong.
For $u\in H^1(\R^3, \C^2) $ we have the following functionals
(the linear momenta and   mass)
which are   conserved in time by \eqref{Eq:NLS}:
\begin{align}
& {\Pi }_a(u)= 2^{-1} \langle \Diamond_a u, {u} \rangle \,, \quad \Diamond_a:=-\mathbbm{i}\sigma_3 \partial_{x_a}\text{ for $a = 1,2, 3$;}\label{eq:charge1}
\\
 & \Pi_4(u)= 2^{-1}\langle \Diamond_4 u, u \rangle\,, \quad \Diamond_4:=\uno (=\text{identity operator}) ;\label{eq:charge}
\end{align}
see \cite[(2.6) and p. 343]{GSS2} for \eqref{eq:charge1}.
We also consider the following functionals $\Pi_j$, $j=5,6,7$:
\begin{equation}\label{eq:charge2} \Pi_j (u):= 2^{-1} \langle \Diamond_j u,u\rangle
\text{ with }
 \Diamond_j:= \left\{\begin{matrix}
\sigma_3\sigma_2
\mathbf{C},&
j=5,
\\
\mathbbm{i} \sigma_3\sigma_2
\mathbf{C},
&
j=6,
\\
\sigma_3,&j=7.
\end{matrix}\right. \end{equation}
 The energy is defined as follows: for $B(0)=0$ and $B'=\beta$ we write
\begin{align}
\label{eq:energyfunctional}
&
 E(u):=E_K(u)+E_P(u),
\\
&
E_K(u):= 2^{-1} \langle
 -\Delta u,{u} \rangle,
\qquad
E_P(u):=-2^{-1}
 \int_{\R^3}B(|u|^2)\,dx.
\nonumber
\end{align}
It is a standard fact which can be proved like for the scalar equation
\eqref{Eq:scalNLS}, for the latter see \cite{MR2002047}, that (H1)--(H2) imply local well--posedness of
 \eqref{Eq:NLS}  in $H^1 (\R^3, \C^2)$.

 \noindent We denote by  $dE$   the Fr\'echet derivative of the energy  $E$, see \eqref{eq:energyfunctional}.
 We define $ \nabla {E}$ by $dEX= \langle \nabla {E}, {X} \rangle$. Notice that
 $ \nabla {E}\in C^1 (H^1(\R^3, \C^2), H^{-1}(\R^3, \C^2))$, that $\nabla {E}(u) =-\Delta u +\beta(|u|^2) u$
and henceforth that \eqref{Eq:NLS} can be written as
\begin{equation}\label{eq:graden-dot}
\dot u =-{\mathbbm{i}}\sigma_3 \nabla {E}(u) = X_{E} (u),
\end{equation}
that is, as a hamiltonian system with hamiltonian $E$. Notice that $\nabla \Pi_j(u) =\Diamond_j u$ for $j=1\le j\le 7$.

\noindent Consider now the column vector $\overrightarrow{e}_1= {^{t}(1,0)} $.
By \eqref{eq:charge1} and (H4),
 $(\omega, v)\mapsto (\Pi_j (e^{ \sigma_3\frac{\mathbbm{i}} 2 v\cdot x }
 \phi_\omega \overrightarrow{e}_1))_{j=1}^{4}$ is a diffeomorphism into an open
subset of $ \R_+\times \R^3$.
We introduce
\begin{equation}\label{def p}
 p=p(\omega, v)\in \R^7 \text{ defined by } p_j(\omega, v)=
\begin{cases}
\Pi_j (e^{ \sigma_3\frac{\mathbbm{i}} 2 v\cdot x }\phi_\omega \overrightarrow{e}_1),
&\quad
1\le j\le 4;
\\
0,
&\quad
j=5,6;
\\
\Pi_j (e^{ \sigma_3\frac{\mathbbm{i}} 2 v\cdot x }\phi_\omega \overrightarrow{e}_1)= p_4(\omega, v),
&\quad
j=7. \end{cases}
\end{equation}
Notice that $ \Pi_j (e^{ \sigma_3\frac{\mathbbm{i}} 2 v\cdot x }\phi_\omega \overrightarrow{e}_1) =0$ for $j=5,6$.  We denote by $\mathcal{P}$ the subset of $ \R^7 $  defined by
\begin{align}\label{def-p}
\mathcal{P}=
\{
p(\omega,v);\,
\omega\in\mathcal{O},
\ v\in\R^3
\}.
\end{align}
For $p=p(\omega, v)\in \mathcal{P}$, we set
\begin{equation}\label{eq:defPhi}
\Phi
_{p}(x):=e^{\frac{\mathbbm{i}} 2v\cdot x } \phi_\omega (x) \overrightarrow{e}_1.
\end{equation}
Obviously $\Phi
_{p(\omega, v)} =\psi _{\omega ,v} (0) $, see \eqref{eq:psi}.
  We will set $\Phi
_{p_1}= \psi _{\omega ^1,0} (0)$ for the function in Theorem \ref{theorem-1.1}.
We have $\Pi _j(\Phi
_{p_1})= 0$ for $j=1,2,3,5,6$. It   is not restrictive
to pick the initial datum s.t.
\begin{equation}\label{eq:normzion}
    \Pi _j(u_0)= 0 \text{  for $j=1,2,3,5,6$.}
\end{equation}
Indeed, by continuity, $  \Pi _j$ for $j=1,2,3,5,6$ take values close to 0 in a neigborhood
of $\Phi
_{p_1} $. By boosts and Lemma \ref{lem:lambdas1}, one can act   on $u_0$ changing it into another nearby initial datum  which satisfies  \eqref{eq:normzion}: we skip the elementary details.  We introduce   \begin{equation}\label{eq:lambda0}
 \lambda (p) =(\lambda_1(p),..., \lambda_7(p)) \in \R^7 \text{ defined by } \lambda_j(p):=
\begin{cases}
-v_j,&
1\le j\le 3;
\\ -\omega - \frac{v^2}4,
\quad
&
j=4;
\\ 0,&
j\ge 5.
\end{cases}
\end{equation}They are Lagrange multipliers, and
an elementary computation shows that
 \begin{equation}\label{eq:solw0}
   e^{- {\mathbbm{i}} \sigma_3 t \lambda (p) \cdot \Diamond}\Phi_p = \psi _{\omega ,v} (t)
 \end{equation}
 and that $\Phi_{ p }$ is a constrained critical value for the energy  satisfying
 \begin{equation}\label{eq:eqphi}
 \nabla E (\Phi_{ p })- \sum_{j=1,...,7}\lambda_j (p) \Diamond_j\Phi_{ p } =0.
 \end{equation}
We consider the representation $T:\mathbf{G}\to B (H^1
 (\R^3, \C^2), H^1
 (\R^3, \C^2))$ defined by
\begin{align}\label{eq:rep1}&
 T(g) u_0 := e^{{\mathbbm{i}} \sigma_3  \tau \cdot  \Diamond } (a +b \sigma_2\mathbf{C})u_0 \text{ for } g= \left (\tau , \begin{bmatrix}
a&b\\-\mathbf{C} b&\mathbf{ C} a
\end{bmatrix} \right) \text{ where}\\& \tau =(\tau _1,\tau_2,\tau _3,\tau _4)\in  \R^3 \times \mathbb{T} \text{ and }\tau \cdot  \Diamond := \sum _{j=1,...,4} \tau _j\Diamond  _j .\nonumber
\end{align}
An elementary but very important fact to us is the following lemma.

\begin{lemma} \label{lem:cons0}
 We have the following facts.
\begin{itemize}
\item[(1)] The action of $\mathbf{G}$ given by \eqref{eq:rep1} preserves the symplectic form $\Omega$ defined in \eqref{eq:Omega}.

\item[(2)] The action \eqref{eq:rep1} preserves the invariants $\Pi_j$ for $1\le j\le 4$ and $E$.

\item[(3)] The functionals $\Pi_j$, $1\le j\le 7$,  and   $E $
 are conserved by the flow of \eqref{Eq:NLS} in $H^1
 (\R^3, \C^2)$.
\end{itemize}

\end{lemma}\proof (1) follows from the
commutation $[{\mathbbm{i}} \sigma_3, a+b\sigma_2\mathbf{C} ]=0$. (2) is a consequence of

\begin{equation*}
\begin{aligned}
 &\abs{ (a+b\sigma_2\mathbf{C})u}^2
=
\Re
{^t(\mathbf{C}u)}
((\mathbf{C} a)+\mathbf{C}\sigma_2(\mathbf{C} b))
(a+b\sigma_2\mathbf{C})u
\\&=
(\abs{a}^2+\abs{b}^2)
\abs{u}^2
+
\Re
{^t(\mathbf{C}u)}
((\mathbf{C} a) b\sigma_2\mathbf{C}+\mathbf{C}\sigma_2 a \mathbf{C} b)
u
=\abs{u}^2.
\end{aligned}
\end{equation*}
The fact that the functionals $\Pi_j$, $1\le j\le 4$, and the energy $E $
are preserved by the flow of \eqref{Eq:NLS} is standard.
To deal with the cases $j=5, 6,7$,
we first recall that the Lie algebra
of $\mathbf{SU}(2)$
can be written as
%%AC $\mathbf{su}(2) =\oplus_{i=1}^{3} \R {\mathbbm{i}} \sigma_i$.
$\mathbf{su}(2) =
\mathop{\mathrm{Span}}
\left(\mathbbm{i}\sigma_i,\,1\le i\le 3\right)$.
We have
\begin{align}\label{diffT}
\frac{d}{dt} \left.T(e^{-{\mathbbm{i}}t \sigma_i})\right |_{t=0}
=
\begin{cases}
\frac{d}{dt} \left. (\cos (t)-{\mathbbm{i}} \sin (t) \sigma_2 \mathbf{C})   \right |_{t=0} =-{\mathbbm{i}} \sigma_2 \mathbf{C}  ,
&
i=1;
\\
\frac{d}{dt} \left. (\cos (t)+ \sin (t) \sigma_2 \mathbf{C})   \right |_{t=0} =  \sigma_2 \mathbf{C}  ,
&
i=2;
\\
\frac{d}{dt} \left. e^{-{\mathbbm{i}}t}   \right |_{t=0} = -{\mathbbm{i}}  ,
&
i=3.
\end{cases}
\end{align}
Like in  \cite[line 5 p.313]{GSS2},
\begin{align*}& \frac{d}{dt}\Pi  _{4+i}(u) = \left \langle \Diamond _{4+i}u, -{\mathbbm{i}}\sigma _3\nabla E(u) \right \rangle = \left \langle {\mathbbm{i}} \sigma _3\Diamond _{4+i}u,  \nabla E(u) \right \rangle \\& = \left .\frac{d}{ds} \left \langle T(e^{ {\mathbbm{i}}s \sigma_i})u,  \nabla E(u) \right \rangle \right |_{s=0} = \left .\frac{d}{ds} E(T(e^{ {\mathbbm{i}}s \sigma_i})u)  \right |_{s=0} =0,
\end{align*}
where the 1st equality holds for sufficiently regular solutions, while the last one follows from (2). By a density argument and well posedness of \eqref{Eq:NLS}, we obtain claim (3).
\qed

\begin{lemma}\label{lem:Tmanifold0}
The following 10 vectors are linearly independent over $\R$:
\begin{equation}\label{eq:Tmanifold}
 \partial_{p_1}\Phi
_{p},\ \partial_{p_2}\Phi
_{p},\ \partial_{p_3}\Phi
_{p}, \ \partial_{p_4}\Phi
_{p}, \ \partial_{x_1}\Phi
_{p},\ \partial_{x_2}\Phi
_{p},\ \partial_{x_3}\Phi
_{p},\ {\mathbbm{i}} \sigma_2 \mathbf{C} \Phi
_{p}, \ \sigma_2 \mathbf{C} \Phi
_{p}, \ {\mathbbm{i}} \Phi
_{p}.
\end{equation}
\end{lemma}The  proof is elementary. \qed

We consider now the ``solitary manifold''
\begin{equation} \label{eq:manifold0} \begin{aligned}
\mathcal{M} :&= \left\{ e^{{\mathbbm{i}} \sigma_3 \tau \cdot  \Diamond} (a +b \sigma_2 \mathbf{\mathbf{C}}) \Phi
_{p}(x) :\; \tau \in \R^3\times \mathbb{T}, \ \begin{bmatrix}
a&b\\- \mathbf{C} b& \mathbf{C } a
\end{bmatrix} \in \mathbf{SU}(2),\ p
\in \mathcal{P} \right\}\,. \end{aligned}
\end{equation}
The vectors in \eqref{eq:Tmanifold} are obtained
computing the partial derivatives in $(0,p,0)$ of
\begin{align}
&
\text{ the function in }  C^\infty (\mathbb{D}_{\C}(0,\varepsilon_0)
\times \mathcal{P}\times \mathbb{T} \times \R^3, \Sigma_k (\R^3,\C^2))\text{ given by }
\nonumber
\\
&
(b,p, \tau) \mapsto e^{{\mathbbm{i}} \sigma_3   \tau \cdot   \Diamond } \s (b) \Phi_{p}, \text{ where }  \s (b):= \sqrt{1-|b|^2} +b\sigma_2 \mathbf{C}
 .
\label{eq:s sym}
\end{align}
Then
Lemma~\ref{lem:Tmanifold0}   implies that for any $k>0$ there is $\varepsilon_0>$
s.t.  \eqref{eq:s sym}
is an embedding and   $\mathcal{M} $ is a manifold.
The $\R$--vector space generated by  vectors  in Lemma \ref{lem:Tmanifold0}   is the tangent space $T_{\Phi
_{p}}\mathcal{M}$.

\noindent Consider the \emph{ linearized}  operator ${\mathcal H}_{p}:= -{\mathbbm{i}} \sigma_3(\nabla^2 E(\Phi_{p})- \lambda (p) \cdot \Diamond) $.
By $\lambda (p(\omega,0)) \cdot \Diamond =-\omega $ we have
\begin{equation} \label{eq:lincom} \begin{aligned} & {\mathcal H}_{p(\omega,0)} \begin{pmatrix} u_1
 \\
 u_2
 \end{pmatrix} = - \begin{pmatrix}{\mathbbm{i}} \mathfrak{L}^{(1)}_{\omega
 }u_1 \\ - {\mathbbm{i}} \mathfrak{L}^{(2)}_{\omega
 }u_2 \end{pmatrix},
\quad \text{ where }\\& \mathfrak{L}^{(1)}_{\omega
 }u_1= - \Delta u_1 + \beta (\phi_{\omega}^2) u_1+ 2 \beta ' (\phi_{\omega}^2)\Re (u_1) + \omega u_1,
 \\& \mathfrak{L}^{(2)}_{\omega
 }u_2 =-\Delta u_2 + \beta (\phi_{\omega}^2) u_2+ \omega u_2.
 \end{aligned}
\end{equation}
It is well known that ${\mathcal H}_{p}$ is $\R$-linear but not
 $\C$-linear, see \cite{MR1199635,MR1835384}. For this reason we interpret $H^1(\R^3, \C^2)$ as a vector space over $\R$. Later,
in  Section~\ref{sec:speccoo}, we perform a complexification.  Recall the  generalized kernel $N_g ({\mathcal H}_{ p }) :=\cup_{j=1}^{\infty} \ker ({\mathcal H}_{ p })^{j} $. The following lemma is very important.

\begin{lemma}\label{lem:kernel}
We have  $N_g ({\mathcal H}_{ p(\omega,0)}) =T_{\Phi
_{ p(\omega,0)}}\mathcal{M} $.
\end{lemma}
\proof First of all $\mathfrak{L}^{(i)}_{\omega
 }$ for $i=1,2$ are decoupled,
so that it is enough to consider them separately. We have the following, which is a  well-known   fact about ground states, see for example \cite[Sect.XIII.12]{RS4}:
\begin{equation*}
 \ker ({\mathbbm{i}} \mathfrak{L}^{(2)}_{\omega
 })=N_g({\mathbbm{i}} \mathfrak{L}^{(2)}_{\omega
 }) =\text{Span} \{ {\mathbbm{i}} \phi_{\omega}, \phi_{\omega} \}.
\end{equation*}
 The following   well-known consequence of (H4)--(H5),  derived in  \cite{W2},   completes the proof
\begin{equation*} \begin{aligned} & \ker ({\mathbbm{i}} \mathfrak{L}^{(1)}_{\omega
 }) =\text{Span} \{ {\mathbbm{i}} \phi_{\omega}, \left. \partial_{x_a}\phi_{\omega} \right |_{a=1}^{3} \},
\\& N_g(\mathfrak{L}^{(1)}_{\omega
 }) = \ker ({\mathbbm{i}} \mathfrak{L}^{(1)}_{\omega
 })^2 = ({\mathbbm{i}} \ker \mathfrak{L}^{(1)}_{\omega
 })\oplus \text{Span} \{   \partial_{p_j} e^{\frac{\mathbbm{i}}2v\cdot x }\phi_{\omega}  |_{j=1}^{4} \}.
 \end{aligned} \end{equation*}
\qed

System  \eqref{Eq:NLS} is  an interesting example for the stability theory    in the classical paper by Grillakis {\it et al.} \cite{GSS2}   because all the examples of systems of NLS's in
Sect. 9 in  \cite {GSS2}     for  $x\in \R^3$ and $u(t,x)\in \R ^4$
have      4--dimensional centralizers, while    for \eqref{Eq:NLS} dimension is 6,  see the following  two remarks.

\begin{remark}\label{rem:repr1} From the identification $\C^2=\R^4$ there is a natural inclusion $ \mathbf{SU}(2)\subseteq \mathrm{SO}(4)$. By  the identification implicit in \eqref{eq:su2} of $ \mathbf{a} \in \mathbf{SU}(2)$ and an element in the unit sphere $\widetilde{ {{\mathbf{a}}}}\in S^3 \subset \R^4$,  the action of $ {\mathbf{a}}\in \mathbf{SU}(2) $ on $ {\mathbf{v}}\in \R^4$ is nothing else
but the product of quaternions, ${\mathbf{v}}\widetilde{ {{\mathbf{a}}}}$. Similarly, by elementary computations it is possible to see that $ (a +b \sigma_2\mathbf{C}){\mathbf{v}} =\widehat{ {{\mathbf{a}}}}{\mathbf{v}} $
(on the r.h.s. multiplication of two quaternions)
for all $\mathbf{v}\in \R^4$ and for an appropriate
$\widehat{ {{\mathbf{a}}}}\in S^3$.
In the framework of \cite{GSS2} when applied to \eqref{Eq:NLS}, a key role is played by the centralizer of the group
$\{e^{ \tau_4{\mathbbm{i}} \sigma_3};\,\tau_4\in\R\}$  inside  $\R^3\times \mathrm{SO}(4)$.
Using \cite[p. 111]{wolf},
   it can be shown  that  $\mathbf{G}$, acting as in \eqref{eq:rep1}, is a connected component of
this centralizer.
 %A simple consequence of Lemma \ref{lem:kernel}, see also Remark \ref{assumption 3} below, is that \cite[Assumption 3]{GSS2}  is satisfied for this centralizer.
\end{remark}

\begin{remark}\label{assumption 3} The key hypothesis in  \cite{GSS2} is Assumption 3 on p. 314, stating $Z=\ker ({\mathcal H}_{ p(\omega,0)})$  for
\begin{equation*}
 Z:=\left \{  \partial _t \left.\widetilde{T}(e^{t\varpi}) \Phi_{ p(\omega,0)}\right |_{t=0} :\; \varpi\in
 \R^3\times so(4) \text{ commutes in $\R^3\times so(4)$ with ${\mathbbm{i}} \sigma_3$} \right \} \ ,
\end{equation*}
   where for $\varpi\in \R^3$ we have $\widetilde{T}(e^{t\varpi})=
{T}(e^{t\varpi}) $ and for $\varpi \in so(4)$  we set $\widetilde{T}(e^{t\varpi}) w = e^{t\varpi}  w $
for any $w\in \R^4$, with the usual product row column $ SO(4)\times \R ^4\to \R ^4$.

\noindent  Always
 $Z\subseteq\ker ({\mathcal H}_{ p(\omega,0)})$, see \cite[Lemma 2.2]{GSS2}.
 Lemma~\ref{lem:kernel} yields the equality.
   Assumption 1, i.e.   local well posedness, is true   and
     Assumption 2, about bound states, is true under our hypothesis (H3). Other hypotheses needed in  \cite{GSS2}, such as   that the centralizer, or at least its connected component containing the unit element   in
 $\R^3\times SO(4)$, acts by symplectomorphisms which leave
 the energy invariant,  follow from Lemma \ref{lem:cons0}.
  So by  \cite{GSS2} the   bound states \eqref{eq:solw0} are $\mathbf{G}$--orbitally stable.
 \end{remark}

\section{Modulation}\label{sec:mod}

 The manifold $\mathcal{M}$ introduced in \eqref{eq:manifold0} is a  symplectic submanifold of $L^2(\R^3, \C^2)$.
This follows from
\begin{equation*}\begin{aligned} &
 \Omega ({\mathbbm{i}} \sigma_2 \mathbf{C} \Phi
_{p}, \sigma_2 \mathbf{C} \Phi
_{p}) =p_4\ ,  \quad  \Omega (\partial_{p_4}\Phi
_{p}, {\mathbbm{i}} \Phi
_{p}) = 2^{-1} \partial_{p_4} \langle {\mathbbm{i}} \sigma_3 \Phi_{p },{\mathbbm{i}} \Phi_{p } \rangle = \partial_{p_4}p_4 =1,
\\&
 \Omega (\partial_{p_a}\Phi
_{p}, \partial_{x_a}\Phi
_{p}) = 2^{-1} \partial_{p_a}\langle \Phi_{p }, \Diamond_{a} \Phi_{p }\rangle = \partial_{p_a}p_a =1 \text{ for $a=1,2,3$},
\end{aligned}
\end{equation*}
and from   symplectic  orthogonality of all    other pairs of vectors in  \eqref{eq:Tmanifold}. We obtain a  bilinear form
\begin{equation*}
 \Omega :\; \mathcal{S}(\R^3, \C^2) \times \mathcal{S}'(\R^3, \C^2) \to \R .
\end{equation*}
Since  $ T_{\Phi
_{p} } \mathcal{M} \subseteq \mathcal{S} (\R^3, \C^2) $,
  we can define the subspace $ T_{\Phi
_{p} }^{\perp_\Omega}\mathcal{M} \subseteq \mathcal{S}'(\R^3, \C^2) $.
$\Omega $ also defines a pairing
\begin{equation*}
 \Omega :\; \Sigma_n (\R^3, \C^2) \times \Sigma_{-n}(\R^3, \C^2) \to \R.
\end{equation*}
This yields the decomposition
\begin{equation}\label{eq:dirsumn} \begin{aligned} &
\Sigma_{-n}(\R^3, \C^2) = T_{ \Phi_{p} }\mathcal{M} \oplus
(T_{\Phi_{p} }^{\perp_\Omega}\mathcal{M} \cap \Sigma_{-n}(\R^3, \C^2)).
\end{aligned}
\end{equation}
We denote by $\widehat{P}_{ p} $ and $ {P}_{{ p} }$
the projections onto the first and second term of the direct sum, respectively:
\begin{align}
\label{def-widehat-p}
&&
\widehat{P}_{ p}:\;
\Sigma_{-n}(\R^3, \C^2)
\to T_{ \Phi_{p} }\mathcal{M},
\\
&&
{P}_{{ p} }:\;
\Sigma_{-n}(\R^3, \C^2)
\to
T_{\Phi_{p} }^{\perp_\Omega}\mathcal{M} \cap \Sigma_{-n}(\R^3, \C^2). \nonumber
\end{align}
A special case of \eqref{eq:dirsumn} is
\begin{equation}\label{eq:dirsum} \begin{aligned} &
 L^2(\R^3, \C^2) = T_{ \Phi
_{p} }\mathcal{M} \oplus (T_{\Phi
_{p} }^{\perp_\Omega}\mathcal{M} \cap L^2(\R^3, \C^2)).
\end{aligned}
\end{equation}
It is easy to see  that the map
$p\mapsto \widehat{{P}}_{{p}}$
is in  $C^\infty(\mathcal{P},B(\Sigma_{-n}(\R^3,\C^2),\Sigma_{n}(\R^3,\C^2)))$ for any $n\in \Z$.
% Indeed, if we denote by $( e_{i})\left . \right |  _{i=1}^{10}$ the 10 vector fields  in \eqref{eq:Tmanifold}, which belong to $C^\infty(\mathcal{P}, \Sigma_{-n}(\R^3,\C^2))$ for any $n$, and if we denote by $A=\{  a _{ij}\} $ the inverse matrix to $ \{  \Omega ( e_i, e_j)  \} $, then\begin{equation*} \widehat{{P}}_{{p}} =\sum _{i,j=1,...,10}e_i \  a _{ij}\ \Omega ( e_j, \cdot ) . \end{equation*}
The following    about the $\s (b)$ in  \eqref{eq:s sym}    is consequence of elementary computations:
\begin{align}& (\s (b)) ^{-1}= (\s (b)) ^{*}= \s (-b) \  ;\label {eq:s sym1}  \\& \s (b) \sigma _j=  \sigma _j\s (-b)  \text{ for all }j=1,2,3 \, ; \nonumber \\& \mathbf{C}\s (b)=\s \left ( -\mathbf{C}b \right )  \, , \quad \s (b)  {\mathbbm{i}} = {\mathbbm{i}}\s (-b)  . \nonumber
\end{align}

%\begin{remark}\label{rem:invariants} Notice that, in the notation of  Lemma \ref{eq:lambdas11}, we have $\Pi_{5}(u)\sim b_R$ and $\Pi_{6}(u)\sim b_I$   for $\psi = \Phi _p $. This follows from    $\Pi _5 (\Phi _p )=\Pi _6 (\Phi _p  )=0$ and  $\Pi _7 (\Phi _p  )=p_4  $.\end{remark}

\begin{lemma}[Modulation]
 \label{lem:modulation} Fix $n_1  \in\N_{0}:=\mathbb{N}\cup \{  0 \} $
	and
   $ p^1 \in \mathcal{P} $. Then $\exists$ an open neighborhood $ \mathcal{U}_{-n_1 } $
of $ \Phi_{ p^1 }$
in $ \Sigma_{-n_1  }(\R^3, \C^2) $
 and functions $p \in C^\infty (\mathcal{U}_{-n_1   }, \mathcal{P})$,
	 $\tau \in C^\infty ({\mathcal{U}}_{-n_1   }, \R^3\times \mathbb{T})$ and $ b \in C^\infty ({\mathcal{U}}_{-n_1   }, \C)$
such that
$p (\Phi_{ p^1 })=p^1 $, $\tau (\Phi_{ p^1 }) =0$,
 $b (\Phi_{ p^1 }) = 0$ and $\vartheta (\Phi_{ p^1 }) = 0$
so that
for any $u\in \U_{-n_1 } $,
\begin{align}\label{eq:ansatz} &
 u = e^{- {\mathbbm{i}} \sigma_3 \tau (u) \cdot \Diamond }   \s (b(u))      (\Phi_{ p (u) } +R(u)),   \text{    with
$R(u)\in T^{\bot_{\Omega}}_{\Phi_{ p (u) }}\mathcal{M}\cap\Sigma_{-n_1 }(\R^3,\C^2)$.}
\end{align}

\end{lemma}
\proof  The proof is standard. For $\mathbf{v}_{\iota}(p)$, $1\le\iota\le 10$  varying among
the $10$ vectors in \eqref{eq:Tmanifold}, set
\begin{equation*}
 F_{\iota}(u,p,\tau,b) := \Omega (e^{ {\mathbbm{i}} \sigma_3 \tau \cdot \Diamond  } \s (-b)u-\Phi_{ p },\mathbf{v}_{\iota}(p)).
\end{equation*}
 Next, setting $\overrightarrow{F}=(F_1,..., F_{10})$, we compute
\begin{align}
& \left.  \overrightarrow{ F} (u,p,\tau,b)  \right |_{u= e^{ -{\mathbbm{i}} \sigma_3\tau \cdot \Diamond } \Phi_{ p }, \ b=0 } =0 \text{ and  the Jacobian matrix is}
\nonumber
\\& \label{eq:Jac}
 \left. \frac{\partial \overrightarrow{F} (u,p,\tau,b)}{\partial (p,\tau,b)} \right |_{u= e^{ -{\mathbbm{i}} \sigma_3\tau \cdot \Diamond } \Phi_{ p }, \ b=0 } = \left [ \varepsilon_{ij} \Omega (\mathbf{v}_i(p),\mathbf{v}_j(p)) \right ]_{i,j},
\qquad
1\le i,j\le 10,
\end{align}
where the numbers $\varepsilon_{ij} $ belong to $\{ 1, -1\}$. Since for each $\mathbf{v}_i(p)$
there is exactly one $\mathbf{v}_j(p)$ such that $\Omega (\mathbf{v}_i(p),\mathbf{v}_j(p))\neq 0$, it follows that all the columns in \eqref{eq:Jac} are linearly independent. We can therefore
apply the implicit function theorem which yields the statement.
\qed

It can be proved, see \cite[Lemma 2.3]{Cu0}, that in  a sufficiently small neighborhood $ \mathcal{V}$ of $p^1$ in $\mathcal{P}$,
  for any   any $k\ge -n_1 $
the projection
\begin{align}\label{def-p-p}
P_{p}:\; T_{\Phi
_{p^1 } }^{\perp_\Omega}\mathcal{M} \cap \Sigma_{k}(\R^3, \C^2)
\longrightarrow T_{\Phi
_{p } }^{\perp_\Omega}\mathcal{M}\cap \Sigma_{k}(\R^3, \C^2)
\end{align}
is an isomorphism.
From Lemma~\ref{lem:modulation} we have
the parametrization
\begin{equation}\label{eq:coordinate}
 \mathcal{P}\times
(\R^3 \times \mathbb{T})
\times  {D}_{\C}(0,\varepsilon_0)\times (T_{\Phi
_{p^1 } }^{\perp_\Omega}\mathcal{M}\cap H^1 (\R^3, \C^2))
\to H^1 (\R^3, \C^2)
\end{equation}
with the \emph{modulation coordinates}
\begin{equation}\label{eq:coordinate-2}
\\ (p,\tau,b, r)
\mapsto u = e^{ -{\mathbbm{i}} \sigma_3 \tau \cdot \Diamond  } \s (b )  (\Phi_{ p } +P_{ p}r).
\end{equation}
We choose $p^0\in \mathcal{P}$ so that
\begin{equation}\label{eq:coordinate-3} \text{$ \Pi _j( u_0)= p^0 _j$ for $j \in I=\{ 1,2,3,4  \}$ }
\end{equation}
(that is $   p^0 _j=0$
for $j=1,2,3 $ and $ \Pi _4( u_0)= p^0 _4$, i.e.  $u_0$ and $\Phi _{p^0}  $    have same charge).

  In terms of these coordinates, system \eqref{Eq:NLS}, which we have also written as $\dot u = X_{E} (u)$, see
 \eqref{eq:graden-dot}, can be expressed in terms of the Poisson brackets as follows,
   see \cite[Lemma 2.6]{Cu0}:
 \begin{equation} \label{eq:SystE1} \begin{aligned} &
\dot p = \{ p, E \} \,, \quad \dot \tau = \{ \tau, E \} \,, \quad \dot b = \{ b, E \} \,, \quad \dot r= \{ r, E
\}. \end{aligned}
\end{equation}
 By the intrinsic definition of partial derivative on manifolds, see  \cite[p.25]{docarmo}, we have the following vector fields (recall $b_R=\Re (b)$ and $b_I=\Im (b)$):
\begin{equation}\label{eq:vectorfields} \begin{aligned}
&
  {\partial _{{\tau_j}}} = -{\mathbbm{i}} \sigma_3 \Diamond_ju
 \text{ for $1\le j\le 4$, }
\\
&
   {\partial _{{p_k}}}= e^{ -{\mathbbm{i}} \sigma_3\tau \cdot \Diamond  }\s (b) (\partial_{p_k}\Phi_{ p } + \partial_{p_k}P_{ p }r)\text{ for $1\le k\le 4$, } \\&   {\partial _{{b_A}}} = e^{ -{\mathbbm{i}} \sigma_3\tau \cdot \Diamond  }\partial_{b_A}\s (b) (\Phi_{p } + P_{p }r) \text{ for $A=R,I$},\end{aligned}
\end{equation}
which are obtained by differentiating by the various coordinates
the r.h.s. of the equality in \eqref{eq:coordinate-2}.

 \noindent By \eqref{eq:vectorfields} we have the elementary and crucial fact that
$ X_{\Pi_j } (u) ={\mathbbm{i}} \sigma_3\nabla \Pi_j (u) = {\mathbbm{i}} \sigma_3 \Diamond_ju$ for $1\le j\le 7$ which corresponds to formulas (2.5)--(2.6) in \cite{GSS2}. In particular we have
 \begin{equation}\label{eq:hampi}
\begin{aligned} & X_{\Pi_j } (u) = {\partial _{\tau_j}} \text{ for } 1\le j\le 4 \end{aligned}
\nonumber \end{equation} which immediately implies \begin{equation*} \label{eq:Ham.VecField1} \{ \Pi_j,\tau_k \} =-\delta_{jk}, \quad \{ \Pi_j,b_A \} = 0, \quad \{ \Pi_j,p_k \} = 0, \quad \{ r, \Pi_j \} = 0 \text{ for } 1\le j\le 4.
\end{equation*}
A natural step, which helps to reduce the number of equations in \eqref{eq:SystE1}
and corresponds to an application of Noether's Theorem to Hamiltonian systems, see \cite[Theorem 6.35, p.402]{Olver},
is to substitute each function $\left. p_j \right |_{j=1}^{4}$
in the coordinate system $(p,\tau,b, r)$ with the functions $\left.\Pi_j\right |_{j=1}^{4}$ and move to coordinates $(\left.\Pi_j\right |_{j=1}^{4},\tau,b, r)$.
Indeed, as in \cite[formula (34)]{Cu0}, we have \begin{equation}\label{eq:variables}
\begin{aligned} &
 \Pi_j = p_j+ \varrho _j + \Pi_j((P_{p} -P_{p^1 }) r) +\langle r,
 \Diamond_j (P_{p} -P_{p^1 }) r\rangle,
\text{ with } \varrho _j :=\Pi_j(r) \text{ and $1\le j\le 4$}.
\end{aligned} \end{equation}
This allows to  move from $(p,\tau,b, r)$ to $(\left.\Pi_j\right |_{j=1}^{4},\tau,b, r)$.
Furthermore,  $
  {\partial _{\tau_k }} \Pi_j(u)\equiv 0 $ for $k\le 4$ implies that the vector fields $ \left. {\partial _{\tau_k} }\right |_{k=1}^{4}$ are the same whether   defined using the coordinates $( p , \left.\tau_j\right |_{j=1}^{4},b, r)$ or the coordinates $(\left.\Pi_j\right |_{j=1}^{4}, \left.\tau_j\right |_{j=1}^{4},b, r)$. %(here we recall that the partial derivative in every single coordinate in \eqref{eq:vectorfields} depends on the other coordinates, but we have just shown that the $ \frac{\partial}{\partial \tau_k }$ for $k\le 4$ remain the same vector field).
 Hence, exploiting  the invariance $E(e^{{\mathbbm{i}} \sigma_3  \tau \cdot \Diamond }u)=E(u)$
 \begin{equation*}\label{eq:poi--1}
 \{ \Pi_j, E \} =- \{ E,\Pi_j \} =-dE X_{\Pi_j} = -dE  {\partial _{\tau_j}}=- \partial _{\tau_j}  E =0\quad \text{for $1\le j\le 4$}.
 \end{equation*}
 By these identities,
\eqref{Eq:NLS} in the new coordinates $(\left.\Pi_j\right |_{j=1}^{4},\tau,b, r)$ becomes
 \begin{align} &
\dot \Pi_j = 0 \text{ for $1\le j\le 4$}, \quad \dot \tau = \{ \tau, E \} \,, \nonumber \\&  \dot b = \{ b, E \} \,, \quad \dot r= \{ r, E
\}.   \label{eq:SystE2}\end{align}
 Notice that  we have produced a Noetherian reduction
 of coordinates because the equations of $b$ and $r$ are independent from the   ones in the 1st line.    We point out that by Lemma \ref{lem:cons0} we have    also
 \begin{equation} \label{eq:SystE22} \begin{aligned} &
\dot \Pi_j =  \{ \Pi_j ,E  \}=0 \text{ for $5\le j\le 7$} . \end{aligned}
\end{equation}

\section{Expansion of the Hamiltonian}
\label{sec:expansion}

 We introduce now  the following  new Hamiltonian,
\begin{equation} \label{eq:K} \begin{aligned} &
 {K}(u):= {E} (u)- {E} \left (\Phi_{ {p}^0}\right
) - \sum_{j=1,...,4} \lambda_j(p) \left (\Pi_j - {p}^0_j
 \right).
 \end{aligned}
\end{equation}
 For solutions $v$ of  \eqref{Eq:NLS} with initial value $v_0$
   satisfying   $\Pi _j( v_0) ={p}^0_j$ for $1\le j \le 4$,
 we have
\begin{equation*} \begin{aligned} & \{ \Pi_j, K \} = \{ \Pi_j, E \} =0
\quad
\text{ for $1\le j\le 7$ },
\\&
\{ b, K \} = \{ b, E \} \,, \qquad \{ r, K \} = \{ r, E \} \,,
\qquad
 \{ \tau_j, K \} = \{ \tau_j, E \} - \lambda_j(p) \text{ for $1\le j\le 4$ }.
\end{aligned}
\end{equation*}
Indeed, for example, since  $ \{ \Pi_j,    \Pi_k \}=0$
for $j \le 7$ and any $k\le 4$ (which follows from $[\Diamond_j,\Diamond_k]=0$ for $j \le 7$ and any $k\le 4$, cf. \eqref{eq:charge1}--\eqref{eq:charge2})
we have by Lemma
\begin{equation*} \begin{aligned} & \{ \Pi_j, K \}(v)  = \{ \Pi_j, {E} \} (v) - \sum_{j=1,...,4} \( \lambda_k   {\{ \Pi_j,    \Pi_k \}}  (v) + {( \Pi_j(v)-p_j^0  )}  \{ \Pi_j,    \lambda_k \} (v) \) =\{ \Pi_j, {E} \}(v),
\end{aligned}
\end{equation*}
where we use $\Pi_j(v)=p_j^0$.
The other Poisson brackets are computed similarly.

By $\partial_{\tau_j} {K} \equiv 0$ for $1\le j\le 4$, the evolution of the variables $\left. (\Pi_j) \right |_{j=1}^{7},b,r$ is unchanged if we consider the following new
Hamiltonian system,
\begin{equation} \label{eq:SystK} \begin{aligned} &
\dot \Pi_j = \{ \Pi_j, K \} =0 \text{ for $1\le j\le 4$, } \quad \dot \tau = \{ \tau, K \} \,, \quad \dot b = \{ b, K \} \,, \quad \dot r= \{ r, {K}
\}, \end{aligned}
\end{equation}
where $\left. (\Pi_j) \right |_{j=1}^{4},\tau,b,r$ is a system of independent coordinates, and where we consider also
\begin{equation} \label{eq:SystK11} \begin{aligned} &
\dot \Pi_j = \{ \Pi_j, K \} =0 \text{ for $5\le j\le 7$. }  \end{aligned}
\end{equation}

Key in our discussion is the expansion of $K(u)$ in terms of the coordinates $(\left. (\Pi_j) \right |_{j=1}^{4},r)$. We consider the expansion, with the canceled term is equal to 0 by \eqref{eq:eqphi} and \eqref{eq:lambda0},
\begin{equation*} \begin{aligned}
K(u)=& K(\Phi_{ p } +P_{p}r) = K (\Phi_{ p }) + \cancel{\langle {\nabla E (\Phi_{ p })- \sum_{j=1,...,4} \lambda_j(p) \nabla \Pi_j(\Phi_{ p })} , P_{p}r \rangle}
\\& + \int_0^1 (1-t)\Big\langle
\Big [ \nabla^2 E (\Phi_{ p } + t P_{p}r)P_{p}r)-\sum_{j=1,...,4}  \lambda_j(p) \nabla^2 \Pi_j(\Phi_{ p } + t P_{p}r)P_{p}r)\Big ] P_{p}r, P_{p}r
\Big\rangle\,dt.
 \end{aligned}
\end{equation*}
The last line equals (cf. \cite[(99)]{Cu0})
\begin{equation*} \begin{aligned} &
 2^{-1}\langle (-\Delta +\sum_{j=1,...,4} \lambda_j(p) \Diamond_j)P_{p}r,P_{p}r\rangle +
 \int_0^1 (1-t)\langle \nabla^2 E_P (\Phi_{ p } + t P_{p}r) P_{p}r, P_{p}r \rangle\,dt=\\& 2^{-1}\langle (-\Delta +\sum_{j=1,...,4} \lambda_j(p) \Diamond_j)P_{p}r,P_{p}r\rangle + \int_{\R^3}dx \int_{[0,1]^2}
 \frac{t^2}{2 } (\partial_t^{2})\at{t=0} \partial_s [B (|s \Phi_{ p } +tP_{p}r |^2) ]\,dt\,ds \\&+ \sum_{ j=2,3} \int_{\R^3}dx\int_{[0,1]^2} \frac{t^j}{j!} (\partial_t^{j+1})\at{t=0} \partial_s [ B (|s \Phi_{ p } +tP_{p}r |^2) ]\,dt\,ds \\& + \int_{\R^3}dx \int_{[0,1]^2}\,dt\,ds \int_0^t \partial_\tau^5 \partial_s [B (|s \Phi_{ p } +\tau P_{p}r |^2)] \frac{(t-\tau)^3}{3!} d\tau + E_P (P_{p}r).
 \end{aligned}
\end{equation*}
The 2nd term in the
2nd line is ${2}^{-1}
\langle \nabla^2 E_P(\Phi_{p})P_{p}r, P_{p}r\rangle$ and so in particular the second line is
\begin{equation*}
 2^{-1}\langle (-\Delta + \nabla^2 E_P(\Phi_{p}) -\sum_{j=1,...,4}\lambda_j(p) \Diamond_j)P_{p}r,P_{p}r\rangle = 2^{-1}\langle {\mathbbm{i}} \sigma_3{\mathcal H}_{p}P_{p}r,P_{p}r\rangle.
\end{equation*}
By \eqref{eq:K},
we have
\begin{equation}\label{eq:expKphi0}
K(\Phi_{p}) = d(p) -d({p}^0)+(\lambda (p) - \lambda ({p}^0)) \cdot {p}^0,
\end{equation}
where
\begin{align}\label{def-d-p}
d(p):= E (\Phi_{p})-\lambda (p)\cdot p.
\end{align}
Since $\partial_{p_j}d(p)=-p\cdot \partial_{p_j} \lambda (p)$,
we conclude $K(\Phi_{p}) =O((p- {p}^0)^2)$.
Furthermore, from \eqref{eq:variables} we have
%% \begin{align}&
%% K(\Phi_{p})
%% = \mathfrak{G} \left (\left.(\Pi_j- {p}^0_j) \right |_{j=1}^{4},
%% \left.(\Pi_j (r)) \right |_{j=1}^{4},
%% \left.(\Pi_j((P_{p} -P_{p^0 }) r) \right |_{j=1}^{4}+ \left.(\langle r,
%% \Diamond_j (P_{p} -P_{p^0 }) r\rangle) \right |_{j=1}^{4}\right))
%% \nonumber
%% \\
%% &\text{with $\mathfrak{G} $ smooth and equal to zero at $(0,0,0,0)$ up to second order.}
%% \label{eq:expKphi-old}
%% \end{align}
%% \ac{can we rewrite the above as follows
%% (above, something is wrong with brackets):}
\begin{align}
K(\Phi_{p})
= \mathfrak{G} \Big(
(\Pi_j- {p}^0_j)|_{j=1}^{4},
\Pi_j (r)|_{j=1}^{4},
\left(
\Pi_j\big((P_{p}-P_{p^1})r\big)
+\big\langle r,\Diamond_j(P_{p}-P_{p^1})r\big\rangle
\right)
\big|_{j=1}^{4}
\Big),
\label{eq:expKphi}
\end{align}
with $\mathfrak{G} $ smooth and equal to zero at $(0,0,0)$ up to second order.
%% \ac{Shouldn't it be $(0,0,0)$? It seems that
%% $\mathfrak{G}:\;(\R^4)^3\to \R $
%% takes three arguments only!}
Summing up, we have the following.
\begin{lemma}
 \label{lem:back}
There is an expansion
\begin{align} \label{eq:Kexpansion0} &
K(u)= K(\Phi_{p}) + 2^{-1}\Omega (\mathcal{H}_p P_{p}r, P_{p}r)+ E_P(P_{p}r)\\& + \sum_{d=3,4}
\langle B_{d }(p), (P_{p}r)^{ d} \rangle
 +\int_{\mathbb{R}^3}
B_5 (x, p, r(x)) (P_{p}r)^{ 5}(x)\,dx, \text{ where  for any $k\in \N$:}  \nonumber
\end{align}

\begin{itemize}

\item $ K(\Phi_{p})$ satisfies \eqref{eq:expKphi0}--\eqref{eq:expKphi};

\item $(P_{p}r)^d(x)$ represents $d$-products of components of $P_{p}r$;

 \item
$B_{d
} \in C^\infty (\mathcal{P},
\Sigma_{k} (\mathbb{R}^3, B (
 (\mathbb{R}^{4 })^{\otimes d},\mathbb{R}))) $ for $3\le d \le 4$;
 \item for
$ \zeta \in \mathbb{R}^{4 }$, $B_5$ depends smoothly on its variables,
so that $\forall$  $i\in\N$, there is a constant $C_i$ s.t.
\begin{equation} \label{eq:B5}\begin{aligned} & \| \nabla_{ p, \zeta }
^iB_5 (\cdot,p, \zeta) \|_{\Sigma_{k}(\mathbb{R}^3,
B(
 (\mathbb{R}^{4 })^{\otimes 5},\mathbb{R}))} \le C_i.
 \end{aligned} \end{equation}
\end{itemize}

\end{lemma}
\qed

We will perform a normal form argument on the expansion \eqref{eq:Kexpansion0},
eliminating some terms from the expansion by means of changes of variables.
The first step in a normal forms argument is the diagonalization of the \textit{homological equation},
see \cite[p. 182]{Arnold}, which is discussed in Section~\ref{sec:Normal form}.

\section{Symbols $\mathcal{R}^{i,j}_{k, m}$, $\mathbf{S}^{i,j}_{k, m}$ and restrictions of $ K $ on submanifolds} \label{sect:restrK}

We begin with the following elementary lemma.
 \begin{lemma}
 \label{lem:lambdas1} Set $u=\s (b) \psi $. Then, for $b_R=\Re (b)$ and $b_I=\Im (b)$, we have
 \begin{equation}\label{eq:lambdas11} \begin{aligned} &
 \Pi_{5}(u)= (1-2b_R^2) \Pi_{5}(\psi) -2b_I b_R \Pi_{6}(\psi) -2 \sqrt{1-|b|^2} b_R \Pi_{7}(\psi),
\\&
 \Pi_{6}(u)= -2b_I b_R \Pi_{5}(\psi) +(1-2b_I^2) \Pi_{6}(\psi) -2 \sqrt{1-|b|^2} b_I \Pi_{7}(\psi),
\\&
 \Pi_{7}(u)= 2 \sqrt{1-|b|^2}b_R \Pi_{5}(\psi) +2 \sqrt{1-|b|^2} b_I \Pi_{6}(\psi) + (1-2|b|^2) \Pi_{7}(\psi).
\end{aligned}
\end{equation}
 \end{lemma}
 \proof We have
 \begin{equation*} \begin{aligned} &
 2\Pi_{5}(u)= \langle \sigma_3 \sigma_2 \mathbf{C}u,u \rangle = \langle \s (-b) \sigma_3 \sigma_2 \mathbf{C}\s (b) \psi,\psi \rangle  = \langle \sigma_3 \sigma_2\s (-b) \s \left ( -\mathbf{C}b \right ) \mathbf{C}\psi,\psi \rangle \\& =
 \langle \sigma_3 \sigma_2 \left [ ({1-|b|^2} + b \sigma_2 \mathbf{C} (\mathbf{C} {b})\sigma_2 \mathbf{C})
 - \sqrt{1-|b|^2} (b+ (\mathbf{C} {b})) \sigma_2 \mathbf{C} \right ] \mathbf{C}\psi,\psi \rangle
 \\& =
 \langle \sigma_3 \sigma_2 \left [ {1-b_R^2-\cancel{b_I^2}} - (b_R^2-\cancel{b_I^2}+2{\mathbbm{i}} b_R b_I) - 2
 \sqrt{1-|b|^2} b_R \sigma_2 \mathbf{C} \right ] \mathbf{C}\psi,\psi \rangle
 \\& = (1-2b_R^2) \langle \sigma_3 \sigma_2 \mathbf{C}\psi,\psi \rangle -2b_R b_I
 \langle {\mathbbm{i}} \sigma_3 \sigma_2 \mathbf{C}\psi,\psi \rangle - 2 \sqrt{1-|b|^2} b_R \langle \sigma_3 \psi,\psi \rangle.
\end{aligned}
\end{equation*}
 This yields the formula for $ \Pi_{5}(u)$. By a similar computation
 \begin{equation*} \begin{aligned} &
 2\Pi_{6}(u)= \langle {\mathbbm{i}} \sigma_3 \sigma_2 \mathbf{C}u,u \rangle = \langle\s (-b) {\mathbbm{i}} \sigma_3 \sigma_2 \mathbf{C}\s (b) \psi,\psi \rangle  = \langle {\mathbbm{i}} \sigma_3 \sigma_2 \s (b)  \s (- \mathbf{C} {b})  \mathbf{C}\psi,\psi \rangle \\& =
 \langle {\mathbbm{i}} \sigma_3 \sigma_2 \left [ ({1-|b|^2} - b \sigma_2 \mathbf{C} (\mathbf{C} {b}) \sigma_2 \mathbf{C})
 + \sqrt{1-|b|^2} (b- (\mathbf{C} {b})) \sigma_2 \mathbf{C} \right ] \mathbf{C}\psi,\psi \rangle
 \\& =
 \langle {\mathbbm{i}} \sigma_3 \sigma_2 \left [ {1-\cancel{b_R^2}-b_I^2 } +\cancel{b_R^2}- {b_I^2}+2{\mathbbm{i}}b_R b_I + 2
 {\mathbbm{i}} \sqrt{1-|b|^2} b_I \sigma_2 \mathbf{C} \right ] \mathbf{C}\psi,\psi \rangle
 \\& = (1-2b_I^2) \langle {\mathbbm{i}} \sigma_3 \sigma_2 \mathbf{C}\psi,\psi \rangle -2b_R b_I
 \langle \sigma_3 \sigma_2 \mathbf{C}\psi,\psi \rangle - 2 \sqrt{1-|b|^2} b_I \langle \sigma_3 \psi,\psi \rangle.
\end{aligned}
\end{equation*}
 This yields the formula for $\Pi_{6}(u)$. Finally,  the formula for $\Pi_{7}(u)$ is obtained from
\begin{equation*} \begin{aligned} &
 2\Pi_{7}(u)= \langle \sigma_3 u,u \rangle = \langle \s (-b) \sigma_3 \s (b) \psi,\psi \rangle  = \langle \sigma_3 \s (b) \s (b) \psi,\psi \rangle \\& =
 \langle \sigma_3 \left [ ({1-|b|^2} + b \sigma_2 \mathbf{C} { b} \sigma_2 \mathbf{C})
 + 2\sqrt{1-|b|^2} b \sigma_2 \mathbf{C} \right ] \psi,\psi \rangle
 \\& =
 \langle \sigma_3 \left [ 1-2|b|^2+ 2
 \sqrt{1-|b|^2} b_R \sigma_2 \mathbf{C} + 2
 {\mathbbm{i}} \sqrt{1-|b|^2} b_I \sigma_2 \mathbf{C} \right ] \psi,\psi \rangle
 \\& = (1-2|b|^2) \langle \sigma_3 \psi,\psi \rangle + 2
 \sqrt{1-|b|^2} b_R
 \langle \sigma_3 \sigma_2 \mathbf{C}\psi,\psi \rangle + 2 \sqrt{1-|b|^2} b_I \langle {\mathbbm{i}} \sigma_3 \sigma_2 \mathbf{C} \psi,\psi \rangle.
\end{aligned}
\end{equation*}
\qed

 We introduce the following spaces
 \begin{align}\label{phasespace1} & \Xi _{k} :=\{ ( \Pi _4, \varrho , r) \in  \R _+  \times \R ^7\times   (T^{\bot_{\Omega}}{\mathcal M}_{ p^1 } \cap \Sigma_{k})  \} \text{ for $k\in \Z$,}
\end{align}
where  $\varrho$ is an auxiliary variable which we will use to represent $\Pi  (r)$. We now introduce   two classes of symbols which will be important in the sequel.

	\begin{definition}\label{def:scalSymb} For $ {A} \subset\R^d$     an open set,
$k\in\N_{0}$,
$\mathcal{A}\subset  \Xi _{-k} $
 an open neighborhood of $(p^1_4, 0,0)$,
we say that $ F \in C^{m}(A\times \mathcal{A},\R)$
 is $\mathcal{R}^{i,j}_{k, m}$
 if there exists $C>0$ and an open neighborhood
 $\mathcal{A}'\subset\mathcal{A}$
 of $(p^1_4, 0,0)$ in $ \Xi _{-k}$
 such that
 \begin{equation}\label{eq:scalSymb}
 |F(a,\Pi _4, \varrho, r)|\le C \| r\|_{\Sigma_{-k}}^j (\| r\|
_{\Sigma_{-k}}+|\varrho | +|\Pi _4- p^1_4|) ^{i} \text{ in $I\times
 \mathcal{A}'$}.
\end{equation}
 We will write also $F=\mathcal{R}^{i,j}_{ n,m}$ or
 $F=\mathcal{R}^{i, j}_{ k,m} (a,\Pi _4, \varrho, r)$.
We say $F=\mathcal{R}^{i, j}_{k, \infty}$ if $F=\mathcal{R}^{i,j}_{k, l}$ for all $l\ge m$.
We say $F=\mathcal{R}^{i, j}_{\infty, m} $ if for all $l\ge k$ the above $F$ is the restriction of an
$F \in C^{m}(A\times \mathcal{A}_{l },\R)$ with $\mathcal{A}_l$ an open neighborhood of $(0,0)$ in
$ \R^ 7\times (T^{\bot_{\Omega}}{\mathcal M}_{p^1 } \cap \Sigma_{-l})$ and
$F=\mathcal{R}^{i,j}_{l, m}$. If $F=\mathcal{R}^{i, j}_{\infty, m} $ for any $m$, we set $F=\mathcal{R}^{i, j}_{\infty, \infty} $.
\end{definition}

\begin{remark}\label{rem:scalSymb}  Above we can have $d=0$, that is $A$ missing.
We will also use the following cases: $d=1$  with   a time parameter;
$A$ an open neighborhood of the origin of $ \R   \times \mathbf{su}(2)$. The last case is used only in Appendix \ref{sec:flowspf}.
\end{remark}

\begin{definition}\label{def:opSymb} $T \in C^{m}(A\times \mathcal{A},\Sigma_{k} (\R^3, \C^{2 }))$, with $A\times \mathcal{A}$
like above,
 is $ \mathbf{{S}}^{i,j}_{k,m} $, and we write as above
 $T= \mathbf{{S}}^{i,j}_{k,m}$ or $T= \mathbf{{S}}^{i,j}_{k,m} (a, \Pi _4,\varrho, r)$,
 if there exists $C>0$ and a smaller open neighborhood
 $\mathcal{A}'$ of $(0,0)$ such that
 \begin{equation}\label{eq:opSymb}
 \|T(a,\Pi _4, \varrho, r)\|_{\Sigma_{k}}\le C \| r \|_{\Sigma
_{-k}}^j (\| r\|_{\Sigma_{-k}}+|\varrho | +|\Pi _4- p^1_4|)^{i} \text{ in
 $I\times \mathcal{A}'$}.
\end{equation} We use notation
$T=\mathbf{{S}}^{i,j}_{k,\infty }$, $T=\mathbf{{S}}^{i,j}_{\infty,m}$ and $T=\mathbf{{S}}^{i,j}_{\infty,\infty}$
as
above.

\end{definition}

\begin{lemma}
 \label{lem:repp} On the manifold $\Pi_j= p_j^0$ for $1\le j\le 4$
there exist functions $\resto^{1,2}_{\infty,\infty} $ such that
\begin{equation} \label{eq:p in r}
\begin{aligned}
&
p_j=p_j^0-\Pi_4(r) + \resto^{1,2}_{\infty,\infty}(p ^0_4,\Pi _j(r)\left .  \right | _{j=1}^{4}, r).
\end{aligned}
\end{equation}
\end{lemma}   \proof   By implicit function theorem to  \eqref{eq:variables}   is elementary. \qed

% \begin{lemma}\label{lem:back1} On the manifold $\Pi_j= p_j^0$ for $1\le j\le 4$ we have the expansion \begin{align} \label{eq:Kexpansion01} & K(u)= \resto^{1,0}_{\infty,\infty}(\Pi (r)) + 2^{-1}\Omega (\mathcal{H}_{p^0} r, r)+ \resto^{1,2}_{\infty,\infty}(\Pi (r), r)+ E_P(r)\\ & \qquad\qquad + \sum_{d=3}^4 \langle B_{d }(\Pi (r), r), r^{ d} \rangle +\int_{\mathbb{R}^3} B_5 (x, \Pi (r), r, r(x)) r^{ 5}(x)\,dx \nonumber \end{align} with:

 %\begin{itemize}  \item $B_{d }(\cdot, \cdot) \in C^\infty \left(\mathcal{U}, \Sigma_k \left(\mathbb{R}^3, B ((\mathbb{R}^{4 })^{\otimes d},\mathbb{R})\right)\right) $ for $3\le d \le 4$ were $ \mathcal{U}$ is an open neighborhood of $(0,0)\in \R^7\times \left(T^{\bot_{\Omega}}{\mathcal M}_{p^0 } \cap \Sigma_{-n} \right)$;
% \item for $ \zeta \in \mathbb{R}^{4 }$, $B_5$ depends smoothly on its variables $(\varrho, r,\zeta)$ varying in an open neighborhood $\mathfrak{U}$ of $(0,0,0)$ in $ \R^7\times T^{\bot_{\Omega}}{\mathcal M}_{p^0 } \cap \Sigma_{-n} \times \R^4$ and for any $i\le n$, \begin{equation} \label{eq:B51}\begin{aligned} & \| \nabla_{ \varrho, r, \zeta } ^iB_5 (\cdot,\varrho, r, \zeta) \|_{\Sigma_k(\mathbb{R}^3, B ((\mathbb{R}^{4 })^{\otimes 5},\mathbb{R})} \le C_i \text{ in $ \mathfrak{U} $.} \end{aligned} \end{equation} \end{itemize} \end{lemma} \proof The elementary proof follows by expanding $P_pr= r+(P_p- P_{p^0})r$ and substituting \eqref{eq:p in r} into \eqref{eq:Kexpansion0}. For details see \cite[Lemma 4.4]{Cu0}. \qed

Inside the space parametrized by $(\left.\Pi_j\right |_{j=1}^{4},\tau,b, r)$
we consider
\begin{align}\label{eq:reductionM}
\mathscr{M}_{ 1}^{6}(p^0 ) \text{ defined by } \left.\Pi _j\right |_{j=1}^{6} =\left.  p^{0}_j  \right |_{j=1}^{6}.
\end{align}
Notice that the intersection of $\mathscr{M}_{ 1}^{6}(p^0 )$ with a  small neighborhood of $\{   e^{\mathbbm{i}\vartheta} \Phi _{p^1} : \vartheta \in \R \}$ is a manifold. Indeed, on
the soliton manifold  $\mathcal{M}$ the differential forms $\left . dp _j \right | _{j=1}^{4}, db_R, db_I$  are linearly independent. In the points  of $\mathcal{M}$  formula \eqref{eq:variables} implies  $dp _j =d\Pi  _j $
for $1\le j\le 4$ while the 1st two lines of \eqref{eq:lambdas11}  imply $d\Pi _5=-2 p_4 db_R$  and $d\Pi _5=-2 p_4 db_I$. Hence, since $\Pi _j\in C^\infty (H^1 (\R^3, \C ^2), \R)$,
it follows that  $\left . d\Pi _j \right | _{j=1}^{6}$ are linearly independent in a neighborhood of $\{   e^{\mathbbm{i}\vartheta} \Phi _{p^1} : \vartheta \in \R \}$. Then since
  $\mathscr{M}_{ 1}^{6}(p^0 )$ is defined by $\Pi_j=p^0_j$ for $j\le 6$   we obtain our claim   on $\mathscr{M}_{ 1}^{6}(p^0 )$ for any $p^0 $ sufficiently close to $p^1 $.

   $\mathscr{M}_{ 1}^{6}(p^0 )$ is invariant by the system \eqref{eq:SystK}.
 The following shows that when we factor  $\mathscr{M}_{ 1}^{6}(p^0 )$ by the action  of $\R ^3\times  \mathbb{T}$, the corresponding manifold is parametrized by
  $r\in T^{\bot_{\Omega}}{\mathcal M}_{p^1 } \cap H^1
 (\R^3, \C^2)$.
\begin{lemma}
 \label{lem:repp1}  There exist functions $\resto^{1,2}_{\infty,\infty}(p^0 _4,\Pi (r), r) $ and functions $\resto^{2,0}_{\infty,\infty}(p^0_4 ,\Pi (r)) $ dependent only on $(p^0 _4,\Pi (r))$ s.t. on  $\mathscr{M}_{ 1}^{6}(p^0 )$
   \begin{equation}\label{eq:lambdas13} \begin{aligned} &
 b_R =  ({2p_4^0}) ^{-1} \Pi_{5}(r)+ \resto^{2,0}_{\infty,\infty}(p^0_4,\Pi (r)) + \resto^{1,2}_{\infty,\infty} (p^0_4,\Pi (r), r) \ , \\&
 b_I = ({2p_4^0}) ^{-1} \Pi_{6}(r)+ \resto^{2,0}_{\infty,\infty} (p^0_4,\Pi (r)) + \resto^{1,2}_{\infty,\infty}(p^0_4,\Pi (r), r).
\end{aligned}
\end{equation}
 \end{lemma}
  \noindent \textit{Proof (sketch)}. Since     $\Pi_5=\Pi_6=0$   by
 the first two equations in \eqref{eq:lambdas11}, by
 $\Pi_j (\Phi_p+P_pr)= \Pi_j (P_pr)$ for $j=5,6$ and by $\Pi_7 (\Phi_p+P_pr)=p_4+ \Pi_7 (P_pr)$
 we have
 \begin{equation}\label{eq:lambdas1--3} \begin{aligned} &
 2 \sqrt{1-|b|^2} b_R (p_4+\Pi_{7}(P_pr))= (1-2b_R^2) \Pi_{5}(P_pr) -2b_I b_R \Pi_{6}(P_pr), \\&
 2 \sqrt{1-|b|^2} b_I (p_4+\Pi_{7}(P_pr))= -2b_I b_R \Pi_{5}(P_pr) +(1-2b_I^2) \Pi_{6}(P_pr).
\end{aligned}
\end{equation}
We consider the following change of coordinates, which defines $x_R$ and $x_I$:
 \begin{equation}\label{eq:lambdas1--4}
    2p^0_4 b_R=\Pi_{5}(r) + x_R \text{  and } 2p^0_4 b_I=\Pi_{6}(r) + x_I.
 \end{equation}
 Substitute  in the l.h.s. of \eqref{eq:lambdas1--3}  both \eqref{eq:lambdas1--4} and   \eqref{eq:p in r},  and write
 $\Pi_{j}(P_pr) = \Pi_{j}( r)+  \resto^{1,2}_{\infty,\infty}(p^0_4 ,\Pi (r), r)   $  everywhere in \eqref{eq:lambdas1--3}. Then
 from  the 1st equation in \eqref{eq:lambdas1--3} we get
  \begin{equation}\nonumber  \begin{aligned} &
(1+O(b^2))      \left [  1-\Pi_4(r) /p^0_4 +\Pi_{7}( r)/p^0_4 + \resto^{1,2}_{\infty,\infty}(p^0_4 ,\Pi (r), r) \right ]    (\Pi_{5}(r) + x_R) \\& =\Pi_{5}( r) + O(b^2 \Pi (r)) + \resto^{1,2}_{\infty,\infty}(p^0_4 ,\Pi (r), r))    .
\end{aligned}
\end{equation}
 So, after an obvious cancelation, we have
  \begin{equation}\nonumber  \begin{aligned} &
(1+O(b^2))      \left [  1-\Pi_4(r) /p^0_4 +\Pi_{7}( r)/p^0_4 + \resto^{1,2}_{\infty,\infty}(p^0_4 ,\Pi (r), r)\right ]    x_R \\& = \resto^{2,0}_{\infty,\infty}(\Pi (r) )   + O\( b^2 \Pi (r)\) + \resto^{1,2}_{\infty,\infty}(p^0 ,\Pi (r), r))    .
\end{aligned}
\end{equation}
 which in turn implies for $A=R$
  \begin{equation}\nonumber  \begin{aligned} &
     x_A   = \resto^{2,0}_{\infty,\infty}(p^0_4,\Pi (r) )   + O(b^2 \Pi (r)) + \resto^{1,2}_{\infty,\infty}(p^0_4 ,\Pi (r), r))
\end{aligned}
\end{equation}
 where the big $O$ is smooth. Since a similar equality holds also for $A=I$, substituting again $b$ by means of \eqref{eq:lambdas1--4} and applying the Implicit Function Theorem, we obtain
 \begin{equation}\nonumber  \begin{aligned} &
     x_A   = \resto^{2,0}_{\infty,\infty}(p^0_4 ,\Pi (r) )    + \resto^{1,2}_{\infty,\infty}(p^0_4 ,\Pi (r), r)) \text{ for $A=R,I$.}
\end{aligned}
\end{equation}

\qed

\begin{lemma}
 \label{lem:repp--1}  In   $\mathscr{M}_{ 1}^{6}(p^0 )$  we have
   \begin{equation}\label{eq:Pi7} \begin{aligned} &
   \Pi _7= p_{4}^0+ \Pi_{7}( r) +   \resto^{2,0}_{\infty,\infty}(p^0 ,\Pi (r)) + \resto^{1,2}_{\infty,\infty} (p^0 ,\Pi (r), r) .
\end{aligned}
\end{equation}
 \end{lemma}
\proof By the   3rd identity in  \eqref{eq:lambdas11} and by the definition of $P_p$, we have
  \begin{align*} &
 \Pi_{7} = 2 \sqrt{1-|b|^2}b_R \Pi_{5}(P_pr) +2 \sqrt{1-|b|^2} b_I \Pi_{6}(P_p r) + (1-2|b|^2) ( p_{4}+ \Pi_{7}(P_pr)).
\end{align*}
Using Lemmas \ref{lem:repp}  and \ref{lem:repp1}, we obtain \eqref{eq:Pi7}.  \qed

\section{Expressing $\Omega$ in coordinates } \label{sect:darboux1}

 Normal forms arguments are crucial in the proof of Theorem~\ref{theorem-1.1}.
It is important to settle on a coordinate system where the homological equations look
manageable.  While the symplectic form $\Omega $ has a very simple definition \eqref{eq:Omega}
in terms of
the hermitian structure of $L^2(\R^3, \C^2)$,  it has a rather complicated representation
in terms of the coordinates $(\left.\Pi_j\right |_{j=1}^{4},\tau,b, r)$.
Eventually we will settle on a coordinate system where the symplectic form is equal to the form $\Omega_0$ to be introduced in Section~\ref{sec:speccoo}. In this section we consider some
preliminary material.

\noindent   We consider $\widetilde{\Gamma}:= 2^{-1}\langle {\mathbbm{i}} \sigma_3 u, \,\cdot\, \rangle $.
Using the definition of exterior differentiation it is elementary to show that $
 d \widetilde{\Gamma } =\Omega $.
 We consider now the function
 \begin{equation*}
\psi (u):= 2^{-1}\langle {\mathbbm{i}} \sigma_3 e^{ -{\mathbbm{i}} \sigma_3  \tau \cdot \Diamond }\s (b)\Phi_{ p},u
 \rangle
 \end{equation*}
and   set
$\Gamma:= \widetilde{\Gamma}-d\psi  +d\sum _{j=1,...,4} \Pi _j \tau _j$. Obviously $d\Gamma=\Omega$. We have the following.

\begin{lemma}\label{lem:tildeB} We have
\begin{equation}\label{eq:defgamma}
   \Gamma =  \sum_{j=1,...,4} \tau_j d{\Pi_j} +2^{-1} \Omega (P_p r, dr)+\sum_{j=1,...,4}2^{-1} \Omega (r,P_p \partial_{p_j}P_p r)dp_j+ \varsigma, \text{ where}
\end{equation}
\[
\varsigma := \left (\Pi_5 \frac{b_Rb_I}{\sqrt{1-|b|^2}}-\Pi_6 \frac{1-b_I^2}{\sqrt{1-|b|^2}} - \Pi_7b_I \right)\,db_R + \left (\Pi_5 \frac{1-b_R^2}{\sqrt{1-|b|^2}}-\Pi_6 \frac{b_Rb_I}{\sqrt{1-|b|^2}} + \Pi_7b_R \right)\,db_I
.
\]
\end {lemma}
\proof
The proof is elementary.  The identity operator is $du$,  which can be expanded
\begin{align*}
  du &=  -\sum_{j=1,...,4}  {\mathbbm{i}} \sigma_3 \Diamond_ju d\tau _j +    \sum_{j=1,...,4}  e^{ -{\mathbbm{i}} \sigma_3\tau \cdot \Diamond}\s (b) \partial _{p_j} ( \Phi_{ p}+ P_pr) d{p_j}  \\&
   +  \sum_{A=,R,I}  e^{ -{\mathbbm{i}} \sigma_3 \tau \cdot \Diamond} \partial _{b_A} \s (b) ( \Phi_{ p}+ P_pr) db_A + e^{ -{\mathbbm{i}} \sigma_3\tau \cdot \Diamond}\s (b)P_pdr.
\end{align*}
Then, inserting this in $\widetilde{\Gamma}$ and after some elementary simplification which uses also \eqref{eq:s sym1},  we obtain
\begin{align}
  &  \widetilde{\Gamma} =  2^{-1}\langle {\mathbbm{i}} \sigma_3 u, du \rangle = -\sum_{j=1,...,4} \Pi _j d\tau _j  +  \sum_{A=,R,I} 2^{-1}\langle {\mathbbm{i}} \sigma_3 \s (b)  ( \Phi_{ p}+ P_pr) , \partial _{b_A}\s (b)  ( \Phi_{ p}+ P_pr) \rangle d{b_A}   \nonumber  \\& +  \sum_{j=1,...,4} 2^{-1}\langle {\mathbbm{i}} \sigma_3   ( \Phi_{ p}+ P_pr) ,  \partial _{p_j} ( \Phi_{ p}+ P_pr) \rangle d{p_j}+2^{-1}\langle {\mathbbm{i}} \sigma_3   ( \Phi_{ p}+ P_pr) ,    P_pdr  \rangle .\label{exgamma1}
\end{align}
We have
\begin{align}\label{exgamma2}
&   \text{2nd line of \eqref{exgamma1}}=
  \sum_{j=1,...,4} 2^{-1}  \langle {\mathbbm{i}} \sigma_3   P_{ p}r ,  \partial _{p_j}  P_{ p}r  \rangle  d{p_j} +2^{-1}\langle {\mathbbm{i}} \sigma_3   P_pr  ,    P_pdr  \rangle +d 2^{-1}\langle {\mathbbm{i}} \sigma_3    \Phi_{ p} ,    P_p r  \rangle ,
\end{align}
where we used what follows:
\begin{align*}
&   \langle {\mathbbm{i}} \sigma_3  P_{ p} r ,  \partial _{p_j}  \Phi_{ p}  \rangle  =0 \text{ from the definition of $P_p$;}\\&   \langle {\mathbbm{i}} \sigma_3    \Phi_{ p} ,  \partial _{p_j}  \Phi_{ p}  \rangle =   \langle {\mathbbm{i}} e^{\frac{\mathbbm{i}} 2v\cdot x } \phi_\omega  ,  \partial _{p_j} e^{\frac{\mathbbm{i}} 2v\cdot x } \phi_\omega  \rangle = 0\text{ from formula \eqref{eq:defPhi}.}
\end{align*}
  Hence, by the definition of $\Gamma$,  $\psi (u)$  and $\s (b)$ we obtain
 \begin{align}
 & \Gamma =  \sum_{j=1,...,4}\tau_j d\Pi_j +\sum_{j=1,...,4} 2^{-1}\langle {\mathbbm{i}} \sigma_3 P_{p} r,\partial_{p_j} P_{p} r\rangle dp_j + 2^{-1}\langle {\mathbbm{i}} \sigma_3 r, P_{p}dr \rangle
 \label{eq:brack1} \\& -
 2^{-1}\sum_{A=R,I}\langle {\mathbbm{i}} \sigma_3 \partial_{b_A} (\sqrt{1-| b |^2} + b \sigma_2 \mathbf{C}) (\Phi_{ p } +P_{ p}r), (\sqrt{1-| b |^2} + b \sigma_2 \mathbf{C}) (\Phi_{ p } +P_{ p}r)
 \rangle\,db_A. \nonumber
\end{align}
For $A=R$ the bracket in the
last line equals
\begin{equation*} \begin{aligned}
 &
 \langle {\mathbbm{i}} \sigma_3 \left (\frac{-b_R}{\sqrt{1-| b |^2}} + \sigma_2 \mathbf{C} \right) (\sqrt{1-| b |^2} - b \sigma_2 \mathbf{C}) u, u
 \rangle =\\& \langle {\mathbbm{i}} \sigma_3 \left [ {-b_R} +\overline{b} + \left (
 \sqrt{1-| b |^2} + \frac{ b_R b}{\sqrt{1-| b |^2}} \right) \sigma_2 \mathbf{C} \right ] u, u
 \rangle =\\& \langle {\mathbbm{i}} \sigma_3 \left [ -{\mathbbm{i}} b_I +
 \frac{ 1-\cancel{b^2_R}-b^2_I + \cancel{b^2_R}+{\mathbbm{i}} b_R b_I}{\sqrt{1-| b |^2}} \sigma_2 \mathbf{C} \right ] u, u
 \rangle =\\&  \langle {\mathbbm{i}} \sigma_3 \left [ -{\mathbbm{i}} b_I +\frac{ 1- b^2_I }{\sqrt{1-| b |^2}} \sigma_2 \mathbf{C} + \frac{ b_R b_I}{\sqrt{1-| b |^2}} {\mathbbm{i}} \sigma_2 \mathbf{C} \right ] u, u
 \rangle   = b_I \Pi_7 + \frac{ 1- b^2_I }{\sqrt{1-| b |^2}} \Pi_6 - \frac{ b_R b_I}{\sqrt{1-| b |^2}}\Pi_5.
\end{aligned}
\end{equation*}
For $A=I$ the bracket in the
last line of \eqref{eq:brack1} equals
\begin{equation*} \begin{aligned}
 &
 \langle {\mathbbm{i}} \sigma_3 \left (\frac{-b_I}{\sqrt{1-| b |^2}} + {\mathbbm{i}} \sigma_2 \mathbf{C} \right) (\sqrt{1-| b |^2} - b \sigma_2 \mathbf{C}) u, u
 \rangle =\\& \langle {\mathbbm{i}} \sigma_3 \left [ {-b_I} +{\mathbbm{i}} \overline{b} + \left (
 {\mathbbm{i}} \sqrt{1-| b |^2} + \frac{ b_I b}{\sqrt{1-| b |^2}} \right) \sigma_2 \mathbf{C} \right ] u, u
 \rangle =\\& \langle {\mathbbm{i}} \sigma_3 \left [ {\mathbbm{i}} { b_R} +
 \frac{ {\mathbbm{i}} (1-b^2_R -\cancel{b^2_I}) + \cancel{{\mathbbm{i}} b^2_I}+ b_R b_I}{\sqrt{1-| b |^2}} \sigma_2 \mathbf{C} \right ] u, u
 \rangle =\\& \langle {\mathbbm{i}} \sigma_3 \left [ {\mathbbm{i}} b_R +\frac{ 1- b^2_R }{\sqrt{1-| b |^2}} {\mathbbm{i}} \sigma_2 \mathbf{C} + \frac{ b_R b_I}{\sqrt{1-| b |^2}} \sigma_2 \mathbf{C} \right ] u, u
 \rangle   =- b_R \Pi_7
 - \frac{ 1- b^2_R }{\sqrt{1-| b |^2}} \Pi_5 + \frac{ b_R b_I}{\sqrt{1-| b |^2}}\Pi_6.
\end{aligned}
\end{equation*}
 This completes the proof of Lemma~\ref{lem:tildeB}. \qed

\begin{lemma}
 \label{lem:repp2} Consider the immersion $ i :  \mathscr{M}_{ 1}^{6}(p^0 )\hookrightarrow H^1(\R^3, \C ^2)$
 and the pullback $i^*\Gamma $, which by an abuse of notation we will still denote by $\Gamma$. We have   \begin{align}  \Gamma = i^*\Gamma &= 2^{-1} \Omega (r, dr)+\langle \mathcal{R}^{0,2}_{\infty,\infty}(p^0_4 ,\Pi (r), r) \cdot \Diamond r + \mathbf{S}^{1,1}_{\infty,\infty}(p^0_4 ,\Pi (r), r), dr \rangle + \Pi_7\varpi  \text{ where} \label{eq:eqB} \\    \varpi  &= (b_R\,db_I-b_I\,db_R)  =\frac{1}{4(p_4^0)^2} (\Pi_5(r)d \Pi_6(r)  -\Pi_6(r) d \Pi_5(r))  + \resto^{2,0}_{\infty,\infty}(p^0_4 ,\Pi (r))d\Pi (r) \nonumber\\&+ \langle \mathbf{ S}^{2,1}_{\infty,\infty} (p^0_4 ,\Pi (r), r),dr \rangle. \label{eq:eqB1}
 \end{align}
 \end{lemma}
\proof   The starting point is formula \eqref{eq:defgamma} for $\Gamma$.  Obviously for the restrictions we have
$  \left . d\Pi _k \right | _{ \mathscr{M}_{ 1}^{6}(p^0 )} =0$
 for $1\le k \le 6$.
So that the 1st summation  in the r.h.s. of \eqref{eq:defgamma} contributes 0.

\noindent Next, notice that
for $1\le j\le 4$ from \eqref{eq:p in r} we obtain
 \begin{equation*} \begin{aligned} &
 d p_j = -\langle \Diamond_jr +\mathbf{S}^{1,1}_{\infty,\infty},dr \rangle +\sum_{k\le 4} \mathcal{R}^{0,2}_{\infty,\infty} dp_k,
\end{aligned}
\end{equation*}
which, solved in terms of the $dp_j$'s, gives
\begin{equation}\label{eq:dpj} \begin{aligned} &
 d p_j = -\sum_{k\le 4} \langle (\delta_{jk} +\mathcal{R}^{0,2}_{\infty,\infty})\Diamond_kr +\mathbf{S}^{1,1}_{\infty,\infty},dr \rangle.
\end{aligned}
\end{equation}
Substituting $dp_j$ by \eqref{eq:dpj}  in \eqref{eq:defgamma}  and using  and $P_pr=r+\mathbf{S}^{1,1}_{\infty ,\infty}(p^0 ,\Pi (r), r)$ on $ \mathscr{M}_{ 1}^{6}(p^0 ) $, we obtain terms like the 2nd in the r.h.s. of \eqref{eq:eqB}.

\noindent Finally, by $\Pi _5=\Pi _6$, we obtain $\varsigma = \Pi _7 \varpi$. To get the r.h.s. in \eqref{eq:eqB1}, we use
 the following formulas,
 \begin{equation}\label{eq:lambdas14} \begin{aligned} &
 db_R = ({2p_4^0}) ^{-1}\langle \sigma_3\sigma_2 \mathbf{C} r, dr\rangle + \resto^{1,0}_{\infty,\infty}(p^0_4 ,\Pi (r))d\Pi (r)+ \langle \mathbf{ S}^{1,1}_{\infty,\infty},dr \rangle,
\\&
 db_I = ({2p_4^0}) ^{-1}\langle {\mathbbm{i}} \sigma_3\sigma_2 \mathbf{C} r, dr\rangle + \resto^{1,0}_{\infty,\infty}(p^0_4 ,\Pi (r))d\Pi (r)+ \langle \mathbf{ S}^{1,1}_{\infty,\infty},dr \rangle,
\end{aligned}
\end{equation}
where $ \resto^{1,0}_{\infty,\infty}(p^0_4 ,\Pi (r))d\Pi (r)$ stands for $ \sum_{j=1,...,7} \resto^{1,0}_{\infty,\infty}(p^0_4 ,\Pi (r))d\Pi_j (r)$
with different real-valued functions
from the class $\resto^{1,0}_{\infty,\infty}(p^0_4 ,\Pi (r))$. Formulas \eqref{eq:lambdas14} are obtained differentiating in \eqref{eq:lambdas13}.

\qed

% \begin{lemma} \label{lem:replace1} Consider in $\mathscr{M}_{ 1}^{6}$ the   differential form  \begin{equation} \label{eq:eqtildeB} \begin{aligned} &\Gamma ' := 2^{-1} \Omega (r, dr)+\langle \mathcal{R}^{0,2}_{\infty,\infty}(\Pi (r), r) \cdot \Diamond r + \mathbf{S}^{1,1}_{\infty,\infty}(\Pi (r), r), dr \rangle + p^{0}_7 \varpi \end{aligned} \end{equation} obtained replacing  $\Pi _7$ with the constant value $ p^0_7$ in the last term of  \eqref{eq:eqB} and consider $ {\Omega}'=d \Gamma '$.    Let $X_K$  resp.  $X_K'$ be  the hamiltonian field of $K$ for $\Omega $ resp. $\Omega '$. Then $X_K=X_K'$ at the points where $\Pi _7=p_7^0$. \end{lemma} \proof We have $  i^*\Gamma =  {\Gamma}' + (\Pi _7-p_7^0)  \varpi$  and so \begin{align*} & i^*\Omega =  {\Omega}' + d \Pi _7  \wedge \varpi +(\Pi _7-p_7^0) d \varpi =  {\Omega}' + d \Pi _7  \wedge \varpi \end{align*} at the points where  $\Pi _7=p_7^0$.  Then at the points where $\Pi _7=p_7^0$ \begin{align*} &  i_{X_K} i^*\Omega  = i_{X_K}  \Omega '  +    \varpi \ i_{X_K} d\Pi _7  - d\Pi _7 \ i_{X_K} \varpi  =  i_{X_K}  \Omega '    - d\Pi _7 \ i_{X_K} \varpi ,  \end{align*} where we used $i_{X_K} d\Pi _7 =\{ \Pi _7, K  \} =0.$

Substituting $\Pi _7$ by   \eqref{eq:Pi7}    in \eqref{eq:eqB}  we obtain
   \begin{align}  &\Gamma   =2^{-1} \Omega (r, dr) +\langle   \mathbf{S}^{1,1}_{\infty,\infty}(p^0_4 ,\Pi (r), r), dr \rangle   \nonumber\\&  +({4p_4^0}) ^{-1}(\Pi_5(r)d \Pi_6(r)  -\Pi_6(r) d \Pi_5(r))   +\( \resto^{2,0}_{\infty,\infty}(p^0_4 ,\Pi (r)) + \mathcal{R}^{0,2}_{\infty,\infty}(p^0_4 ,\Pi (r), r)  \)d\Pi (r) . \label{eq:B0pr}\end{align}

\section{Spectral coordinates associated to $\mathcal{H}_{ p^1}$}
\label{sec:speccoo}

    By assumption  $p^1=p(\omega^1,0)$.
Recall that the operator $ \mathcal{H}_{ p^1} $   defined in $L^ 2(\R^3, \C^2)$ is not $\C$-linear (because of $ \mathfrak{L}^{(1)}_{\omega ^1 }$), but rather
$\R$-linear. To make it $\C$-linear,
%% we introduce a new imaginary unit $\im$ and consider
%% \begin{equation*}
%%  \C =\{ a +\im b:\;a,\,b\in \R \}.
%% \end{equation*}
%% \ac{$\tilde\C$??}
%% Then
we consider the complexification
\[
L^ 2(\R^3, \C^2)\otimes_\R \C.
\]
To avoid the confusion between $\C$ in the left factor and
$\C$ on the right, we will use $\im$ to denote the imaginary unit
in the latter space; that is,
given
$u\in L^ 2(\R^3, \C^2)$,
we will have
$u\otimes(a+\im b)\in L^ 2(\R^3, \C^2)\otimes_\R \C$.
Notice that the domain of
$ \mathcal{H}_{ p^1} $ in $L^ 2(\R^3, \C^2) $ is $ H^2(\R^3, \C^2) $;
we extend it to $L^ 2(\R^3, \C^2)\otimes_\R \C$ with the domain $H^2(\R^3, \C^2)\otimes_\R \C$
by setting $ \mathcal{H}_{ p^1}(v\otimes z)=(\mathcal{H}_{ p^1} v) \otimes z $.

\noindent We extend the bilinear form $\langle\,,\,\rangle $ defined in \eqref{eq:hermitian}
to a $\C$-bilinear form on $L^ 2(\R^3, \C^2)\otimes_\R \C$ by
\begin{equation*}
\langle u\otimes z, v\otimes \zeta\rangle = z\zeta \langle u, v \rangle,
\qquad
u,\,v\in L^ 2(\R^3, \C^2),
\quad
z,\,\zeta\in\C.
\end{equation*}
We also extend $\Omega$ onto $L^2(\R^3,\C^2)\otimes_\R\C$,
setting $\Omega (X,Y)=\langle {\mathbbm{i}} \sigma_3 X,Y \rangle$.
Then the decomposition \eqref{eq:dirsum} extends into \begin{equation}\label{eq:dirsumc} \begin{aligned} &
 L^2(\R^3, \C^4) \otimes_\R \C = (T_{ \Phi_{p^1} }^{\perp_\Omega}\mathcal{M}\otimes_\R \C) \oplus (T_{ \Phi_{p^1} }\mathcal{M}\cap H^1(\R^3, \C^2))  \otimes_\R \C.
\end{aligned}
\end{equation}
Note that the extention of
$ \mathcal{H}_{ p^1} $ onto $L^ 2(\R^3, \C^2)\otimes_\R \C$
is such that its action preserves the decomposition \eqref{eq:dirsumc}.
 The complex conjugation on $ L^ 2(\R^3, \C^2)\otimes_\R \C$
is defined by
$ \overline{v \otimes z} := v \otimes \overline{z} $.

\noindent Notice that if $\im \mathcal{H}_{ p^1} \xi_l =\mathbf{e}_l  \xi_l $ with $\mathbf{e}_l  >0$, then by complex conjugation we obtain
$\im \mathcal{H}_{ p^1} \overline{\xi}_l   =-\mathbf{e}_l  \overline{\xi}_l  $.

\noindent By Weyl's theorem,
$\sigma_e(\im \mathcal{H}_{ p^1}) = (-\infty, - \omega^1]\cup [ \omega^1, \infty)$.  We assume spectral stability, i.e.
$\sigma_e(\im \mathcal{H}_{ p^1}) \subset \R$.
We   assume the set of eigenvalues
$\sigma_p( \im \mathcal{H}_{ p^1}
) \subset (- \omega^1, \omega^1)$, $\pm \omega^1$ are not resonances and the following.

\begin{itemize} \item[(H6)] For any   $ \mathfrak{e} \in \sigma_p(\im \mathcal{H}_{p^1}) \backslash \{ 0 \}$, algebraic and geometric multiplicities
 coincide and are finite.

 \item [(H7)]
%% There is a number $\mathbf{n}\ge 1$
%% and positive numbers $0<\mathbf{e}
%%_1 \le \mathbf{e}_2 \le...\le\mathbf{e}_\mathbf{n} <\omega^0$
%% such that $\sigma_p(\mathcal{H}_{ p^0})$ consists exactly of the numbers
%% $\pm \im \mathbf{e}_\ell $ and $0$.
%% \ac{New version:}
There is a number $\mathfrak{N}\in\N$
and positive numbers $0<\mathbf{e}_1 < \mathbf{e}_2 <\ldots<\mathbf{e}_\mathfrak{N} <\omega^1$
such that $\sigma_p(\mathcal{H}_{ p^1})$ consists exactly of the numbers
$\pm \im \mathbf{e}_\ell $ and $0$. Furthermore, the points $\pm  \im \omega^1$ are not resonances ( that is, if  $\mathcal{H}_{ p^1} \Theta = \pm  \im \omega^1 \Theta$
for one of the two signs, and if $\langle x\rangle  \Theta \in L^\infty, $ then $\Theta =0$).

Denote $d_\ell:=\dim\ker (\mathcal{H}_{ p^1} -\im \mathbf{e}_\ell)$
and let
\[
\mathbf{n}:=\sum\sb{\ell=1,..., \mathfrak{N}}d_\ell.
\]

\item[(H8)]
  %% Denote by $\mathbf{n}_0=0< \mathbf{n}_1<...<\mathbf{n}_{l_0}=\mathbf{n}$
  %% the integers such that
  %% $\mathbf{e}_\ell = \mathbf{e}_i $ exactly for $i$ and $\ell$
  %% both in $(\mathbf{n}_l, \mathbf{n}_{l+1}]$ for some $l\le l_0$.
  %% It follows that
  %% $d_\ell:=\dim\ker (\mathcal{H}_{ p^0} -\im \mathbf{e}_\ell)
  %% =\mathbf{n}_{\ell+1}-\mathbf{n}_\ell$.
 We define
\begin{equation}\label{eq:h8}
 \mathbf{N}:=\sup_\ell
\inf \{ n\in\N :\; n \mathbf{e}_\ell \in \sigma_e(\im\mathcal{H}_{ p^1}) \}-1.
\end{equation}
If $\mathbf{e}_{\ell_1} <...<\mathbf{e}_{\ell_i} $ are distinct and $\mu\in \Z^i$ satisfies
 $|\mu|:=\sum_{j=1}^{i}\mu_{j}\leq 4\mathbf{N} +4$,
 we assume that
$$
\mu_1\mathbf{e}_{\ell_1} +\dots +\mu_k\mathbf{e}_{\ell_i} =0 \quad\iff\quad \mu=0.
$$

\end{itemize}

 It is easy to prove  the symmetry of $\sigma_p(\im \mathcal{H}_{ p^ 1})\subset\R$ around $0$.  We have
 \begin{align*}
    \ker (\im \mathcal{H}_{ p^ 1} \mp \mathbf{e}_l
) \big) \subset \mathcal{S}(\R^3, \C^2) \otimes _\R \C
 \end{align*}
 and using $\Omega$ we consider the set $X_c\subset \mathcal{S}'(\R^3, \C^2) \otimes _\R \C$ defined by
  \begin{align}\label{eq:defXc}
   X_c:=\left [  \( T\mathcal{M}_{ \Phi_{p^1} }   \otimes_\R \C\)  \oplus  _\pm
\oplus_{l=1}^{N}
\left (\ker (\im \mathcal{H}_{ p^ 1} \mp \mathbf{e}_l
) \big) \right)   \right ] ^{\perp_\Omega}.
 \end{align}
 It is possible to prove the   following decomposition:
\begin{align} \label{eq:spectraldecomp} &
 (T_{ \Phi_{p^1} }^{\perp_\Omega}\mathcal{M}\cap L^2(\R^3, \C^2))  \otimes_\R \C
= \big (\oplus_{\pm}\oplus_{l=1}^{N} \ker (\im \mathcal{H}_{ p^ 1} \mp \mathbf{e}_l
) \big) \oplus \(X_c   \cap \(  L^2(\R^3, \C^2)   \otimes_\R \C \)    \)  .
\end{align}
The decomposition in \eqref{eq:spectraldecomp} is $ \mathcal{H}_{p^ 1}$-invariant.

\noindent
Consider now
the coordinate $r\in T_{ p^1 }^{\perp_\Omega}\mathcal{M} \cap L^ 2(\R^3, \C^2)$
from the coordinate system \eqref{eq:coordinate};
it corresponds to the second summand in \eqref{eq:dirsumc}.
Then, considered as an element from $L^ 2(\R^3, \C^2)\otimes_\R \C$, it can be
 decomposed into
 \begin{equation}
 \label{eq:decomp2}
 r (x) =\sum_{l=1,..., \mathbf{{n}}}z_l \xi_l (x) +
\sum_{l=1,..., \mathbf{{n}}}\overline{z}_l \overline{{\xi }}_l (x)
+ f (x), \quad f \in X_c \text{ with $f= \overline{{f}} $ },
\end{equation}
with $\xi_l$ eigenfunctions
of $\mathcal{H}_{ p^ 1}$
corresponding to $\im \mathbf{e}_l $.
We claim that it is possible to choose them
so that
\begin{equation} \label{eq:norm1} \begin{aligned} &
 \langle {\mathbbm{i}} \sigma_3 \xi_i, \xi_l\rangle = \langle {\mathbbm{i}} \sigma_3 \xi_i, f\rangle = 0\text{ for all $i,l$ and for all $f\in X_c $},
\\& \langle {\mathbbm{i}} \sigma_3 \xi_i, \overline{\xi}_l \rangle=-\im \delta_{il} \text{ for all $i,l$}.
\end{aligned}\end{equation}
 To see the second line, observe that on one hand for $\Theta \in (T_{ \Phi_{p^ 1} }^{\perp_\Omega} \mathcal{M}\otimes_\R \C) \backslash \{ 0 \} $ we have
$ \langle {\mathbbm{i}} \sigma_3 \mathcal{H}_{ p^ 1}\Theta , \overline{\Theta } \rangle>0$. Indeed, for $\Theta  =(\Theta _1, \Theta _2) $ we have
\begin{equation*} \begin{aligned} & \langle {\mathbbm{i}} \sigma_3 \mathcal{H}_{ p^ 1}\Theta, \overline{\Theta} \rangle= \langle {\mathbbm{i}} \sigma_3 \mathcal{H}_{ p^ 1}\Theta, \overline{{\Theta}} \rangle = \langle \mathfrak{L}_{\omega^1}^{(1)}\Theta _1, \Theta _1^*\rangle + \langle \mathfrak{L}_{\omega^1}^{(2)}\Theta _2, \overline{\Theta}_2 \rangle
\end{aligned}\end{equation*}
 with $\langle \Theta _2, \phi_{\omega^1} \rangle =0$, which implies $ \langle \mathfrak{L}_{\omega^1}^{(2)}\Theta _2, \overline{\Theta }_2 \rangle >c_0 \|\Theta _2 \|^2_{L^2}$ and with $\langle \Theta_1, \partial_a\phi_{\omega^1} \rangle =\langle \Theta _1, x_a\phi_{\omega^1} \rangle =\langle \Theta _1, {\mathbbm{i}} \phi_{\omega^1} \rangle =0$ which implies $\langle \mathfrak{L}_{\omega^1}^{(1)}\Theta _1, \overline{{ \Theta }}_1 \rangle >c_0 \| \Theta _1\|^2_{L^2}$, for a fixed $c_0>0$.
 On the other hand,
 \begin{equation*} \begin{aligned} & 0< \langle {\mathbbm{i}} \sigma_3 \mathcal{H}_{ p^ 1}\xi_i,\overline{\xi}_i \rangle= \im \mathbf{e}_i \langle {\mathbbm{i}} \sigma_3 \xi_i,\overline{\xi}_i \rangle.
\end{aligned}\end{equation*}

\noindent It is then possible to choose $\xi_i$ so that \eqref{eq:norm1} is true. Notice that
 \eqref{eq:norm1} means that the nonzero eigenvalues have positive \textit{Krein signature}. This proves the second line of \eqref{eq:norm1}. The proof of the 1st line is elementary.

\noindent By \eqref{eq:decomp2} and \eqref{eq:norm1},
 we have
\begin{equation}
 \label{eq:H2}
  {2}^{-1} \langle {\mathbbm{i}} \sigma_3 \mathcal{H}_{ p^1 } r, r\rangle = \sum_{l=1,...,\mathbf{n}} \mathbf{e}_l|z_l|^2
+ {2}^{-1}\langle {\mathbbm{i}} \sigma_3 \mathcal{H}_{ p^1 } f, f\rangle =:H_2.
\end{equation}
 In terms of
 $(z,f)$, the Fr\'echet derivative $dr$ can be expressed as
\begin{equation}
 \label{eq:decomp20}
 dr=\sum_{l=1,...,\mathbf{n}}(dz_l \xi_l+d\overline{z}_l \overline{{\xi}}_l) +df,
\end{equation}
and by \eqref{eq:norm1} we have \begin{equation}
 \label{eq:OmegaCoo}
 2^{-1} \langle {\mathbbm{i}} \sigma_3 r, dr\rangle = 2^{-1} \im \sum_{l=1,...,\mathbf{n}} (\overline{z}_l d {z}_l -z_l d \overline{z}_l)
+ 2^{-1} \langle {\mathbbm{i}} \sigma_3 f, df\rangle.
\end{equation}
Notice now that, in terms of \eqref{eq:decomp2} and \eqref{eq:decomp20},
\begin{equation*}
 d\Pi_j(r)
= \langle \Diamond_j(z\xi + \overline{z}\ \overline{\xi} +f),\ \xi\,dz + \overline{\xi}\,d\overline{z}+df\rangle
= \sum_{l=1,...,\mathbf{n}} (\resto^{0,1}_{\infty,\infty}\,dz_l+\resto^{0,1}_{\infty,\infty}\,d\overline{z}_l)+\langle \Diamond_jf + \mathbf{S}^{0,1}_{\infty,\infty}, df\rangle.
\end{equation*}
Hence, we obtain from \eqref{eq:B0pr}:
\begin{equation} \label{eq:B0pr2}\begin{aligned} &\Gamma  = \Gamma_0+ \sum_{l=1,...,\mathbf{n}} (\resto^{1,1}_{\infty,\infty}\,dz_l+\resto^{1,1}_{\infty,\infty}\,d\overline{z}_l)+\langle \sum _{j\le 7} \resto^{0,2}_{\infty,\infty}\Diamond_jf +\mathbf{S}^{1,1}_{\infty,\infty}, df\rangle,
\quad
\text{where}\\ &\Gamma_0:= 2^{-1} \im \sum_{l=1,...,\mathbf{n}} (\overline{z}_l d {z}_l -z_l d \overline{z}_l)
+ 2^{-1} \langle {\mathbbm{i}} \sigma_3 f, df\rangle +\sum_{j\le 7} \resto^{1,0}_{\infty,\infty}(p^0 ,\Pi (f)) \langle\Diamond_jf, df\rangle
.
 \end{aligned}
\end{equation}
Then
\begin{equation}\label{eq:modomega0}
 \Omega_0:=d \Gamma_0= - \im \sum_{l=1,...,\mathbf{n}} \,dz_l\wedge d \overline{z}_l
+ \langle {\mathbbm{i}} \sigma_3 df, df\rangle + \sum_{j,k} \resto^{0,0}_{\infty,\infty}(p^0 ,\Pi (f)) \langle\Diamond_kf, df\rangle \wedge \langle\Diamond_jf, df\rangle,
\end{equation}
 and, schematically, and using in the last line
$\partial_\rho \left. \mathbf{S}^{1,1}_{\infty,\infty} \right |_{(p^0_4,\rho,z, f)=(p^0_4,\Pi (f),z, f)} = \mathbf{S}^{0,1}_{\infty,\infty} $  and defining
\begin{equation}\label{eq:part f mod}
    \widehat{\nabla}_fF  (\Pi (f),f):={\nabla}_fF - \partial _{\Pi (f)}F \cdot  {\nabla}_f\Pi (f),
\end{equation}
\begin{equation} \label{eq:omedgamega0}\begin{aligned} \Omega   - \Omega_0=& \resto^{1,0}_{\infty,\infty} \,dz \wedge \,dz + \langle \widehat{\nabla}_f \mathbf{S}^{1,1}_{\infty,\infty}  df, df \rangle  \\&+dz \wedge \langle \sum _{j\le 7} \resto^{0,1}_{\infty,\infty}\Diamond_jf+\mathbf{S}^{1,0}_{\infty,\infty}, df \rangle + d\Pi (f) \wedge \langle \mathbf{S}^{0,1}_{\infty,\infty}, df \rangle .
\end{aligned}
\end{equation}
  We will transform $\Omega  $ into $\Omega_0 $ by means of
  the Darboux Theorem, performed in a non-abstract way, to make sure that the coordinate transformation is as in Lemma~\ref{lem:ODE}.

\section{Flows } \label{sect:flows}

% In analogy to the space $ \Xi _{k}$ in \eqref{phasespace1} we introduce  \begin{align}\label{phasespace2} & \widetilde{\Xi} _{k} :=\{ (\Pi _4, \varrho ,z, f) \in  \R _+  \times \R ^7\times \C ^{\mathbf{n}}\times   (X_c\cap \Sigma_{k})  \} \text{ for $k\in \Z$ } \end{align}that is simply obtained considering the splitting  \eqref{eq:decomp2} of $r$ in terms of $(z,f)$.
The  following lemma is a consequence of of Lemma \ref{lem:ODEapp}  in Appendix \ref{sec:flowspf}.

\begin{lemma} \label{lem:ODE} For
$n,M,M_0, {s}, {s}',k,l\in  \N_0$ with $1\le l\le M$, for $a\in A  $   a parameter,  with $A$ an open subset in $\R ^d$, $\Pi _4$ another parameter   and for $\widetilde{\varepsilon}_0>0$,
consider
\begin{equation} \label{eq:ODE}
 \left\{\begin{array}{l} \dot z (t)= \resto^{0,M_0 }_{n,M}(t,a,\Pi _4,  \Pi (f),z,f)
 \\[1ex]   \dot f (t)
	= {\mathbbm{i}} \sigma_3\sum_{j\le 7} \resto^{0,M_0+1}_{n,M}(t,a,\Pi _4, \Pi (f),z,f)  \Diamond_j f +
\mathbf{S}^{i,M_0 }_{n,M}(t,a,\Pi _4, \Pi (f),z,f),\end{array}\right .  \end{equation}
with the coefficients defined for $|t|<5$, $|\Pi (f)| <\widetilde{\varepsilon}_0$, $|z|<\widetilde{\varepsilon}_0$,  $\| r\|_{\Sigma_{-n}} <\widetilde{\varepsilon} $ and $|\Pi _4- p^1_4|\le \widetilde{\varepsilon}
_0$.

\noindent  Let $k\in \Z\cap [0,n-(l+1) ]$ and set for
$s''\ge 1$ and $\varepsilon >0$ \begin{equation} \label{eq:domain0}
\begin{aligned} \U_{\varepsilon,k}^{ {s}''}
:= &\{ (z,f)  \in  \C ^ \mathbf{n} \times \( X_c \cap    \Sigma_{ {s}''}  \) \
 :\; |z|+  \| f\|_{\Sigma_{-k }} + |\Pi (f)| \le \varepsilon \}.
\end{aligned} \end{equation}  Let $a_0\in A$.
Then for $\varepsilon >0$ small enough, \eqref{eq:ODE} 	 defines a flow
$(z^t,f^t) =\mathfrak{F}_{t } (z,f)$ with
\begin{align} &  z ^t= \resto^{0,M_0  }_{ n- l-1, l}   (*)    \text{ , where }  *=(t,a, \Pi _4,\Pi (f),z,f) \ ,\label{eq:ODE1}
\\& f^t
= e^{{\mathbbm{i}} \sigma_3  \sum _{j=1}^4\resto^{0,M_0+1 }_{ n- l-1, l}(*)
   \Diamond _j }  T(  e^{   \sum _{i=1}^3  \resto^{0,M_0+1 }_{ n- l-1, l}(*)
 {\mathbbm{i}} \sigma _i  }  )   \( f+\mathbf{{S}}^{i,M_0 }_{ n- l-1,l}(*)\), \nonumber
\end{align} where for
\begin{equation}\label{eq:index1}
 \text{$n- l-1 \ge  {s}' \ge  {s}+l \ge l $ and $k\in
\Z\cap [0,n- l-1 ]$}
\end{equation} and
 for $ \varepsilon_1> \varepsilon_2> 0$ sufficiently small we have
 \begin{equation} \label{eq:reg1}\begin{aligned} &\mathfrak{F}_{t }
 \in C^l((-4,4)\times D _{\R^d}(a_0, \varepsilon_2 ) \times \U_{\varepsilon_2,k}^{ {s}'}, \U_{\varepsilon_1,k}^{ {s}}
)
. \end{aligned} \end{equation}
 	
 \end{lemma}

 \qed

\noindent In \eqref{eq:reg1} the $C^l$ regularity comes at the cost of a loss
of $l$ derivatives in the space $\Sigma_{ {s}''}$, which is accounted
by $s'\ge  {s}+l$.
In Proposition~\ref{th:main} we will need the following elementary technical lemma.
\begin{lemma} \label{lem:ODEbis} Consider two systems for $\ell =1,2$:
\begin{equation} \label{eq:ODEbis}\begin{aligned} &
 \left\{\begin{array}{l} \dot z (t)=  \mathcal{B}^{(\ell)} (t,a, \Pi _4,\Pi (f),z,f)
 \\[1ex]   \dot f (t)
	= {\mathbbm{i}} \sigma_3\sum_{j\le 7}  \mathcal{A }_j^{(\ell)}(t,a,\Pi _4, \Pi (f),z,f)  \Diamond_j f +
\mathcal{D}^{(\ell)}(t,a, \Pi _4,\Pi (f),z,f),\end{array}\right .
\end{aligned} \end{equation}
with the hypotheses of Lemma~\ref{lem:ODE} satisfied, and suppose that
\begin{equation} \label{eq:ODEbis1}\begin{aligned} &  \mathcal{B}^{(1)}(t,a, \Pi _4,\Pi (f),z,f)-\mathcal{B}^{(2)} (t,a,\Pi _4, \Pi (f),z,f) =  \mathcal{R} ^{0,M_0+1}_{n,M}(t,a,\Pi _4, \Pi (f),z,f) \\&
 \mathcal{D}^{(1)}(t,a,\Pi _4, \Pi (f),z,f)-\mathcal{D}^{(2)}(t,a, \Pi _4,\Pi (f),z,f)=\mathbf{S}^{0,M_0+1}_{n,M}(t,a,\Pi _4, \Pi (f),z,f).
\end{aligned} \end{equation}
Let $(z,f)  \mapsto (z^t _{(\ell)} , f^t _{(\ell)})  $ with $\ell =1,2$ be the two   flows. Then for $ {s}, {s}'$ as in Lemma \ref{lem:ODE}
\begin{equation} \label{eq:ODE1bis}\begin{aligned} &
| z^1_{(1)}- z^1_{(2)} |+ \| f^1_{(1)}- f^1_{(2)}\|_{\Sigma_{- {s}' }} \le C \( |z|+ \| f\|_{\Sigma_{ - {s} }} \) ^{M_0+1}.
	 \end{aligned} \end{equation}
\end{lemma}

\proof  For the proof see Lemma \ref{lem:ODEbisapp}.\qed

\begin{lemma} \label{lem:ODEtris} Under the hypotheses and notation of Lemma \ref{lem:ODEbis} we have
\begin{equation} \label{eq:ODE1tris}\begin{aligned} &
 \Pi _j( f^1_{(1)}) -  \Pi _j( f^1_{(2)}) = \resto ^{0,M_0+2}_{ n- l-3,l}(a,\Pi _4,\Pi(f),z,f) \text{ for }j=1,2,3,4.
	 \end{aligned} \end{equation}
\end{lemma}
{\it Proof (sketch)} For $\ell =1,2$  and $j=1,2,3,4$  we have
\begin{equation} \label{eq:1ODEtris} \begin{aligned} &
 \Pi _j( f^1_{(\ell)}) =   \Pi _j(  f + \mathbf{S} ^{(\ell)})  =  \Pi _j(  f  )
 +\langle f , \Diamond _j \mathbf{S} ^{(\ell)} \rangle + \Pi _j(   \mathbf{S} ^{(\ell)})
	 \end{aligned} \end{equation}
where, the r.h.s.'s equal to the terms of \eqref{eq:ODE1} for $t=1$ for each of the two flows,        \begin{equation} \nonumber
\begin{aligned} & \mathbf{S} ^{(\ell)} =\mathbf{{S}}^{i,M_0 }_{ n- l-1,l}(a, \Pi _4, \Pi (f),z, f ).
\end{aligned} \end{equation}
   Hence $\Pi _j(   \mathbf{S} ^{(\ell)}) =
 \resto ^{i,2M_0 }_{ n- l-2,l}$,  and this term can be absorbed in the r.h.s. of \eqref{eq:ODE1tris}.

\noindent Next, observe that $\mathbf{S} ^{(\ell)}$ is the integral $\int _0^1  \mathcal{D} ^{(\ell)}dt$
of the terms $\mathcal{D} ^{(\ell)}$   of Lemma \ref{lem:ODEbis}. Formula
\eqref{eq:ODEbis1} implies
\begin{equation*}
     \mathbf{S} ^{(1)} -\mathbf{S} ^{(2)}= \mathbf{S}^{0,M_0+1}_{ n- l-2,l}(a,\Pi _4,\Pi(f),z,f),
\end{equation*}
as can be seen by elementary computations, and this  in turn implies
\begin{equation*}
    \langle r ,\Diamond _j \( \mathbf{S} ^{(1)} -\mathbf{S} ^{(2)} \) \rangle = \resto ^{0,M_0+2}_{ n- l-3,l}(a,\Pi _4,\Pi(f),z,f).
\end{equation*}

 \qed

We consider  $f\in   X_c \cap \Sigma_{ \mathbf{N}_0}$ for $\mathbf{N}_0 $ a large number.
We can pick $\mathbf{N}_0 >2\mathbf{N}+2$ where $ \mathbf{N}$ is defined in  \eqref{eq:h8}.
Notice that \eqref{eq:SystE2}  preserves this space. We have the following, which is  proved as in \cite{Cu0},
and which we discuss in Appendix \ref{sec:pullb}.

\begin{lemma} \label{lem:ODE1} Consider
 $\mathfrak{F} =\mathfrak{F}^1 \circ \cdots \circ \mathfrak{F}^L$
 with $ \mathfrak{F}^j= \mathfrak{F}^j_{t=1}$ transformations as in Lemma~\ref{lem:ODE}
 on the manifold  $\mathscr{M}_{ 1}^{6}(p^0 )$.
 Suppose that for any $ \mathfrak{F}^j$ the $M_0 $ in Lemma~\ref{lem:ODE}
 equals $m_j$,
 where $1= m_1\le...\le m_L$ with the constant $i$ in
Lemma~\ref{lem:ODE} (ii)   equal to $1$
when $m_j=1$.
Fix   $M, k$ with $n_1\gg k\ge \mathbf{N}_0$ ($n_1$ picked in Lemma \ref{lem:modulation}).
Then there is a $ n=n(L,M,k)$ such that
if the  assumptions of Lemma~\ref{lem:ODE}
apply to each of operators $ \mathfrak{F}^j$  for  $(M,n)$,
  there exist
$ {{{\psi}}}(p_4,\varrho)\in C^\infty$ with $ {{\psi}}((p_4,\varrho)=O(|\varrho|^2)$
 and a small $\varepsilon >0$
 such that in $\U^{s}_{\varepsilon,k}$  for $s\ge n-(M+1)$      we have the expansion
\begin{equation} \label{eq:ExpH11} \begin{aligned} & {K} \circ \mathfrak{F}= {
\psi} (p^0_4, \Pi (f) ) +H_2  +\textbf{R},
\end{aligned} \end{equation}
where $ \psi (p^0_4 ,\varrho) = O(\varrho^2)$ is $C^{m}$ and with what follows.\begin{itemize}
\item[(1)]
 We have
\begin{equation} \label{eq:ExpH2} H_2  = \sum_{\substack{ |\mu +\nu |=2 \ , \
\mathbf{e} \cdot (\mu -\nu)=0}}
 g_{\mu \nu} (p^0_4 ,\Pi (f)) z^\mu
\overline{z}^\nu +  {2}^{-1} \langle {\mathbbm{i}} \sigma_3 \mathcal{H}_{p^1
 } f, f\rangle.
\end{equation}

\item[(2)]
Denote
$
\varrho =\Pi (f).
$
There is the expansion
$\textbf{R} = \sum_{j=-1,...,3} {\textbf{R}_{ j }} + \resto^{1,2}_{k,m}(p^0_4,\varrho, f) $,
\begin{equation*} \begin{aligned} &{ {\mathbf{R}}_{ - 1 }}=
 \sum_{\substack{ |\mu +\nu |=2\ , \
\mathbf{e} \cdot (\mu -\nu)\neq 0 }} g_{\mu \nu } (p^0_4 ,\varrho)z^\mu
\overline{z}^\nu + \sum_{|\mu +\nu | = 1} z^\mu \overline{z}^\nu \langle
{\mathbbm{i}} \sigma_3 G_{\mu \nu }(p^0_4 ,\varrho),f\rangle ;\\&
|\resto^{1,2}_{k,m}(\Pi _4,\varrho, f)|\le C \| f\|_{\Sigma_{-k}}^2 (\| f\|
_{\Sigma_{-k}}+|\varrho | +|\Pi _4 -p^1_4 | +|z|) ;
\end{aligned}
\end{equation*} for $ \mathbf{N}$ as in (H8),
 \begin{equation*} \begin{aligned} & {\textbf{R}_0}= \sum_{|\mu
 +\nu |= 3,...,2 {\mathbf{N}}+2} z^\mu
\overline{z}^\nu g_{\mu \nu }(p^0_4 ,\varrho) ; \quad {\textbf{R}_1}= \im
\sum_{|\mu +\nu |= 2,...,2 {\mathbf{N}}+1 } z^\mu \overline{z}^\nu \langle {\mathbbm{i}} \sigma_3
G_{\mu \nu }(p^0_4 ,\varrho), f\rangle ;
\end{aligned}\nonumber \end{equation*}
  \begin{align} & {\textbf{R}_2}
= \sum_{|\mu +\nu |= 2\mathbf{N}+3} z^\mu \overline{z}^\nu g_{\mu \nu }(p^0_4 ,\varrho,
z,f) -\sum_{ |\mu +\nu |= 2\mathbf{N}+2} z^\mu
\overline{z}^\nu \langle {\mathbbm{i}} \sigma_3 G_{\mu \nu }(p^0_4 ,\varrho, z,f),
f\rangle ;\nonumber \\ & {\textbf{R}_3}= \sum
_{d=2,3,4} \langle B_{d } (p^0_4 ,\varrho, z,f), f^{ d} \rangle
 +\int_{\mathbb{R}^3}
B_5 (x, p^0_4 ,\varrho, z,f, f(x)) f^{ 5}(x)\,dx + E_P (f) \nonumber \\& \text{ with
$B_2(p^1,0,0,0)=0$.}  \label{eq:canc B2}\end{align}
Above, $f^d(x)$ schematically represents $d$-products of components of $f$.

\item[(3)] For $\delta_j\in\N_{0}^m$
the vectors defined
in terms of the Kronecker symbols by
$\delta_j:=(\delta_{1j},..., \delta_{mj}) $,
\begin{equation} \label{eq:ExpHcoeff1}
 \begin{aligned} & g_{\mu \nu } = \resto^{1,0}_{k,m}
\quad
\text{ for $|\mu +\nu
 | = 2$ \ for \ $(\mu, \nu)\neq (\delta_j, \delta_j)$, \ $1\le j\le m$;}
\\&
 g_{\delta_j \delta_j } =\mathbf{e}_j + \resto^{1,0}_{k,m},
\text{\ $1\le j\le m$};
\quad
 G_{\mu \nu } =\mathbf{S}^{1,0}_{k,m} \text{ for $|\mu +\nu | = 1$};
 \end{aligned} \end{equation}
$g_{\mu \nu }$
and $G_{\mu \nu }$
satisfy
symmetries analogous to \eqref{eq:symm}.
\item[(4)] All the other $g_{\mu \nu } $ are $\resto^{0,0}_{k,m} $ and all the other $G_{\mu \nu } $ are $\mathbf{S}^{0,0}_{k,m} $.

\item[(5)] $B_{d
}(p^0 , \varrho, z,f) \in C^{m } (\U_{-k},
\Sigma_k (\mathbb{R}^3, B (
 (\mathbb{R}^{4 })^{\otimes d},\mathbb{R}))) $ for $2\le d \le 4$ with $\U_{-k}\subset \R ^{8} \times \C^{\mathbf{n}}\times (X_c\cap \Sigma_{-k})$ an open neighborhood of $( p^1_4\ ,  \varrho , z, f)= (0,0,0,0)$.

\item[(6)] Let $ \zeta \in \mathbb{C}^{2 }$. Then for
 $B_5(\cdot, \varrho,z,f, \zeta)$ we have for fixed constants $C_l$ (the
 derivatives are not in the holomorphic sense)
\begin{equation} \label{H5power2}\begin{aligned} &\text{for $|l|\le m $,}
\quad
\| \nabla_{p^0, \varrho, z,f,\zeta }^lB_5(p^0_4 ,\varrho,z,f,\zeta
) \|_{\Sigma_k(\mathbb{R}^3, B (
 (\mathbb{C}^{2 })^{\otimes 5},\mathbb{R})} \le C_l.
 \end{aligned}  \end{equation}

\end{itemize}
\end{lemma}

\qed

\section{Darboux Theorem}
\label{sec:darboux2}

   Recall that we have introduced a model symplectict form $\Omega _0$ in $\mathscr{M}_{ 1}^{6}(p^0 )$ by formula \eqref{eq:modomega0}.
   Now we transform $\Omega  $ into $\Omega_0 $ by means of
  the Darboux Theorem, performed in a non-abstract way, to make sure that the coordinate transformation is as in Lemma~\ref{lem:ODE}.

 \begin{lemma}
 \label{lem:vectorfield01} For $n_1$ the constant in Lemma \ref{lem:modulation} and  $\varepsilon
_2>0$
	consider the set
\[
\U_{2}=\left\{(z,f)\in\C^{\mathbf{n}} \times (X_c \cap H^1) :\quad
\|f\|_{\Sigma_{-{n_1}}}\le \varepsilon
_2,\quad
| \Pi
(f)|\le \varepsilon_2,\quad |z|\le \varepsilon_2\right\}.
\]
 Then for $ \varepsilon
_2>0$ small enough
 there exists a unique vector field
		$ \mathcal{Y}^{t } $ in $\U_{2}$ such that $
 i_{\mathcal{Y}^{t }} (\Omega_0+ t(\Omega -\Omega_0)) =\Gamma _0-\Gamma
 $
 for $|t|<5$ with components, where $\Pi _4=p^0_4$,
\begin{equation*}\label{eq:quasilin10} \begin{aligned} &
 (\mathcal{Y}^{ t })_{z_j}= \mathcal{R}^{1,1}_{{n_1},\infty}(\Pi _4,\Pi (f),z,f) \ , \
 (\mathcal{Y}^{ t })_{f} = {\mathbbm{i}} \sigma_3 \mathcal{R}^{0,2}_{n,\infty}(\Pi _4,\Pi (f),z,f)\cdot \Diamond f +
 \mathbf{S}^{1,1}_{{n_1},\infty }(\Pi _4,\Pi (f),z,f).\end{aligned}
\end{equation*}
\end{lemma}
\proof The proof is essentially the same as that of \cite[Lemma 3.4]{Cu0}.
The first step is to consider a field $Z$ such that $i_{Z}\Omega_0 =\Gamma _0-\Gamma   $. We claim that
\begin{equation*}
 (Z)_{z}= \mathcal{R}^{1,1}_{\infty,\infty}(\Pi _4,\Pi (f),z,f) \,, \qquad (Z)_{f} = {\mathbbm{i}} \sigma_3 \mathcal{R}^{0,2}_{\infty,\infty}(\Pi _4,\Pi (f),z,f)\cdot \Diamond f +
 \mathbf{S}^{1,1}_{\infty,\infty }(\Pi _4,\Pi (f),z,f).
\end{equation*}
Schematically, the equation for $Z$ is of the form
 \begin{equation*}
 (Z)_{z}\,dz+\langle \left [ {\mathbbm{i}} \sigma_3 (Z)_{f} + \mathcal{R}^{0,0}_{\infty,\infty}\langle \Diamond f, (Z)_{f} \rangle \right ] \Diamond f, df \rangle = \mathcal{R}^{1,1}_{\infty,\infty}\,dz +\langle {\mathbbm{i}} \sigma_3 \mathcal{R}^{0,2}_{\infty,\infty}\cdot \Diamond f+\mathbf{S}^{1,1}_{\infty,\infty}, df\rangle.
 \end{equation*}
This immediately yields $(Z)_{z}= \mathcal{R}^{1,1}_{\infty,\infty}$. The equation for $ (Z)_{f}$ is of the form
\begin{equation}\label{eq:quasilin100}
(Z)_{f} + \mathcal{R}^{0,0}_{\infty,\infty}\langle \Diamond f, (Z)_{f}
 \rangle {\mathbbm{i}} \sigma_3 \Diamond f ={\mathbbm{i}} \sigma_3 \mathcal{R}^{0,2}_{\infty,\infty}\cdot \Diamond f+\mathbf{S}^{1,1}_{\infty,\infty},
\end{equation}
with a solution in the form
$\displaystyle (Z)_{f}= \sum_{i=0}^{\infty}(Z)_{f}^{(i)}$, with $(Z)_{f}^{(0)}={\mathbbm{i}} \sigma_3 \mathcal{R}^{0,2}_{\infty,\infty}\cdot \Diamond f+\mathbf{S}^{1,1}_{\infty , \infty} $
and
\begin{equation*}
 (Z)_{f}^{(i+1)}=\mathcal{R}^{0,0}_{\infty,\infty}\langle \Diamond f, (Z)_{f}^{(i)} \rangle {\mathbbm{i}} \sigma_3 \Diamond f = (\mathcal{R}^{0,0}_{\infty,\infty})^{i+1}\langle \Diamond f, {\mathbbm{i}} \sigma_3 \Diamond f \rangle^{i} \langle \Diamond f, (Z)_{f}^{(0)} \rangle {\mathbbm{i}} \sigma_3 \Diamond f,
\end{equation*}
where by direct computation $\langle \Diamond_jf, {\mathbbm{i}} \sigma_3 \Diamond_k f \rangle$ is a bounded bilinear form
in $X_c\cap L^2(\R^3, \C^4)$
for all $j,k$. This implies that the series defining $(Z)_{f}$ is convergent and that
$(Z)_{f}$ is as in \eqref{eq:quasilin100}.

\noindent The next step is to define an operator
$\mathcal{K} $    by $ i_{X}  (\Omega  -\Omega _0) =
 i_{\mathcal{K} X} \Omega_0$. We claim that
 \begin{equation}\label{eq:quasilin101} \begin{aligned} (\mathcal{K}X)_{z }&= \resto^{1,0}_{\infty,\infty} (X)_{z } + \langle \mathcal{R}^{0,2}_{\infty,\infty} \Diamond f  +  \mathbf{S}^{1,0}_{\infty,\infty}, (X)_{f } \rangle \\
 (\mathcal{K}X)_{f } &= {\mathbbm{i}} \sigma_3 \langle \mathbf{S}^{0,1}_{\infty,\infty}, (X)_{f } \rangle \Diamond f +
 \partial_f\left. \mathbf{S}^{1,1}_{\infty,\infty} \right |_{(\rho,z, f)=(\Pi (f),z, f)} (X)_{f }
 \\& (X)_{z }\mathcal{R}^{0,1}_{\infty,\infty}\Diamond f +(X)_{z } \mathbf{S}^{1,0}_{\infty,\infty} + \langle \Diamond f, (X)_{f }\rangle \mathbf{S}^{0,1}_{\infty,\infty}.
\end{aligned}
\end{equation}
 From  \eqref  {eq:modomega0}--\eqref  {eq:omedgamega0}
 we have
  schematically
 \begin{equation*} \begin{aligned} & \im (\mathcal{K}X)_{z }d {z} + \langle \left [{\mathbbm{i}} \sigma_3 (\mathcal{K}X)_{f } + \mathcal{R}^{0,0}_{\infty,\infty} \langle \Diamond f, (\mathcal{K}X)_{f } \rangle \Diamond f\right ], df \rangle \\& = \left (\resto^{1,0}_{\infty,\infty} (X)_{z } + \langle \resto^{0,1}_{\infty,\infty}\Diamond f+\mathbf{S}^{1,0}_{\infty,\infty}, (X)_{f } \rangle \right) \,dz  + \langle \left [ \partial_f\left. \mathbf{S}^{1,1}_{\infty,\infty} \right |_{(\rho,z, f)=(\Pi (f),z, f)} (X)_{f } \right . \\& \left .
+(X)_{z }( \resto^{0,1}_{\infty,\infty}\Diamond f+   \mathbf{S}^{1,0}_{\infty,\infty})+\langle \mathbf{S}^{0,1}_{\infty,\infty},(X)_{f }\rangle
\Diamond f+\langle \Diamond f,(X)_{f }\rangle \mathbf{S}^{0,1}_{\infty,\infty}\right ],df \rangle
\end{aligned}
\end{equation*}
which yields immediately the first equation in \eqref{eq:quasilin101}.
We have $\displaystyle (\mathcal{K}X)_{f } =\sum_{i=0}^{\infty}(\mathcal{K}^{(i)}X)_{f } $ with
\begin{equation*} \begin{aligned}
{\mathbbm{i}} \sigma_3 (\mathcal{K}^{(0)}X)_{f } =
 \partial_f\left. \mathbf{S}^{1,1}_{\infty,\infty} \right |_{(\rho,z, f)=(\Pi (f),z, f)} (X)_{f }
 +(X)_{z }( \resto^{0,1}_{\infty,\infty}\Diamond f+   \mathbf{S}^{1,0}_{\infty,\infty})\\   + \langle \mathbf{S}^{0,1}_{\infty,\infty}, (X)_{f } \rangle
 \Diamond f + \langle \Diamond f, (X)_{f }\rangle \mathbf{S}^{0,1}_{\infty,\infty}
\end{aligned}
\end{equation*}
and
\begin{equation*} \begin{aligned} &
(\mathcal{K}^{(i+1)}X)_{f } = \mathcal{R}^{0,0}_{\infty,\infty} \langle \Diamond f, (\mathcal{K}^{(i)}X)_{f } \rangle {\mathbbm{i}} \sigma_3 \Diamond f = (\mathcal{R}^{0,0}_{\infty,\infty})^{i+1}\langle \Diamond f, {\mathbbm{i}} \sigma_3 \Diamond f \rangle^{i} \langle \Diamond f, (\mathcal{K}^{(0)}X)_{f } \rangle {\mathbbm{i}} \sigma_3 \Diamond f.
\end{aligned}
\end{equation*}
Then the series defining $(\mathcal{K} X)_{f }$ converges and we get in particular the second equation in \eqref{eq:quasilin101}.
 Now the equation defining $\mathcal{Y}^{t}$ is equivalent to $(1+t \mathcal{K}) \mathcal{Y}^t= Z$. So we have \begin{equation*} \begin{aligned} & (\mathcal{Y}^t)_{z } +t \resto^{1,0}_{\infty,\infty} (\mathcal{Y}^t)_{z } + t \langle \mathcal{R}^{0,2}_{\infty,\infty} \Diamond f  +\mathbf{S}^{1,0}_{\infty,\infty}, (\mathcal{Y}^t)_{f } \rangle = \mathcal{R}^{1,1}_{\infty,\infty} \\& (\mathcal{Y}^t)_{f } +{\mathbbm{i}}t \sigma_3 \langle \mathbf{S}^{0,1}_{\infty,\infty}, (\mathcal{Y}^t)_{f } \rangle \Diamond f +
 t\partial_f\left. \mathbf{S}^{1,1}_{\infty,\infty} \right |_{(\rho,z, f)=(\Pi (f),z, f)} (\mathcal{Y}^t)_{f }
 \\&+t(\mathcal{Y}^t)_{z }(  \mathcal{R}^{0,2}_{\infty,\infty} \Diamond f  + \mathbf{S}^{1,0}_{\infty,\infty} )+ t \langle \Diamond f, (\mathcal{Y}^t)_{f }\rangle \mathbf{S}^{0,1}_{\infty,\infty} = {\mathbbm{i}} \sigma_3 \mathcal{R}^{0,2}_{\infty,\infty}\cdot \Diamond f +
 \mathbf{S}^{1,1}_{\infty,\infty }.
\end{aligned}
\end{equation*}
Solving this we get the desired formulas for $(\mathcal{Y}^{ t })_{z_j}$ and $(\mathcal{Y}^{ t })_{f}$.
\qed

We can apply Lemma~\ref{lem:ODE}
to the flow $ {\mathfrak{F}}_t:\;(z,f) \mapsto (z^t, f^t)$ generated by $ \mathcal{Y}^t$.
In terms of the decomposition \eqref{eq:decomp2} of $r$ formula \eqref{eq:ODE1}
becomes for $n= {n_1}$
\begin{align}
 &z^t= z+\mathcal{R}^{1,1}_{{n_1}- l-1, l}(t,\Pi _4
 ,\Pi (f), z,  f), \label{eq:tildef}\\&
 f^t= e^{{\mathbbm{i}} \sum _{j=1}^{4}\sigma_3\resto^{0,2 }_{ {n_1}- l-1, l}(t,\Pi _4,\Pi (f), z,  f)  \Diamond _j }
 T(e^{\sum _{a=1}^{3}\resto^{0,2 }_{ {n_1}- l-1, l}(t,\Pi _4,\Pi (f), z,  f) {\mathbbm{i}}\sigma _a    })\nonumber\\& \times (f+
\mathbf{{S}}^{1,1 }_{ {n_1}- l-1,l}
 (t,\Pi _4,\Pi (f), z,  f)).\nonumber\end{align}
 Classically the Darboux Theorem follows by $ i_{\mathcal{Y}
^{t }} \Omega_t  =\Gamma _0 -\Gamma   $,  where $ \Omega _t:= \Omega _0 +t(\Omega -\Omega _0)$, and by
\begin{equation}\label{eq:fdarboux} \begin{aligned} & \partial_t
 (\mathfrak{F}_{ t }^*\Omega_t ) = \mathfrak{F}_{ t
 }^* (L_{\mathcal{ Y}^{ t }}\Omega_t  +\partial_t\Omega_t ) =
 \mathfrak{F}_{ t }^*
 (d i_{\mathcal{Y}^{ t }}\Omega '_t +d(\Gamma  -\Gamma_0 ))
=0 \end{aligned} \end{equation}
 with $L_X$ the Lie derivative, whose definition is not needed here.
Since this $\mathfrak{F}_{ t }$ is not a differentiable flow on any given manifold, \eqref{eq:fdarboux} is
 formal. Still, \cite[Sect. 3.3 and Sect. 7]{Cu0} (i.e. a regularization and a limit argument for $ \mathfrak{F}_{t }$) yield the following, which we state without proof.

\begin{lemma} \label{lem:darboux} Consider \eqref{eq:ODE} defined by the field
${\mathcal X}^{t }$ and  indexes and notation of Lemma~\ref{lem:ODE} (in particular $M_0=1$ and $i=1$; $n$ and $M$ can be
arbitrary as long as we fix ${n_1}$ large enough). Consider $l$, $ {s}'$,$ {s}$ and $k$ as  in \eqref{eq:index1}.
Then for $\mathfrak{F}_{1 } \in C^l(\U_{\varepsilon
_2,k}^{\mathbf{s}'}, \U_{\varepsilon
_1,k}^{\mathbf{s}})$
derived from \eqref{eq:tildef},
 we have $ \mathfrak{F}_{1 }^* \Omega  =\Omega_0 $.  \qed
\end{lemma}

We now turn to the analysis of the hamiltonian vector fields in the new coordinate system.
For a function $F$ let us decompose $X_F$ according to the spectral decomposition \eqref{eq:decomp2}: for $(X_F)_{f} \in X_c $,
\begin{equation}\label{eq:XFext}
\begin{aligned}
 X_F &= \sum_{j=1,...,\mathbf{n}}(X_F)_{z_j} \xi_j (x) +
\sum_{j=1,...,\mathbf{n}}(X_F)_{\overline{z}_j} {\xi }_j^*(x)
+ (X_F)_{f}. \end{aligned}
\end{equation}
By \eqref{eq:modomega0} and by $i_{X_F}\Omega_0=dF$ we have, schematically (recall also that here and below, $\Pi _4=p^0_4$,
\begin{equation}
\begin{aligned}& - \im (X_F)_{z_l} d \overline{z}_l+ (X_F)_{\overline{z}_l} \,dz_l
+ \langle \left [ {\mathbbm{i}} \sigma_3 (X_F)_{f} +\resto^{0,0}_{\infty,\infty}(\Pi _4,\Pi (f))\langle \Diamond f, (X_F)_{f} \rangle \Diamond f \right ], df\rangle
 \\& = \partial_{z_l} F d {z}_l +\partial_{\overline{z}_l} F d {\overline{z}}_l + \langle \nabla_{f} F, df\rangle. \end{aligned}\nonumber
\end{equation}
and so, schematically,
\begin{equation}
\begin{aligned}& (X_F)_{z_l} =\im \partial_{\overline{z}_l} F \ , \quad (X_F)_{\overline{z}_l} =-\im \partial_{z_l} F\\&
 (X_F)_{f} +\resto^{0,0}_{\infty,\infty}(\Pi _4,\Pi (f))\langle \Diamond f, (X_F)_{f} \rangle P_c {\mathbbm{i}} \sigma_3 \Diamond f = -{\mathbbm{i}} \sigma_3\nabla_{f} F. \end{aligned}\nonumber
\end{equation}
We set
\begin{align}\label{eq:X_Fsplit}
 &X_F =X_F^{(0)}+X_F^{(1)} \text{ with}\\&  (X_F^{(0)})_{z_l} =\im \partial_{\overline{z}_l} F \quad, \quad (X_F^{(0)})_{\overline{z}_l} =-\im \partial_{z_l} F\quad, \quad
 (X_F^{(0)})_{f} = -{\mathbbm{i}} \sigma_3\nabla_{f} F \label{eq:X_F0}
\end{align}
and where the remainder is of the form $(X_F^{(1)})_{z_l} =(X_F^{(1)})_{\overline{z}_l} =0$
\begin{equation}\label{eq:X_F1}
\begin{aligned}&
 (X_F^{(1)})_{f} =
 %\resto^{0,0}_{\infty,\infty}(\Pi  _4,\Pi (f))\langle \Diamond f, (X_F^{(0)})_{f} \rangle P_c {\mathbbm{i}} \sigma_3 \Diamond f =
  \resto^{0,0}_{\infty,\infty}(\Pi _4,\Pi (f))\langle \Diamond f, {\mathbbm{i}} \sigma_3\nabla_{f} F \rangle P_c {\mathbbm{i}} \sigma_3 \Diamond f. \end{aligned}
\end{equation}
Indeed, $(X_F^{(1)})_{f}$ has to satisfy an equation of the form
\begin{equation}
\begin{aligned}&
 (X_F^{(1)})_{f} +\resto^{0,0}_{\infty,\infty}(\Pi _4,\Pi (f))\langle \Diamond f, (X_F)_{f}^{(1)} \rangle P_c {\mathbbm{i}} \sigma_3 \Diamond f = \resto^{0,0}_{\infty,\infty}(\Pi _4,\Pi (f))\langle \Diamond f, {\mathbbm{i}} \sigma_3\nabla_{f} F \rangle P_c {\mathbbm{i}} \sigma_3 \Diamond f. \end{aligned}\nonumber
\end{equation}
This can be solved like in the proof of Lemma~\ref{lem:vectorfield01} by writing
$\displaystyle (X_F)_{f}^{(1)} = \sum_{i=0}^{\infty} X_i$ with
\begin{equation*}\begin{aligned} X_0&= \resto^{0,0}_{\infty,\infty}(\Pi _4,\Pi (f))\langle \Diamond f, {\mathbbm{i}} \sigma_3\nabla_{f} F \rangle P_c {\mathbbm{i}} \sigma_3 \Diamond f \text{ and}\\
 X_{i+1}&=\resto^{0,0}_{\infty,\infty}(\Pi _4,\Pi (f))\langle \Diamond f, X_{i} \rangle P_c {\mathbbm{i}} \sigma_3 \Diamond f \\&= (\mathcal{R}^{0,0}_{\infty,\infty})^{i+1}\langle \Diamond f, {\mathbbm{i}} \sigma_3 \Diamond f \rangle^{i} \langle \Diamond f, {\mathbbm{i}} \sigma_3\nabla_{f} F \rangle P_c {\mathbbm{i}} \sigma_3 \Diamond f\end{aligned}
\end{equation*}
which yields \eqref{eq:X_F1}.
For two functions $F$ and $G$
we have the Poisson brackets \begin{equation} \label{eq:poiss}\begin{aligned} & \{ F,G \} :=dF(X_G) = \partial_{z_{l} }F (X_G)_{z_{l} } + \partial_{\overline{z}_{l} }F (X_G)_{\overline{z}_{l} } +\langle \nabla_fF, (X_G)_{f } \rangle =
 \{ F,G \}_{(0)} + \{ F,G \}_{(1)},
\end{aligned}\end{equation}
where $ \{ F,G \}_{(i)}:=dF(X_G^{(i)}) $ and
where
 \begin{equation} \label{eq:poiss0}\begin{aligned} &
 \{ F,G \}_{(0)}= \im (\partial_{z_l } F \partial_{z_l } G - \partial_{\overline{z}_l } F \partial_{\overline{z}_l } G) -\langle \nabla_fF, {\mathbbm{i}} \sigma_3\nabla_f G \rangle
\end{aligned}\end{equation}
and, schematically,
\begin{equation} \label{eq:poiss1}\begin{aligned}
 \{ F,G \}_{(1)}& = \resto^{0,0}_{\infty,\infty}(\Pi _4,\Pi (f)) \langle \nabla_fF,\Diamond f \rangle \langle \Diamond f, {\mathbbm{i}} \sigma_3\nabla_{f} G \rangle.
\end{aligned}\end{equation}
Compared to   \cite{Cu0}, where the Poisson bracket equals \eqref{eq:poiss0}, here   we have an additional term contributed by
\eqref{eq:poiss1}, which however is of higher order and   harmless, as we will see later.

\section{ Birkhoff normal forms}
\label{sec:Normal form}
 We will reduce now to \cite[Sect. 6]{Cu0}.
 We set, for the $\mathbf{e}_{j}$'s in (H6), see Section~\ref{sec:speccoo},
 \begin{equation*}
 \mathbf{e} :=(\mathbf{e}_{1},..., \mathbf{e}_{\mathbf{n}}).
 \end{equation*}
    In the sequel, $\Pi _4=p^0_4.$

\begin{definition}

\label{def:normal form} A function $Z(\varrho, z,f)$ is in normal form if
$ Z=Z_0+Z_1$,
where $Z_0$ and $Z_1$ are finite sums of the following type:
\begin{equation}
\label{e.12a}Z_1= \sum_{\mathbf{e} \cdot(\nu-\mu)\in \sigma_e(\mathcal{H}_{p^1})}
z^\mu \overline{z}^\nu \langle {\mathbbm{i}} \sigma_3 G_{\mu \nu}(p ^0_4,\varrho),f\rangle
\end{equation}
with $G_{\mu \nu}(x, p _4,\varrho)\in C^{m} (U,\Sigma_{k }(\R^3, \C^{4}))$ for fixed $k,m\in\N$ and $U\subseteq \R^{8} $ an open neighborhood of $(p^0_4,0) $,
\begin{equation}
\label{e.12c}Z_0= \sum_{ \mathbf{e} \cdot(\mu-\nu)=0} g_{\mu \nu}
(p ^0_4,\varrho)z^\mu \overline{z}^\nu,
\end{equation}
with $g_{\mu \nu} (p_4,\varrho)\in C^{m} (U,
\mathbb{C})$.
We assume furthermore the symmetries $ \overline{{g}}_{\mu \nu} =g_{\nu\mu }$ and $ \overline{{G}}_{\mu \nu} =G_{\nu \mu }$.
 \end{definition}

\begin{lemma}
 \label{lem:chi} For $i\in \{ 0,1\}$ fixed and $n,M\in\N$
sufficiently large and for $m\le M-1$ let
	 \begin{equation*} \label{eq:chi1}\chi =\sum_{|\mu +\nu |=M_0 +1}
c_{\mu\nu} (p ^0_4,\Pi (f)) z^{\mu} \overline{z}^{\nu} +\im \sum_{|\mu +\nu
|=M_0 } z^{\mu} \overline{z}^{\nu}
 \langle {\mathbbm{i}} \sigma_3 C_{\mu \nu
}(p ^0_4,\Pi (f)), f \rangle,
\end{equation*} with $ c_{\mu\nu}(p^0,\varrho)= \resto^{i,0}
_{n,M}(p^0,\varrho)$ and $ C_{\mu\nu}(p^0,\varrho)= \mathbf{S}^{i,0}_{n,M}(p^0,\varrho)$ and with \begin{equation} \label{eq:symm}
 \overline{{c}}_{\mu\nu} = {c}_{\nu\mu}, \quad \overline{{C}}_{\mu \nu }
=-C_{\nu\mu}\end{equation}
(so that $\chi$ is real-valued for
$f= \overline{{f}} $). Then we have what follows.
\begin{itemize}
\item[(1)] For $\phi^{t }$ the flow of $X_\chi $, see
 Lemma~\ref{lem:ODE}, and $(z^t,f^t)= (z,f)\circ \phi^{t
 }$,
\begin{align}   z^t =& z +
\resto^{0,M_0}_{n-m-1,m-1}(t, ,\Pi _4, \Pi (f), z, f),\nonumber \\
   f^t =&e^{{\mathbbm{i}} \sigma_3\sum  _{j=1}^{4}\resto^{0,M_0+1}_{n-m-1,m-1}(t, \Pi _4, \Pi (f), z, f)  \Diamond  _j}
T(e^{ \sum  _{i=1}^{3}\resto^{0,M_0+1}_{n-m-1,m-1}(t, \Pi _4, \Pi (f), z, f)
 {\mathbbm{i}} \sigma  _i})\label{eq:quasilin51}\\&
  \quad  \quad\quad\quad\quad  \times  (f+ \mathbf{{S}}^{0,M_0}_{n-m-1,m-1}(t, \Pi _4, \Pi (f), z, f))\nonumber.
\end{align}

\item[(2)] For $n- m-1 \ge s' \ge s+m-1 \ge m-1 $ and $k\in
\Z\cap [0,n- m-1 ]$
and
 for $ \varepsilon_1> \varepsilon_2> 0$ sufficiently small,
 $\phi :=\phi^{1 }
 \in C^{m-1} (\U_{\varepsilon_2,k}^{s'}, \U_{\varepsilon_1,k}^{s}
)
 $
 satisfies
 $\phi^{ \ast }\Omega_0=\Omega_0$.
\end{itemize}
\end{lemma}
\proof This result is a simple corollary of Lemma~\ref{lem:ODE}.
For the proof that $\phi^{ \ast }\Omega_0=\Omega_0$, which is obvious
in the standard setups, see the comments in \cite[Lemma 5.3]{Cu0}.
\qed

Then we have the following result on Birkhoff normal forms.

\begin{proposition}
\label{th:main} For any integer $2\le \ell \le 2\mathbf{N}+2$
there are transformations
$\mathfrak{F}^{(\ell)} = \mathfrak{F}_1 \circ  \phi_2\circ...\circ \phi_\ell $,
with $\mathfrak{F}_1$ the transformation  in \eqref{eq:tildef} and with
the $\phi_j$'s
 like in Lemma~\ref{lem:chi}, such that
 the conclusions of Lemma~\ref{lem:ODE1} hold,
 that is such that we have the following expansion, for $\Pi _4=p^0_4$,
\begin{equation*}
 H^{(\ell)}:=K\circ \mathfrak{F}^{(\ell)} ={ \psi} (p ^0_4,\Pi  (f) )+ H_2  + \resto^{1,2}_{k,m}(\Pi _4,\Pi (f), f) + \sum_{j=-1,...,3}  \textbf{R}_{j}^{(\ell)},
\end{equation*}
with $H_2'$ defined in  \eqref{eq:ExpH2} and with the following additional properties:
\begin{itemize}
\item[(i)] $\textbf{R}_{-1}^{(\ell)} =0$;

\item[(ii)] all the nonzero terms in $\textbf{R}_0^{(\ell)} $ with $|\mu +\nu |\le \ell $ are in normal form,
that is $\mathbf{e} \cdot (\mu -\nu)=0$;

\item[(iii)] all the nonzero terms in $\textbf{R}_1^{(\ell)} $ with $|\mu +\nu |\le \ell -1$ are in normal form,
that is $\mathbf{e}\cdot (\mu -\nu)\in \sigma_e(\mathcal{H}_{p^0})$.
\end{itemize}
\end{proposition}
\proof The proof of the analogue of Proposition~\ref{th:main}   in \cite{Cu0} involves the simpler    symplectic form
\begin{align*}
 \Omega_0^{(0)}:= - \im \sum_{l=1,...,\mathbf{n}} \,dz_l\wedge d \overline{z}_l
+ \langle {\mathbbm{i}} \sigma_3 df, df\rangle.
\end{align*}
In \eqref{eq:ExpH11},
we replace
 $\Pi (f)$ with $\varrho$; then $\mathbf{h}=H^{(\ell)}(p^0,\varrho, z,f)$ is
	$C^{2\mathbf{N}+2}$
	near $(0,0,0) $ in $(\varrho, z,f)\in \R^7\times \C \times (X_c \cap \Sigma_{ k}) $
and the statement of Proposition~\ref{th:main}
 is about the fact that some of the following derivatives vanish: \begin{align}
 \label{eq:derivmain1} & g_{\mu \nu}(p^0,\varrho)
=\frac{1}{\mu !\nu ! } \partial_z^\mu \partial_{\overline{z}}^\nu
 \mathbf{h}\at{(\varrho, z,f)= (\varrho,0,0)}, \quad |\mu +\nu |\le 2\mathbf{N}+2,\\& \label{eq:derivmain2}{\mathbbm{i}} \sigma_3 G_{\mu \nu}(p^0,\varrho)
=\frac{1}{\mu !\nu ! } \partial_z^\mu \partial_{\overline{z}}^\nu \nabla_f
 \mathbf{h}\at{(\varrho, z,f)= (\varrho,0,0)}, \quad |\mu +\nu |\le 2\mathbf{N}+1.
\end{align}
The proof is iterative and consists in assuming   the statement correct for a given $\ell$ and proving it for $\ell +1$,
by picking an unknown $\chi$ as in \eqref{eq:chi1} such that $H^{(\ell)}
\circ\phi $
satisfies the conclusions for $\ell +1$, where $\phi = \phi ^1 $, for $\phi ^t$ the flow for the Hamiltonian vector field
of $\chi$.

\noindent Now, let us pick $\chi$ provided by \cite[Theorem 6.4]{Cu0} when we use the symplectic form
$ \Omega_0^{(0)}$. We will show that this same $\chi$  works here.

\noindent Let $ {\phi}^{(0)} $ be the   $t=1$   flow generated by  $X_\chi^{(0)}$. Notice that $ {\phi}^{(0)} $ is a
symplectomorphism for $ \Omega_0^{(0)}$. Set
\begin{equation}\label{eq:mainblock-1}
    \widehat{H}^{(\ell)}= \psi  (p^0,\Pi (f) ) +H_2  + \sum_{j=-1,0,1}  \mathbf{{R}}_{j}^{(\ell)}.
\end{equation}
Noticing  that  here $\psi  (p^0,\Pi  (f) )$
yields 0 because it is ${ \psi} (p^0,\varrho  )$ with $\varrho$  an auxiliary independent variable,
\begin{equation}\label{eq:mainblock-}
   \begin{aligned}
&   \partial_z^\mu \partial_{\overline{z}}^\nu
 {H}^{(\ell)}\at{(\varrho, z,f)= (\varrho,0,0)}= \partial_z^\mu \partial_{\overline{z}}^\nu
 \widehat{{H}}^{(\ell)}\at{(\varrho, z,f)= (\varrho,0,0)}, \quad 2\le |\mu +\nu |\le 2\mathbf{N}+2,\\& \partial_z^\mu \partial_{\overline{z}}^\nu \nabla_f
  {{H}}^{(\ell)}\at{(\varrho, z,f)= (\varrho,0,0)}=\partial_z^\mu \partial_{\overline{z}}^\nu \nabla_f
  \widehat{{H}}^{(\ell)}\at{(\varrho, z,f)= (\varrho,0,0)}, \quad 1\le |\mu +\nu |\le 2\mathbf{N}+1
\end{aligned}
\end{equation}
since all the other terms of ${H}^{(\ell)}$  not contained in $\widehat{{H}}^{(\ell)}$  are higher order
in some of the variables, for example order 2 or higher in $f$. As we pointed out,
 $\psi  (p^0,\Pi (f))$ contributes nothing to \eqref{eq:mainblock-}.
 The same is true of the term $\frac{1 }{2} \langle {\mathbbm{i}} \sigma_3 \mathcal{H}_{p^1
 } f, f\rangle$  inside  $H_2 '$, see \eqref {eq:ExpH2} (however, the pullbacks of these terms are significant in the
 formulas below). So the only contributors of  \eqref{eq:mainblock-1} to \eqref{eq:mainblock-}
 are very  regular
  functions in $(\varrho , z, f)$, where  $\varrho =\Pi (f)$
  is as before treated as auxiliary variable and $f\in  (X_c
  \cap \Sigma_{- k}) $.  This yields  the useful result that while the l.h.s.'s in
\eqref{eq:mainblock-} require $f$ quite regular, for example $f\in \Sigma_{  k} $  for a sufficiently
 large   $k$,
 the r.h.s.'s are defined for $f\in \Sigma_{  -k} $ for a large preassigned $k$. This because the only term
 in  $ \widehat{H}^{(\ell)}(p^0,\varrho , z,f  )$ that requires some regularity in $f$ to make sense, that is the
 $\frac{1 }{2} \langle {\mathbbm{i}} \sigma_3 \mathcal{H}_{p^1
 } f, f\rangle$  hidden inside $ H'_2$, see \eqref{eq:ExpH2}, does not contribute to \eqref{eq:mainblock-}.

\noindent Furthermore,    by Lemma   \ref{lem:chi} we have
\begin{equation}\label{eq:mainblock}
   \begin{aligned}
&   \partial_z^\mu \partial_{\overline{z}}^\nu
 {H}^{(\ell)}\circ {\phi}^{(0)}  \at{(\varrho, z,f)= (\varrho,0,0)}= \partial_z^\mu \partial_{\overline{z}}^\nu
 \widehat{{H}}^{(\ell)}\circ {\phi}^{(0)}\at{(\varrho, z,f)= (\varrho,0,0)}, \quad 2\le |\mu +\nu |\le 2\mathbf{N} +1,\\& \partial_z^\mu \partial_{\overline{z}}^\nu \nabla_f
   {{H}}^{(\ell)}\circ {\phi}^{(0)}\at{(\varrho, z,f)= (\varrho,0,0)}= \partial_{\overline{z}}^\nu \nabla_f
  \widehat{{H}}^{(\ell)}\circ {\phi}^{(0)}\at{(\varrho, z,f)= (\varrho,0,0)}, \quad 1\le |\mu +\nu |\le 2\mathbf{N}
\end{aligned}
\end{equation}
since the pull backs of the terms of ${H}^{(\ell)}$  not contained in $\widehat{{H}}^{(\ell)}$
have zero derivatives  because are higher order either in $z$ or in $f$, as can be seen considering that
${\phi}^{(0)}$ acts like \eqref{eq:quasilin51}  for $M_0=\ell$.  Since $\phi$ too has this structure,
  \eqref{eq:mainblock} is true also with ${\phi}^{(0)}$ replaced by $\phi$.
  Set now
\begin{equation}\label{eq:pullcheck1}
    \widehat{{H}}^{(\ell)}\circ {\phi}
= \psi (p^0,\varrho  ) +  F \text{ with } F:=  \widehat{{H}}^{(\ell)}\circ {\phi}
-\psi (p^0,\varrho ) .
\end{equation}
We have $d F  \at{(\varrho, z,f)= (\varrho,0,0)}=0$, since
by Lemma \ref{lem:ODE1} we see that
is at least quadratic in $(z,f)$ .
Lemma \ref{lem:ODEbis} is telling us that $\phi ^{-1} \circ  {\phi}^{(0)}$    is the identity
up to a  zero of order $\ell+1$
at $(z,f)=(0,0)$    in  $\C ^{\mathbf{n}}\times  (X_c \cap \Sigma_{- k}) $.
Then by an elementary application of the chain rule
\begin{equation*}
   \begin{aligned}
&   \partial_z^\mu \partial_{\overline{z}}^\nu
 F     \at{(\varrho, z,f)= (\varrho,0,0)}= \partial_z^\mu \partial_{\overline{z}}^\nu
 F\circ \phi ^{-1} \circ  {\phi}^{(0)} \at{(\varrho, z,f)= (\varrho,0,0)}, \quad 2\le |\mu +\nu |\le \ell +1,
 \\& \partial_z^\mu \partial_{\overline{z}}^\nu \nabla_f
   F  \at{(\varrho, z,f)= (\varrho,0,0)}= \partial_{\overline{z}}^\nu \nabla_f
 F\circ \phi ^{-1} \circ  {\phi}^{(0)}\at{(\varrho, z,f)= (\varrho,0,0)}, \quad 1\le |\mu +\nu |\le \ell .
\end{aligned}
\end{equation*}
On the other hand, by Lemma \ref{lem:ODEtris}  we have that $\psi (p^0,\varrho)$
and $\psi (p^0,\varrho)\circ \phi ^{-1} \circ  {\phi}^{(0)}$  differ by a zero of order $\ell +2$ in
$(\varrho,0,0)$. Summing up, we conclude
\begin{equation*}
   \begin{aligned}
&   \partial_z^\mu \partial_{\overline{z}}^\nu
 \widehat{{H}}^{(\ell)}\circ {\phi}   \at{(\varrho, z,f)= (\varrho,0,0)}= \partial_z^\mu \partial_{\overline{z}}^\nu
 \widehat{{H}}^{(\ell)}\circ {\phi}^{(0)}\at{(\varrho, z,f)= (\varrho,0,0)}, \quad 2\le |\mu +\nu |\le \ell +1,\\& \partial_z^\mu \partial_{\overline{z}}^\nu \nabla_f
   \widehat{{H}}^{(\ell)}\circ {\phi} \at{(\varrho, z,f)= (\varrho,0,0)}= \partial_{\overline{z}}^\nu \nabla_f
  \widehat{{H}}^{(\ell)}\circ {\phi}^{(0)}\at{(\varrho, z,f)= (\varrho,0,0)}, \quad 1\le |\mu +\nu |\le \ell .
\end{aligned}
\end{equation*}
   Hence we have shown that \cite[Theorem 6.4]{Cu0} implies Proposition~\ref{th:main}. \qed

\section{Formulation of the system}
\label{sec:system}
So we consider the Hamiltonian $ {H}:= H^{(2\mathbf{N}+1)}$ and the reduced system \begin{equation} \label{eq:SystK1} \begin{aligned} &
 \dot z= \{ z, {H}
\} \,, \quad \dot f= \{ f, {H}
\}. \end{aligned}
\end{equation}
Recall that
\begin{equation}
 \label{eq:partham} {H}=\psi (p ^0_4,\Pi  (f) ) + H _2+Z_0+Z_1+\resto,
\end{equation}
with $H'_2$ like \eqref{eq:ExpH2}, $Z_0$ like \eqref{e.12c}, $Z_1$ like \eqref{e.12a},
and $\resto =\sum_{j=2,3} \mathbf{R}_j +\resto^{1,2}_{k,m}(\Pi _4,\Pi (f), f).$

\noindent We  recall that,
in the context of Strichartz estimates,
a pair $(p,q)$ is called \emph{admissible} if
\begin{equation}\label{admissiblepair}  2/p+3/q= 3/2,
\qquad 2\le q\le 6,
\qquad p\ge 2.
\end{equation}

\begin{theorem}\label{thm:mainbounds} For the constants $0<\epsilon <\epsilon_0 $ of Theorem~\ref{theorem-1.1}, there is a fixed
$C >0$ s.t.
%% for $I= [0,\infty)$ we have:
\begin{align}
& \| f \|_{L^p_t(\R_{+},W^{ 1,q}_x)}\le
 C \epsilon \text{ for all admissible pairs $(p,q)$,}
 \label{Strichartzradiation}
\\& \| z^\mu \|_{L^2_t(\R_{+})}\le
 C \epsilon \text{ for all multi-indexes $\mu$
 with $ \mathbf{{e}}\cdot \mu >\omega_1 $,} \label{L^2discrete}\\& \| z \|
_{W^{1,\infty}_t (\R_{+})}\le
 C \epsilon.\label{L^inftydiscrete}
\end{align}
Furthermore, we have $
\lim_{t\to +\infty} z(t)= 0. $
\end{theorem}

By standard arguments that we skip, such as a  simpler version of \cite[Sect. 7]{CM1}, Theorem~\ref{thm:mainbounds}
is a consequence of the following continuity argument.

\begin{proposition}\label{prop:mainbounds} For the constants $0<\epsilon <\epsilon_0 $ of Theorem~\ref{theorem-1.1}, there exists a constant $ \kappa >0$ s.t.
for any $C_0>\kappa$ there is $ \epsilon_0 >0 $ s.t. if the inequalities \eqref{Strichartzradiation}--\eqref{L^inftydiscrete}
hold for $I=[0,T]$ for some $T>0$
and for $C=C_0$,
then in fact the inequalities \eqref{Strichartzradiation}--\eqref{L^inftydiscrete} hold
for $I=[0,T]$ for $C=C_0/2$.
\end{proposition}

We now discuss the proof of Proposition~\ref{prop:mainbounds}, which is similar to
the proof for the scalar NLS, see for example \cite{CM1} or \cite{Cu3}. We have, see \eqref{eq:X_Fsplit}, $\dot f=(X_H^ {(0)})_f + (X_H^ {(1)})_f.$

\noindent In \cite{Cu3}, the equation was $\dot f= (X_H^ {(0)})_f$.
Given multi-indexes $\Theta ', \Theta \in \N_{0}^{m} $ we write
$\Theta ' < \Theta $ if $\Theta ' \neq \Theta $ and $\Theta '_l \le \Theta_l$,
$1\le l\le m$.
We now introduce
\begin{equation}\label{def-m-0}
\mathbf{M}_0
=\left\{\mu\in\N_{0}^{\mathbf{n}} :\; |\mathbf{e} \cdot\mu|> \omega^1 \,,
\quad
|\mu|\leq 2\mathbf{N} +2,
\quad
|\mathbf{e} \cdot\mu '|< \omega^1 \text{ if $ \mu ' < \mu $} \right\},
\end{equation}
\begin{equation} \label{eq:defM}\begin{aligned} \mathbf{M}
=\{ (\mu, \nu) \in
\N_{0}^{2\mathbf{n}}& :\; |\mathbf{e} \cdot(\mu-\nu)|> \omega^1 \,,
\quad
|\mu+\nu|\leq 2\mathbf{N} +2 \text{ and }\\& |\mathbf{e} \cdot(\mu '-\nu ')|< \omega^1 \text{ if $ (\mu ', \nu ') < (\mu, \nu) $} \}.
\end{aligned}\end{equation}
%% Notice that
%% \begin{eqnarray}
%% \mathbf{M}=
%% \{(\mu,0)\in\N_0^{\mathbf{n}}\times\N_0^{\mathbf{n}}:\;\mu\in\mathbf{M}_0\}
%% \cup
%% \{(0,\nu)\in\N_0^{\mathbf{n}}\times\N_0^{\mathbf{n}}:\;\nu\in\mathbf{M}_0\}.
%% \end{eqnarray}
Notice that
\begin{equation} \label{eq:defM1}\begin{aligned} \text{if $(\mu, \nu) \in \mathbf{M}$ we have either
$\mu =0$ and $\nu\in\mathbf{M}_0$,
or
$\nu =0$ and  $\mu\in\mathbf{M}_0$}.
\end{aligned}\end{equation}
In \cite{CM1, Cu3} it is shown that  for $G_{\mu \nu}^0:=G_{\mu \nu }(p^0,0)$ we have
\begin{equation} \label{eq:equation f}\begin{aligned} & (X_H^ {(0)})_f =
\mathcal{H}_{p^1}f + \sum_{j=1,...,7} (\partial_{ \Pi_j(f)} H) P_c
{\mathbbm{i}} \sigma_3 \Diamond_j f - \sum_{(\mu, \nu) \in \mathbf{M} }
z^\mu \overline{z}^\nu
 G_{\mu \nu}^0 + R_1+ R_2,
\end{aligned}\end{equation} $P_c$ the projection on $X_c$ in
 \eqref{eq:spectraldecomp},
and there is a constant $ C(C_0
)$ independent of $\epsilon$
s.t.
\begin{equation}
\label{bound1:z1} \| R_1 \|_{L^1_t([0,T],H^{ 1 })}+\|
 R_2 \|_{L^{2
 }_t([0,T],W^{ 1, \frac{6}{5}})}\le C(C_0
) \epsilon^2.
\end{equation}
We sketch briefly this point. With $\widehat{\nabla} _f$ defined in \eqref{eq:part f mod}, we define
\begin{align*} &  R_2=\sum_{(\mu, \nu) \in \mathbf{M} }
z^\mu \overline{z}^\nu
 \(   G_{\mu \nu}^0 -G_{\mu \nu}\)   -  {\mathbbm{i}} \sigma_3\widehat{\nabla} _f \mathbf{R}_2-  {\mathbbm{i}} \sigma_3 B_2 f, \end{align*}
where the last term is defined schematically from $\widehat{\nabla} _f \< B_2,f^2  \> \sim \< \widehat{\nabla} _f B_2,f^2  \> + B_2 f. $ Then
 the desired estimate on $R_2$ in \eqref{bound1:z1} is elementary.
 For example
 \begin{equation*}
   \|
 B_2 f\|_{L^{2
 } ([0,T],W^{ 1, \frac{6}{5}})} \le  \|
 B_2  \|_{L^{\infty
 } ([0,T],L^{ 3/2})} \|
   f\|_{L^{2
 } ([0,T],W^{ 1,   {6} })} \lesssim \epsilon  \|
   f\|_{L^{2
 } ([0,T],W^{ 1,   {6} })}  \lesssim \epsilon ^2
 \end{equation*}
by \eqref{eq:canc B2} and \eqref{Strichartzradiation} in $[0,T]$.   $R_1$ is formed by the other terms and it is standard to show that it satisfies the bound \eqref{bound1:z1}. For example for $2\le d \le 4$
\begin{align*}
   \| \< \widehat{\nabla} _f B_d,f^d  \> \| _{L^{1
 }_t  H^1_x}  \le   \|   \sup _{\| g \| _{H^{-1}} =1}\< \widehat{\nabla} _f B_d g,f^d  \> \| _{L^{1
 }_t  }   &\le \left \|  \sup _{\| g \| _{H^{-1}} =1} \|  \widehat{\nabla} _f B_d g \| _{\Sigma _k} \| f ^d\|  _{L ^{6/d}_x} \right \| _{L^{1
 }_t  } \\& \lesssim  \| f \| _{L^2_t L^6_x}^2 \| f \| _{L^\infty_t H^1_x} ^{d-2} \lesssim   \epsilon ^d
 \end{align*}
 and for $d=3,4$, for  $(d-1, q_d)$ admissible,
 \begin{align}\label{eq:d=3,4}
   \|   B_d f ^{d-1}   \| _{L^{1
 }_t  H^1_x}  \lesssim   \|    f     \| _{L^{\infty
 }_t  H^1_x}    \|    f     \| ^{d-1
 } _{L^{d-1
 }_t  L ^{q_d}_x}   \lesssim \epsilon ^{d-1
 }.
 \end{align}
 The $d=5$ term can be treated similarly, but has an additional part, when the
 $f$ derivative is applied to the $\zeta $ variable in   \eqref{H5power2}.
 But the resulting term is like \eqref{eq:d=3,4}  for $d=6$. Finally,
 $  \|  \nabla E_P(f)   \| _{L^{1
 }_t  H^1_x} \lesssim \epsilon ^{2
 }$ by hypotheses (H1)--(H2).
Having discussed \eqref{bound1:z1},
  by  \eqref{eq:X_F1}  we get
\begin{equation} \label{eq:equation f-1}\begin{aligned} X_H^ {(1)} = \resto^{0,0}_{\infty,\infty}(\Pi (f))
& \left [
\langle \Diamond f, \mathcal{H}_{p^1}f \rangle + (\partial_{ \Pi (f)} H)\langle \Diamond f, {\mathbbm{i}} \sigma_3 \Diamond f\rangle + \langle \Diamond f, R_1+ R_2 \rangle \right. \\& - \sum_{(\mu, \nu) \in \mathbf{M}}
z^\mu \overline{z}^\nu
 \left. \langle \Diamond f, G_{\mu \nu}^0 \rangle \right ] P_c {\mathbbm{i}} \sigma_3 \Diamond f.
\end{aligned}\end{equation}
Then, for $\mathbf{v}$  obtained summing contributions from   \eqref{eq:equation f-1}
and   the $\sum _{j=1,...,7}$ in \eqref{eq:equation f}, we obtain
\begin{equation} \label{eq:equation f1}\dot f  -\left (
\mathcal{H}_{p^1}f + P_c
{\mathbbm{i}} \sigma_3 \mathbf{v}\cdot \Diamond f\right )=  - \sum_{(\mu, \nu) \in \mathbf{M}}
z^\mu \overline{z}^\nu
 G_{\mu \nu}^0 + R_1+ R_2.
\end{equation}
It is easy to see,
 from \eqref{Strichartzradiation}--\eqref{L^inftydiscrete}
and
\eqref{bound1:z1}, that
 \begin{equation} \label{eq:bound v}\begin{aligned} & \| \mathbf{v} \|_{L^1 ([0, T], \R^7) +L^\infty ([0, T], \R^7)}\le C(C_0
) \epsilon.
\end{aligned}\end{equation}
Strichartz and smoothing estimates on $f$ are a consequence of well-known estimates for the group $e^{t\mathcal{H}_{p^1}} P_c$
which resemble those valid for $e^{{\mathbbm{i}}t \Delta}$, see \cite{MR2805462} for references.

 To deal with the term $ P_c
{\mathbbm{i}} \sigma_3 \mathbf{v}\cdot \Diamond f$, where the operator $ P_c
{\mathbbm{i}} \sigma_3 \mathbf{v}\cdot \Diamond  $
which does not commute with $\mathcal{H}_{p^1}$ we  adopt an ideas by Beceanu \cite {MR2831875}. We consider the system $\dot f =
{\mathbbm{i}} \sigma_3 \mathbf{v}\cdot \Diamond f$,
writing it in the form
 \begin{equation} \label{eq:int fa} \begin{aligned} & \dot f = A(t)f+B(t)f \text{ with } A(t):= \sum_{j=1,...,4}{\mathbbm{i}} \sigma_3 \mathbf{v}_j(t) \Diamond_j \text{ and }B(t):= \sum_{j=5,6,7}{\mathbbm{i}} \sigma_3 \mathbf{v}_j (t)\Diamond_j.
\end{aligned} \end{equation}
Since $A(t)$ and $B(t)$ commute and the terms of the sum defining $A(t)$ commute,
if we denote by $W (t,s)$ the fundamental solution of the system
\eqref{eq:int fa},
that is,
 \begin{equation} \label{eq:fu sol} \begin{aligned} & \partial_t W (t,s)= (A(t)+ B(t)) W (t,s) \text{ with } W (s,s)=I,
\end{aligned} \end{equation}
and by $W_A(t,s)= e^{\int_s^tA(t')\,dt'}$ (resp. $W_B(t,s)$)
the fundamental solution of $\dot f = A(t)f$ (resp. $\dot f = B(t)f$),
then we have $ W (t,s)= W_A (t,s) W_B (t,s)$.
%%and $W_A (t,s)= e^{\int_s^tA(t')\,dt'}$.
\begin{lemma}
 \label{lem:error1} Let $M>5/2$ and $\alpha\in [0,1/2)$.
Then there exists a   constant $C>0$ dependent only on $M$ such that for all $s<t$ in $[0,T]$
 \begin{align} & \| \langle x \rangle^{-M}\left (W (t,s) -1\right) e^{ {{\mathbbm{i}} \sigma_3 (\Delta -\omega^{1}) (t-s)}} \langle x \rangle^{-M } \|_{B (L^2, L^2)}\le C \psi_\alpha (t -s) \| \mathbf{v}\|^\alpha_{L^1 ([s,t]) +L^\infty ([s,t])} \nonumber\\& \text{ with $\psi_\alpha (t) =
\langle t\rangle^{- \frac{3}{2} +\alpha }$ for $t\ge 1$ and $\psi_\alpha (t) =
 t^{- \alpha }$ for $t\in (0,1)$.}\label{eq:psialpha} \end{align}

 \end{lemma}
\proof We have
\begin{equation}\label{eq:error11} \begin{aligned} & W (t,s) -1= \left [(W_A(t,s) -1) W_B(t,s)\right ] +\left [ W_B(t,s) -1\right ].
 \end{aligned} \end{equation}
In the 1st term in the r.h.s.  $W_B(t,s)$ commutes with the other  operators   and is an isometry in $L^2$:
\begin{equation} \begin{aligned} & \| \langle x \rangle^{-M}W_A(t,s) -1) W_B(t,s) e^{ {{\mathbbm{i}} \sigma_3 (\Delta -\omega^{1}) (t-s)}} \langle x \rangle^{-M } \|_{B (L^2, L^2)} \\& = \| \langle x \rangle^{-M}W_A(t,s) -1) e^{ {{\mathbbm{i}} \sigma_3 (\Delta -\omega^{1}) (t-s)}} \langle x \rangle^{-M } \|_{B (L^2, L^2)}.
 \end{aligned}\nonumber \end{equation}
Then the desired estimate of this
is that of \cite[Lemma 9.4]{CM1}.
We next consider the 2nd term in the r.h.s. of \eqref{eq:error11}.
By the commutation properties of $W_B(t,s)$ we are reduced to bound
\begin{equation} \begin{aligned} & \| \langle x \rangle^{-M} e^{ {{\mathbbm{i}} \sigma_3 (\Delta -\omega^{1}) (t-s)}} \langle x \rangle^{-M } \|_{B (L^2, L^2)} \left (\int_s^t \| B(t') W_B(t',s)dt' \|_{B (L^2, L^2)} \right)^{\alpha}.
 \end{aligned}\nonumber \end{equation}
The first factor is bounded by $c_0 \langle t-s\rangle^{-\frac{3}{2}}$ while the second by $ |t-s|^{\alpha} \| B \|_{L^\infty ((s,t), B (L^2, L^2)) }^{\alpha}$,
where the last factor is bounded by $\| \mathbf{v} \|_{L^\infty ((s,t), \R^7) }^{\alpha}$.
\qed

  \begin{proposition}\label{thm:strich}
Let $F(t)$ satisfy $P_cF(t)=F(t)$
Consider the equation
 \begin{equation} \label{eq:strich1} \dot u -
\mathcal{H}_{p^1}u - P_c
{\mathbbm{i}} \sigma_3 \mathbf{v}\cdot \Diamond u=F.
 \end{equation}
 Then there exist fixed $\sigma >3/2 $, and an $\epsilon_0>0$ such that if $ \epsilon \in (0,\epsilon_0)$ then  we have
\begin{equation}\label{eq:strich3} \| u \|
_{L^p([0,T],W^{1,q})}\le
 C (\| P_cu(0) \|_{H^1 } + \| F \|
_{L^2([0,T],H^{1, \sigma})+ L^1([0,T],H^{1 })}) \text{ $\forall$ admissible pairs $ (p,q)$}.
\end{equation}
\end{proposition}
Before the proof,  we observe that Proposition~\ref{thm:strich} implies
the following.
\begin{corollary}\label{prop:conditional4.2} Under the hypotheses of Theorem~\ref{thm:mainbounds}
there exist two constants $c_0$ and $\epsilon_0>0$ such that if $ \epsilon \in (0,\epsilon_0)$ then
\begin{equation}\label{4.5}
 \| f \|_{L^p_t([0,T],W^{ 1,q}_x)}\le c_0 \epsilon +c_0 \sum_{(\mu, \nu) \in \mathbf{M}} \| z^{\mu +\nu }\|_{L^2(0,T)} \text{ for any admissible pair $(p,q)$.}
\end{equation}
\end{corollary}

For the  elementary proof of this corollary see for instance \cite[Lemma 8.1]{CM1}.

\noindent \textit{Proof of Prop.~\ref{thm:strich}}. We follow \cite{MR2831875}.
Denote $u_0=P_cu(0)$.
We set $P_d:=1-P_c$, fix $\delta >0$ and consider
\begin{equation}\label{ast}
\begin{aligned} & \dot Z -\mathcal{H}_{p^1}P_cZ
 -P_c{\mathbbm{i}} \sigma_3 \mathbf{v}\cdot \Diamond P_cZ=F - \delta P_d Z
 \,, \quad Z(0) =u_0.\end{aligned}
 \end{equation}
Notice that, see \eqref{eq:lincom},
\begin{equation}\label{eq:potlin}
 \text{$ \mathcal{H}_{p^1}= {\mathbbm{i}} \sigma_3 (-\Delta + \omega^{1})+ V $ with $V\in \mathcal{S}(\R^3, B(\C^2, \C^2))$;}
\end{equation}
we then rewrite
\eqref{ast}
as
\begin{equation*} \begin{aligned} & \dot Z - {\mathbbm{i}} \sigma_3 (\Delta - \omega^{1})Z
 -{\mathbbm{i}} \sigma_3 \mathbf{v}\cdot \Diamond Z=F +V_1V_2 Z -\widetilde{P}_d (\mathbf{v}) Z \text{ with } Z(0) =u_0, \end{aligned}
 \end{equation*}
 $\widetilde{P}_d (\mathbf{v}):= P_d {\mathbbm{i}} \sigma_3 \mathbf{v}\cdot \Diamond + {\mathbbm{i}} \sigma_3 \mathbf{v}\cdot \Diamond P_d$ and $V_1V_2 =V -\mathcal{H}_{p^1}P_d
 - \delta P_d$ with $V_2(x)$ a
smooth exponentially decaying and invertible matrix, and with
the multiplication operator $V_1:H^{k,s'}\to H^{k,s}$
bounded for all
$k$, $s$ and $s'$. We have:
\begin{align}
 &Z (t) = W(t,0)e^{ {\mathbbm{i}} \sigma_3 (-\Delta + \omega^{1}) t}Z (0)\label{eq:423} \\& + \int_0^t
 e^{ {\mathbbm{i}} \sigma_3 (-\Delta + \omega^{1}) (t-t')}W(t,t')
 \left [ F(t')+ V_1 V_2 Z (t ') - \widetilde{P}_d (\mathbf{v}(t ')) Z(t ') \right ]\,dt'. \nonumber\end{align}
 For arbitrarily fixed pairs $(K,S)$ and
$(K',S')$ there exists a constant $C$ such that we have
 \begin{equation*}
\| \widetilde{P}_d (\mathbf{v}) V_2^{-1 } \|_{B (H^{-K',-S'},H^{ K, S })} \le C \epsilon.
 \end{equation*}
 By picking $\epsilon$ small
enough, we can assume that
 the related operator norm is small. We have
 \begin{equation*} \label{eq:4231} \begin{aligned}
 & \| Z \|_{L^p_t W^{ q,1}
 \cap L^2_t H^{ k,-\tau_0} }\le
 C\| Z(0) \|_{H^{ 1 }}+C \| F \|_{L^1_t H^{ 1 }
 +L^2_t H^{ 1, \tau_0} } \\& +
 \| V_1 - \widetilde{P}_d (\mathbf{v}(t ')) V_2^{-1 } \|_{
 L^\infty_t(B(H^{ 1 }, H^{ 1, \tau_0})) } \| V_2 Z (t)
 \|_{ L^2_t H^{ 1 } }. \end{aligned}
 \end{equation*}
 For $\widetilde{T}_0 f (t)= V_2\int_0^t
e^{ {\mathbbm{i}} \sigma_3 (-\Delta + \omega^{1}) (t-t')}W(t,t') V_1 f(t')
dt',$ by \eqref{eq:423}, we obtain:
\begin{equation}\begin{aligned}
 & (I- \widetilde{T}_0) V_2Z (t)=
 V_2 W(t,0)e^{ {\mathbbm{i}} \sigma_3 (-\Delta + \omega^{1}) t}Z (0) \\& - V_2\int_0^t
 e^{ {\mathbbm{i}} \sigma_3 (-\Delta + \omega^{1}) (t-t')}W(t,t')
 \left [ F(t') - \widetilde{P}_d (\mathbf{v}(t ')) Z(t ') \right ]\,dt' \end{aligned}\nonumber \end{equation}
We then obtain the desired result
if we can show that \begin{equation} \label{ns1}\| (I-
 \widetilde{{T}}_0)^{-1} \|_{ L^2([0,T), H^1(\R^3)) \righttoleftarrow} <C_1,
\end{equation} for $\epsilon C_1 $ smaller than a fixed number. Thanks to Lemma~\ref{lem:error1} it is enough to prove \eqref{ns1} with
$\widetilde{T}_0$ replaced by
\begin{equation}\begin{aligned}
 & {T}_0 f (t)= V_2\int_0^t
e^{ {\mathbbm{i}} \sigma_3 (-\Delta + \omega^{1}) (t-t')} V_1 f(t')
dt'. \end{aligned}\nonumber
\end{equation}
Set
\begin{equation}\begin{aligned}
 & {T}_1 f (t)= V_2\int_0^t
 e^{ (-\mathcal{H}_{p^1}P_c +\delta P_d)(t'-t)} V_1 f(t')\,dt'.
 \end{aligned}\nonumber \end{equation}
By \cite{MR1835384} we have $\| T_1 \|_{L^2([0,T), H^1(\R^3)) \righttoleftarrow} <C_2$ for a fixed $C_2$. By elementary arguments, see \cite{NakanishiSchlag},
\begin{equation}\begin{aligned}
 & (I- {T}_0) (I+ {T}_1)= (I+ {T}_1)
 (I- {T}_0)=I.
 \end{aligned}\nonumber \end{equation}
This yields \eqref{ns1} with $ \widetilde{{T}}_0$ replaced by $
{T}_0$ and with $C_1=1+C_2$. \qed

Now we turn to the equations $\dot z_l = \im \partial_{\overline{z}_l}H$.
We will prove the following.
%\ac{will we actually prove this?? THIS INEQUALITY IS MORE OR LESS PROVED IN THE REMAINDER OF THE SECTION}

%\ac{where do we use the following proposition?IT IS CRUCIAL. First of all it proves {proposition} \ref{prop:mainbounds}
%in the case of inequality \eqref{L^2discrete}. Then it does the same for inequality \eqref{Strichartzradiation}
%thanks to Corollary \ref{prop:conditional4.2}. As for {proposition} \ref{prop:mainbounds}
%in the case of inequality \eqref{L^inftydiscrete}, this follows directly from \eqref{prop:est z}
%and the equation satisfied by $z$}

\begin{proposition}\label{prop:est z} There exists a fixed $c_0>0$ and a constant $ \epsilon_0>0$ which depends on $C_0$ such that
\begin{equation}\label{eq:est z1}\begin{aligned}
 &
\sum_l | z
_l(t)|^2 + \sum_{(\mu, \nu) \in \mathbf{M}} \| z^{\mu +\nu }\|_{L^2(0,t)}^{2}\le c_0 (1+C_0)\epsilon^2,
\qquad\forall t\in [0,T],\quad\forall\epsilon\in (0, \epsilon_0).
 \end{aligned} \end{equation}
 \end{proposition}
Proposition~\ref{prop:est z} allows to conclude the proof of Proposition~\ref{prop:mainbounds}.
The proof of Proposition~\ref{prop:est z} follows a series of standard steps, and is
basically the same as the analogous proof in \cite{MR2805462}, or in \cite{MR2843104}.

 The first step in the proof of Proposition~\ref{prop:mainbounds} consists in splitting $f$ as follows: \begin{equation}
 \label{eq:g variable}
g=f+ Y \,, \quad Y:=-\im \sum_{(\mu, \nu) \in \mathbf{M} }
z^\mu \overline{z}^\nu
 R^{+}_{\im \mathcal{H}_{p^1
 } } (\mathbf{{e}} \cdot(\nu-\mu))
 G_{\mu \nu}^0,
\end{equation}
where $R^{+}_{\im \mathcal{H}_{p^ 1
 } }$ is extension from above of the resolvent and makes sense because the theory of Jensen and Kato  \cite{JK} holds also for these operators, see for example Perelman \cite[Appendix 4]{perelman}.

\noindent The part of $f$ that acts   effectively on the variables $z$ will be shown to be $Y$, while $g$ is small, thanks to the following lemma.
\begin{lemma}\label{lemma:bound g} For fixed $s>1$ there exist
a fixed $c$
such that if $\epsilon_0$ is sufficiently small we have $\| g
\|_{L^2((0,T), H^{0,-s} (\R^3, \C^4)}\le c \epsilon $.
\end{lemma}
\proof In the same way as the proof of Proposition~\ref{thm:strich}
(which we wrote explicitly) is
similar to analogous proofs valid for the scalar NLS \eqref{Eq:scalNLS}, the proof of
Lemma~\ref{lemma:bound g} is analogous to the proof of \cite[Lemma 8.5]{CM1} contained in
\cite[Sect. 10]{CM1} and is skipped here. The only difference between \cite{CM1} and the present situation
is notational,
in the sense that
inside \eqref{eq:strich1}
one has
$
{\mathbbm{i}} \sigma_3 \mathbf{v}\cdot \Diamond u = {\mathbbm{i}} \sigma_3 \sum_{j\le 7}\mathbf{v}_j \Diamond_ju$,
as opposed to \cite[(10.1)]{CM1},
where the corresponding terms are
$ {\mathbbm{i}} \sigma_3 \sum_{j=1}^{4}\mathbf{v}_j \Diamond_ju$. But this does not make any difference in the proof because
what matters is simply that each $\Diamond_j$ commutes with $-\Delta + \omega^1$, which was
used to get \eqref{eq:423}.
\qed

 Now we examine the equations on $z$. We have
\begin{equation*}\label{eq:FGR0} \begin{aligned} & -\im \dot z_j
=
\partial_{\overline{z}_j}(H_2 +Z_0+Z_1+ \mathcal{R}).
\end{aligned} \end{equation*}
When we substitute \eqref{eq:g variable} and we set $R_{\mu \nu }^+:=R_{ \im \mathcal{H}_{p^1} }^+ (\mathbf{e}\cdot (\nu -\mu
))$ we obtain
\begin{equation}\label{eq:FGR0-1} \begin{aligned} & -\im \dot z_l-\partial_{\overline{z}_l} H_2
=
\partial_{\overline{z}_l} Z_0  +\im \sum_{ \substack{(\alpha, \beta), (\mu, \nu) \in \mathbf{M}
 }} \nu_l\frac{z^{\mu +\alpha } \overline{{z }}^ { {\nu}
+\beta}}{\overline{z}_l} \langle R_{ \alpha \beta}^+G^0_{ \alpha
\beta },{\mathbbm{i}} \sigma_3G_{\mu \nu }\rangle
 \\ & + \sum_{(\mu, \nu) \in \mathbf{M}} \nu_l\frac{z^\mu
 \overline{ {z }}^ { {\nu} } }{\overline{z}_l} \langle g,
 {\mathbbm{i}} \sigma_3 G
_{\mu \nu }\rangle +
\partial_{ \overline{z}_l} \resto.
\end{aligned} \end{equation}
Using \eqref{eq:defM1}, we rewrite this as
\begin{align} \label{equation:FGR1}& -\im \dot z_j-\partial_{\overline{z}_j} H_2
=
\partial_{\overline{z}_j} Z_0 + \sum_{
 (\mu, \nu) \in \mathbf{M}} \nu_j\frac{z^\mu
 \overline{ {z }}^ { {\nu} } }{\overline{z}_j} \langle g,
 {\mathbbm{i}} \sigma_3 G
_{\mu \nu }\rangle +
 \mathcal{E}_j
\\ & \label{equation:FGR12} +\im\sum_{ \substack{
\beta,\nu \in \mathbf{M}_0
 }} \nu
_j\frac{ \overline{{z }}^ {\nu +\beta } }{\overline{z}_j}\langle
R_{ 0 \beta}^+
 { G}_{ 0\beta }^0, {\mathbbm{i}} \sigma_3G^0_{0 \nu }\rangle
 \\ & \label{equation:FGR13} +\im \sum_{ \substack{
\alpha,\nu\in \mathbf{M}_0
 }} \nu
_j\frac{z^{ \alpha } \overline{{z }}^ {\nu }
}{\overline{z}_j}\langle R_{ \alpha 0 }^+
 G_{ \alpha 0}^0, {\mathbbm{i}} \sigma_3G^0_{0 \nu }\rangle.
\end{align}
Here the elements in \eqref{equation:FGR12} can be eliminated
through a new change of variables that we will see momentarily and
$\mathcal{E}_j$ is a remainder
term defined by
\begin{equation} \label{eq:def e}\begin{aligned} & \mathcal{E}_j:=
\sum_{ (\mu, \nu) \in \mathbf{M}} \nu_j\frac{z^\mu
 \overline{ {z }}^ { {\nu} } }{\overline{z}_j} \langle g,
 {\mathbbm{i}} \sigma_3 G
_{\mu \nu }\rangle +
\partial_{ \overline{z}_j} \resto -\text{\eqref{equation:FGR12}}- \text{\eqref{equation:FGR13}}.
\end{aligned}
\end{equation}
Set $\zeta_l =z_l +F_l(z,\overline{z})$  with
\begin{equation*}\label{equation:FGR2} \begin{aligned} &
F_l(z,\overline{z})= \sum_{ \substack{
\beta,\nu \in \mathbf{M}_0
 }}
\frac{ \im \nu_l\overline{{z
}}^ {\nu +\beta }  }{\mathbf{e} \cdot (\beta +\nu) \overline{z}_l}  \langle R_{ 0 \beta}^+
 { G}_{ 0 \beta }^0, {\mathbbm{i}} \sigma_3G_{0 \nu }^0\rangle
- \sum_{ \substack{
\alpha, \nu \in \mathbf{M}_0
\\
\mathbf{e}\cdot \alpha \neq
 \mathbf{e}\cdot \nu
 }}\frac{ \im \nu_lz^{ \alpha } \overline{ z}^ { \nu
}}{\mathbf{e}
 \cdot (\alpha - \nu) \overline{z}_l}   \langle R_{ \alpha 0 }^+ G^0_{ \alpha 0},{\mathbbm{i}} \sigma_3 G_{0 \nu }^0\rangle.
\end{aligned}
\end{equation*}
This change of variables is such that, setting $F=(F_1,...,F_\mathbf{n})$, we get
\begin{equation*}  \begin{aligned}   \mathfrak{L}_j(z,f) | _{f=0}:&=\sum _{l=1,...,\mathbf{n}} \(    \partial _{ \overline{{z}}_l}F_j (z,\overline{z}) \partial _{  {{z}}_l}H_2(z,0) - \partial _{ {z}_l}F_j (z,\overline{z}) \partial _{ \overline{{z}}_l}H_2(z,0)\)   \\& =\partial _{ \overline{{z}}_l}H_2(F (z,\overline{z}),0)+ \text{\eqref{equation:FGR12}}+ \text{\eqref{equation:FGR13}}.
\end{aligned}
\end{equation*}
Furthermore,   by
$\nu\in \mathbf{M}_0$,
which implies $\nu \cdot \mathbf{e}>\omega^1$, we have $| {\nu} |>1$. Then by \eqref{L^2discrete}--\eqref{L^inftydiscrete}
\begin{equation} \label{equation:FGR3} \begin{aligned} & \| \zeta -
 z \|_{L^2(0,T)} \le C \epsilon \sum_{\substack{
\alpha\in \mathbf{M}_0
}} \| z^{\alpha }\|_{L^2(0,T)}
\le C(C_0) \epsilon^2,
\quad  \| \zeta -
 z \|_{L^\infty (0,T)} \le C(C_0)\epsilon^3.
\end{aligned}
\end{equation}
In the new  $\zeta$ variables, \eqref{equation:FGR1} takes the form
\begin{equation} \label{equation:FGR4} \begin{aligned} &
 -\im \dot \zeta
_j=
\partial_{\overline{\zeta}_j}H_2 (\zeta, f) +
\partial_{\overline{\zeta}_j}Z_0 (\zeta, f)+ \mathcal{D}_j
+ \im \sum_{ \substack{
\alpha, \nu \in \mathbf{M}_0
\\
\mathbf{e}\cdot \alpha =
 \mathbf{e}\cdot \nu
 }}\nu_j \frac{\zeta^{
\alpha } \overline{ \zeta}^ { \nu }}{\overline{\zeta}_j} \langle R_{
\alpha 0 }^+ G^0_{ \alpha 0},{\mathbbm{i}} \sigma_3 G^0_{0 \nu
}\rangle,
\end{aligned}
\end{equation}
with for $A_l=$r.h.s. of \eqref{eq:FGR0-1},
\begin{equation} \label{eq:eq D} \begin{aligned} &
   \mathcal{D}_j = \mathcal{E}_j+ \mathfrak{L}_j(z,0)-\mathfrak{L}_j(z,f)  +\sum _{l=1,...,\mathbf{n}} \( \partial _{ {z}_l}F_j (z,\overline{z}) A_l  - \partial _{ \overline{{z}}_l}F_j (z,\overline{z}) \overline{A} _l\)  .
\end{aligned}
\end{equation}
 From these equations by
$\sum_l \mathbf{e}_l
\big(\overline{\zeta}_l
\partial_{\overline{\zeta}_l}(H_2+Z_0) - {\zeta}_j
\partial_{ {\zeta}_l}(H_2+Z_0)\big) =0$ we get
\begin{equation} \label{eq:FGR5} \begin{aligned}
 &\partial_t \sum_{l=1,..., \mathbf{n} } \mathbf{e}_l
 | \zeta_l|^2 = 2 \sum_{l=1,..., \mathbf{n} }\mathbf{e}_l{\Im} \left (
\mathcal{D}_l\overline{\zeta}_l \right) +2 \sum_{ \substack{ \alpha,\nu\in \mathbf{M}_0\\
\mathbf{e}\cdot \alpha =
 \mathbf{e}\cdot \nu
 }} \mathbf{e}\cdot \nu
\Re \Big(\zeta^{ \alpha } \overline{\zeta }^ { \nu } \langle
R_{ \alpha 0}^+ G_{ \alpha 0}^0, {\mathbbm{i}} \sigma_3 G^0_{0\nu
}\rangle \Big).
\end{aligned}
\end{equation}
\begin{lemma}
\label{lemma:FGR1} Assume inequalities \eqref{Strichartzradiation}--\eqref{L^inftydiscrete}. Then for a
fixed constant $c_0$ we have
\begin{equation}\label{eq:FGR7}
 \sum_{j=1,...,\mathbf{n}}\|{\Im} \left (
\mathcal{D}_j\overline{\zeta}_j \right)\|_{
L^1[0,T]}\le (1+C_0)c_0 \epsilon^{2}
.
\end{equation}
\end{lemma}
\textit{Proof (sketch).} For a detailed proof  we refer to   \cite[Appendix B]{MR2843104}: here we
give a sketch.
First of all, we consider the contribution of  $\mathcal{E}_j$. This, in turn, is a sum of various terms. For the terms originating from $\mathbf{R} _3$, cf. Lemma \ref{lem:ODE1}, we have
\begin{align*}  \|  \langle \partial _{\overline{z}_j}B_{d } (p^0_4 ,\Pi (f), z,f), f^{ d} \rangle  {\zeta}_j\|_{
L^1_t}\le      \|  f \| _{L_t^d L ^{p_d}_x}    ^d \| \zeta \| _{L^\infty _t}\lesssim \epsilon ^{d+1},
\end{align*}
with $(d,p_d)$ admissible, and for $d=2,3,4,5$. For the following term,  we claim
\begin{align}\label{contrR2}  \|   \partial _{\overline{z}_j}\mathbf{R} _2  {\zeta}_j\|_{
L^1_t} \lesssim \epsilon ^{3}.
\end{align}
From Lemma \ref{lem:ODE1} we know that $\mathbf{R} _2$ is basically
a sum of degree $2\mathbf{N}+3$  monomials in $(z,\overline{z},f)$, which are at most degree 1 in $f$. Let us take a term which is degree 0 in $f$. Then its  $\partial _{\overline{z}_j}$
derivative is in absolute value bounded above by a term $ | z^{\mu +\nu }  | $
with $|\mu |+|\nu | \ge 2\mathbf{N}+2$.
So we can write it as $  | z^{\alpha  +\beta +\gamma }  |      $ with $|\alpha |\ge \mathbf{N}+1$, $|\beta |\ge \mathbf{N}+1$.  But then $\alpha \cdot \mathbf{e} > \omega ^1$
$\beta \cdot \mathbf{e} > \omega ^1$. Then
 \begin{align*}  \|   z^{\alpha  +\beta +\gamma }  {\zeta}_j\|_{
L^1_t} \le \|   z^{\alpha   }   \|_{
L^2_t} \|   z^{ \beta   }   \|_{
L^2_t}\|      {\zeta}_j\|_{
L^\infty_t}  \lesssim  \epsilon ^{3}.
\end{align*}
Terms degree 1 in $f$ can be treated similarly, yielding \eqref{contrR2}.  We claim also
\begin{align}\label{contrR2-}  \|     \nu _j  \frac{z^{\mu +\alpha } \overline{z}^{\nu +\beta } }{\overline{z}_j}  {\zeta}_j\|_{
L^1_t} \lesssim \epsilon ^{3} \text{ for $| (\mu -\nu )\cdot \mathbf{e}|>\omega ^1$ and $(\mu ,\nu ) \not \in \mathbf{M}$. }
\end{align}
In this case we can write $ z^\mu \overline{z}^\nu  =  z^{\mu  '    } \overline{z}^ {\nu '}    z^{ \gamma } \overline{z}^{\delta}$
with $(\mu ' , \nu ' ) \in \mathbf{M}$ and   $|\gamma |+|\delta |>0$. Then we consider
\begin{align*}   \nu _j  \frac{z^{\mu +\alpha } \overline{z}^{\nu +\beta }  }{\overline{z}_j}  \overline{{\zeta}}_j =  \nu _jz^{\mu '+\alpha } \overline{z}^{\nu '+\beta }     z^{ \gamma } \overline{z}^{\delta} +\nu _j \frac{z^{\mu '+\alpha } \overline{z}^{\nu '+\beta } }{\overline{z}_j}  (\overline{{\zeta}}_j-\overline{z}_j).
\end{align*}
By \eqref{L^2discrete}--\eqref{L^inftydiscrete}
\begin{align*}   \|   z^{\mu '+\alpha } \overline{z}^{\nu '+\beta }     z^{ \gamma } \overline{z}^{\delta}  \| _{L^1_t} \lesssim  \|   z^{\mu '  } \overline{z}^{\nu '  }      \| _{L^2_t}\|   z^{ \alpha } \overline{z}^{ \beta }    \| _{L^2_t}\|     z \| _{L^\infty_t}^{|\gamma |+|\delta |} \lesssim \epsilon ^{3}.
\end{align*}
and by \eqref{equation:FGR3}
\begin{align*}   \|  \nu _j \frac{z^{\mu  +\alpha } \overline{z}^{\nu  +\beta } }{\overline{z}_j}  (\overline{{\zeta}}_j-\overline{z}_j) \| _{L^1_t} \lesssim   \|   z^{ \alpha } \overline{z}^{ \beta }    \| _{L^2_t}\|     z -\zeta \| _{L^2_t}  \lesssim \epsilon ^{3}.
\end{align*}
 This yields \eqref{contrR2-}. By similar arguments, one can prove
\begin{align*}   \|     \nu _j  \frac{z^{\mu   } \overline{z}^{\nu   } }{\overline{z}_j}  \< g, {\mathbbm{i}} \sigma_3G _{\mu \nu}   \> {\zeta}_j\|_{
L^1_t} \lesssim \epsilon ^{3} \text{ for $| (\mu -\nu )\cdot \mathbf{e}|>\omega ^1$ and $(\mu ,\nu ) \not \in \mathbf{M}$. }
\end{align*}
We next consider the following, see   \cite[Lemma B.1]{MR2843104},
\begin{align}\label{contrR2-}  \|    \overline{ \partial}_{ j}(Z_0 (\zeta, f) -Z_0(z,f) ) {\zeta}_j\|_{
L^1_t} \lesssim \epsilon ^{3}  .
\end{align}
   Is enough to consider   $  z ^{\alpha}
\frac{  \overline{z} ^{ {\beta} }}{ \overline{z }_j}   \overline{\zeta} _j   - \zeta
^{\alpha} \frac{  \zeta ^{ {\beta} }}{  \overline{\zeta} _j} \overline{\zeta} _j   $ with $\mathbf{e} \cdot \alpha =
\mathbf{e} \cdot \beta  $ and $\beta _j>0$. By Taylor expansion these
are
\begin{equation*}\label{other_ta}  \sum  _k\partial _k
\left (\frac{z  ^{\alpha}   \overline{z}  ^{ {\beta} }}{  \overline{z}
_j}\right ) (\zeta  _k-z  _k)  \overline{\zeta} _j + \sum
_k\overline{\partial} _k \left (\frac{z ^{\alpha}   \overline{z} ^{
{\beta} }}{  \overline{z} _j}\right ) (  \overline{\zeta}  _k-  \overline{z} _k)  \overline{\zeta}
_j +\overline{\zeta} _j O(|z -\zeta |^2).
\end{equation*}
The remainder term is the easiest, the other two terms similar.
Substituting  the definition of $\zeta$, a typical term in the first
summation is $ \frac{z ^{\alpha +A}   \overline{z} ^{B+\beta}    }{|z _k|^2} $, with   $\alpha \cdot \mathbf{e}>\omega ^1$, $\beta\cdot \mathbf{e}>\omega ^1$, $A\cdot \mathbf{e}>\omega ^1$ and
$B\cdot \mathbf{e}>\omega ^1$.
and with $\alpha _k\neq 0\neq B_k$.   By
(H8),  $\mathbf{e}
\cdot \alpha =\mathbf{e} \cdot \beta$ implies that
   there is at least
one index $\beta _\ell \neq 0$ such that $\mathbf{e} _\ell = \mathbf{e} _k$.
Then, by the fact that   monomials $z ^\alpha   \overline{z} ^\beta $ in $Z_0$
are such that $|\alpha |=|\beta |\ge 2$,
\begin{equation}\label{eq:delic} \left \| \frac{z ^{\alpha}  \overline{z} ^{\beta}  z ^A \overline{
z} ^{B} }{|z _k|^2}\right \| _{L^1_t}\le \left \| z ^{ A}
\right  \| _{L^2_t} \left \| \frac{ z ^{  {B} }  z _\ell}{z
_k} \right \| _{L^2_t} \left \| \frac{z ^{ {\alpha}  }   \overline{z} ^{
{\beta}  }}{z _\ell z _k}\right \| _{L^\infty_t}\lesssim
C_2^2\epsilon ^{ |\alpha|+|\beta |}\le C_2^2\epsilon ^{4}.
\end{equation}
  Other contributions from \eqref{other_ta}
can be treated similarly, yielding \eqref{contrR2-}.

The main contribution to   the l.h.s. of \eqref{eq:FGR7} is originated from the following terms
\begin{equation}\label{eq:main cont}
    \|  \nu_j\frac{z^\mu
 \overline{ {z }}^ { {\nu} } }{\overline{z}_j} \langle g,
 {\mathbbm{i}} \sigma_3 G
_{\mu \nu }\rangle  \overline{\zeta} _j \| _{L^1_t}\le c_1 C_0 \epsilon ^2 \text{  for $(\mu, \nu) \in \mathbf{M}$}
\end{equation}
with $c_1$ a fixed constant. Indeed the term to bound equals
\begin{align*}   \nu_j z^\mu
 \overline{ {z }}^ { {\nu} }   \langle g,
 {\mathbbm{i}} \sigma_3 G
_{\mu \nu }\rangle  +  \nu_j\frac{z^\mu
 \overline{ {z }}^ { {\nu} } }{\overline{z}_j} \langle g,
 {\mathbbm{i}} \sigma_3 G
_{\mu \nu }\rangle  (\overline{\zeta} _j -\overline{z}_j).
\end{align*}
By Lemma \ref{lemma:bound g}, the 1st term has $L^1_t$ norm bounded by
\begin{align*}   \| G
_{\mu \nu } \| _{L^\infty _t  H^{0, s} }  \| z^\mu
 \overline{ {z }}^ { {\nu} } \| _{L^2_t}  \|  g \| _{L^2_t  H^{0,-s} } \le \| G
_{\mu \nu } \| _{L^\infty _t  H^{0, s} } C _0\epsilon c \epsilon \le c_1 C _0\epsilon ^2
\end{align*}
   for a fixed $c_1$. The 2nd term has $L^1_t$ norm bounded by the following, which yields \eqref{eq:main cont},
\begin{align*}  \| \nu_j\frac{z^\mu
 \overline{ {z }}^ { {\nu} } }{\overline{z}_j} \|_{
L^\infty_t}    \| G
_{\mu \nu } \| _{L^\infty _t  H^{0, s} }     \|  g \| _{L^2_t  H^{0,-s} }  \| \zeta -z\| _{L^2 _t}\lesssim \epsilon ^{ 4}.
\end{align*}
We estimated
the contribution  to the l.h.s. of \eqref{eq:FGR7} of   $\mathcal{E}_j$.
There are further terms in \eqref{eq:eq D} to estimate.  We claim
\begin{align}\label{eq:diff-1}  \|  (\mathfrak{L}_j(z,0)-\mathfrak{L}_j(z,f)) \overline{\zeta}_j   \|_{
L^\infty_t}  \lesssim \epsilon ^{ 4}.
\end{align}
A typical contribution to the l.h.s. is
\begin{align*}   \( g(\Pi (f))- g(\Pi (0)\) \frac{   \nu_j\overline{{z
}}^ {\nu +\beta }  }{  \overline{z}_j}  (\overline{z}_j  +(\overline{\zeta}_j-\overline{z}_j)) \text{ with } \alpha, \nu \in \mathbf{M}_0,
\end{align*}
with $g\in C^1(\R ^7, \C )$. We can bound its  $L^1_t$ norm using
\begin{align*}   \| f \| _{L^\infty H^1}^2  \| {z
} ^ {\nu } \|   _{L^2} \| {z
} ^ { \beta } \|   _{L^2} \lesssim \epsilon ^{ 4} \end{align*}
and using the argument that leads to \eqref{eq:delic}.  For the discussion of the
bound  for the contribution  originating from the $\sum _{l=1,..., \mathbf{n}}$  term in \eqref{eq:eq D}, which is also higher order,   see \cite{MR2843104}.

\qed

The 2nd term in the r.h.s. of \eqref{eq:FGR5}
equals, using  $G_{\mu \nu }^0= \overline{G^0}_{ \nu\mu} $,
\begin{equation} \label{eq:FGR8} \begin{aligned} & 2\sum_{\kappa  \in \mathfrak{K}  } \kappa
 \Re \left \langle R_{\im
\mathcal{H}_{p^1}}^+ (-\kappa)\sum_{ \substack{ \alpha\in \mathbf{M}_0 \, , \ \mathbf{e}\cdot \alpha =\kappa }}\zeta^{
\alpha } G_{ \alpha 0}^0, {\mathbbm{i}} \sigma_3\sum_{ \substack{\nu\in \mathbf{M}_0 \ , \ \mathbf{e}\cdot \nu =\kappa }}\overline{\zeta}^{ \nu } G^0_{0\nu
 } \right \rangle =\\&  {\pi}^{-1}\sum_{\kappa  \in \mathfrak{K} } \kappa \Re \left \langle R_{\im \mathcal{H}_{p^1}}^+ (- \kappa)\mathbf{G},
{\mathbbm{i}} \sigma_3\overline{\mathbf{G} }\right \rangle \text{ for }\mathbf{G}:=\sqrt{2\pi}\sum_{ \substack{ \alpha\in \mathbf{M}_0 \ , \ \mathbf{e}\cdot \alpha =\kappa }}\zeta^{ \alpha } G_{ \alpha 0}^0,
\end{aligned}
\end{equation}
where
$\mathfrak{K}=\{ k\in \R : \exists \ \nu\in \mathbf{M}_0 \text{ s.t. }\kappa =\mathbf{e}\cdot \nu\}$. Notice that $\kappa \in \mathfrak{K} \Rightarrow \kappa >\omega ^1.$

\noindent As in \cite[Lemma 10.5]{MR2805462},
there exist  ${L}_{ \alpha 0} \in W^{k,p}(\R ^3, \C ^4)$
 for all $k\in \R$ and $p\ge 1$
 s.t.
the r.h.s. of \eqref{eq:FGR8}
 is equal to
\begin{equation*} \label{eq:FGR81} \begin{aligned} & \sum_{\kappa  \in \mathfrak{K}} \kappa \Lambda(\kappa,\zeta) \text{ for }\Lambda(\kappa,\zeta)
=\frac{1}{\pi}\Re \left \langle R_{\im{\mathbbm{i}} \sigma_3 (-\Delta + \omega^1)}^+ (- \kappa)\mathbf{L}(\zeta),
{\mathbbm{i}} \sigma_3\overline{\mathbf{L} }\right \rangle \text{ and }
\mathbf{L}(\zeta):=\sqrt{2\pi}\sum_{ \substack{
\alpha\in \mathbf{M}_0
\\ \mathbf{e}\cdot \alpha =\kappa }}\zeta^{ \alpha } L_{ \alpha 0}^0.
\end{aligned}
\end{equation*}
We claim that each term in the above summation is non-negative. Observe that $ \Lambda(\kappa,\zeta)=\Lambda_1(\kappa,\zeta)+\Lambda_2 (\kappa,\zeta)$,
$\mathbf{L}(\zeta) = {^t (\mathbf{L}_1(\zeta),\mathbf{L}_2(\zeta))}$,
with
%\ac{like this???}
\begin{equation*}
 \Lambda_i(\kappa,\zeta)=  {\pi}^{-1}(-1)^{i+1}\Re \left \langle R_{\im{\mathbbm{i}} (-1)^{i+1} (-\Delta + \omega^1)}^+ (- \kappa)\mathbf{L}_i,
{\mathbbm{i}} \overline{\mathbf{L} }_i\right \rangle.
\end{equation*}
Introduce now
\begin{equation*}
 U=\frac{1}{\sqrt{2}}\begin{pmatrix} 1 &
1 \\
\im & -\im
 \end{pmatrix} \text{ such that }U^{-1}{\mathbbm{i}} U =-\im \sigma_3,
\end{equation*}
with $\sigma_3$ the Pauli matrix \eqref{pauli-matrices}.
%% that is acting on $\C^2$.
Taking the complex conjugation,
$\overline{U}^{-1}{\mathbbm{i}} \overline{U} =\im \sigma_3$.
Then, using ${^t \overline{U}}= {U}^{-1}$, we have, for $ U^{-1}\mathbf{L}_{i}= {^t (\mathbf{L}_{i1},\mathbf{L}_{i2})}$:
\begin{equation*} \begin{aligned} & \pi \Lambda_i(\kappa,\zeta)
= (-1)^{i+1}\Re \left \langle U^{-1} R_{\im{\mathbbm{i}} (-1)^{i+1} (-\Delta + \omega^1)}^+ (- \kappa)U U^{-1}\mathbf{L}_i,
\overline{U}^{-1} {\mathbbm{i}} \overline{U} \overline{U}^{-1}\overline{\mathbf{L} }_i\right \rangle \\& = (-1)^{i+1}\Re \left \langle R_{ (-1)^{i+1} \sigma_3(-\Delta + \omega^1)}^+ (- \kappa) U^{-1}\mathbf{L}_i,
\im \sigma_3 \overline{U}^{-1}\overline{\mathbf{L} }_i\right \rangle \\& = (-1)^{i+1}\Re \left \langle R_{ (-1)^{i+1} (-\Delta + \omega^1)}^+ (- \kappa) \mathbf{L}_{i1},
\im \overline{\mathbf{L} }_{i1}\right \rangle -(-1)^{i}\Re \left \langle R_{ (-1)^{i} (-\Delta + \omega^1)}^+ (- \kappa) \mathbf{L}_{i2},
\im \overline{\mathbf{L} }_{i2}\right \rangle.
\end{aligned}
\end{equation*}
Using Plemelj formula we have:
\begin{equation*} \begin{aligned} & \Lambda_1(\kappa,\zeta)= \left \langle \im \delta (\Delta- \omega^1+ \kappa) \mathbf{L}_{12},
\im \overline{\mathbf{L} }_{12}\right \rangle =- \left \langle \delta (\Delta- \omega^1+ \kappa) \mathbf{L}_{12},
 \overline{\mathbf{L} }_{12}\right \rangle \le 0;
\end{aligned}
\end{equation*}
\begin{equation*} \begin{aligned} & \Lambda_2(\kappa,\zeta)= \left \langle \im \delta (\Delta - \omega^1+ \kappa) \mathbf{L}_{21},
\im \overline{\mathbf{L} }_{21}\right \rangle =- \left \langle \delta (\Delta - \omega^1+ \kappa) \mathbf{L}_{21}, \overline{\mathbf{L} }_{21}\right \rangle \le 0.
\end{aligned}
\end{equation*}
The Fermi Golden Rule consists in two parts. The 1st part consists in
showing  that $ \Lambda ( \kappa , \zeta )$ are negative quadratic forms  for the vector
$(\zeta^\alpha )_{\alpha\in \mathbf{M}_0 \ s.t.\ \alpha\cdot \omega ^1=\kappa}$. This was proved here.  The 2nd part is that
the  $ \Lambda ( \kappa , \zeta )$ are strictly negative quadratic forms.  This   is  expected to be generically true    (as a similar statement was expected to be true in \cite{MR1334139,SW3}).    We don't know how to prove   this. For a proof on a different problem, see \cite[Proposition 2.2]{MR3053771}.
 For specific systems the strict negative condition   ought to be checked numerically. Here we assume it as an hypothesis:

\begin{itemize}
\item[(H9)] (\textit{Fermi Golden Rule})  the l.h.s. of \eqref{eq:FGR}, proved above to be negative, is strictly negative,
 that is
 for some fixed constants and
for any vector $\zeta \in \mathbb{C}^\mathbf{n}$ we have
\begin{equation} \label{eq:FGR} \begin{aligned} &
 \sum_{\in \mathfrak{K}} \kappa \Lambda(\kappa,\zeta)
 \approx -\sum_{
\alpha\in \mathbf{M}_0
} | \zeta^\alpha |^2.
\end{aligned}
\end{equation}

\end{itemize}

  By (H9) we have
\begin{equation} \label{eq:FGR10} \begin{aligned} &
2\sum_{l=1,...,\mathbf{n} }\mathbf{e }_l\Im \left (\mathcal{D}_l \overline{\zeta}_l
\right)\gtrsim \partial_t \sum_{l=1,...,\mathbf{n} }\mathbf{e }_l| \zeta_l|^2 + \sum_{
\alpha\in \mathbf{M}_0
} | \zeta^\alpha |^2.
\end{aligned}
\end{equation}
Then, for $t\in [0,T]$ and assuming Lemma~\ref{lemma:FGR1}, we have
\begin{equation*}
 \sum_{l=1,...,\mathbf{n} } \mathbf{e }_l | \zeta_l(t)|^2+ \sum_{
\alpha\in \mathbf{M}_0
 } \| \zeta^\alpha \|_{L^2(0,t)}^2\lesssim
\epsilon^2+ C_0\epsilon^2.
\end{equation*}
By \eqref{equation:FGR3} this implies $|z|_{L^\infty(0,t)}^2+\sum_{
\alpha\in \mathbf{M}_0
} \| z^\alpha \|
_{L^2(0,t)}^2\lesssim \epsilon^2+ C_0\epsilon^2$   and  yields Proposition~\ref{prop:est z}.
\qed

In the course of the proof we have shown that  $\| z^\alpha \|_{L^2(0,t)}^2\lesssim
 C_0^2\epsilon^2$ and \eqref{def:epsilon}  together imply
$\| z^\alpha \|_{L^2(0,t)}^2\lesssim C_0\epsilon^2$.
This means that we can take
$C_0\approx 1$.  With Corollary \ref{prop:conditional4.2} this completes the proof of Proposition \ref{prop:mainbounds}. \qed

\section{Proof of Theorem~\ref{theorem-1.1}}
\label{sec:end}

\begin{lemma}\label{lem:scattf} There is $f_+ \in H^1(\R^3, \C^4)$ such that
$f(t)$ from \eqref{Strichartzradiation} satisfies
\begin{equation}\label{eq:scattering1}
 \begin{aligned} & \lim_{t\to + \infty } \| f(t) - W(t,0)e^{ {\mathbbm{i}} \sigma_3 (-\Delta + \omega^{1}) t}f_+ \|_{H^1 }=0,
\end{aligned}
\end{equation}
where $W(t,s)$ is the fundamental solution from \eqref{eq:fu sol}.
\end{lemma}
\proof Starting from \eqref{eq:equation f} and using \eqref{eq:potlin}, we obtain
the following analogue of \eqref{eq:423}:
\begin{equation*} \begin{aligned} & f (t) = W(t,0)e^{ {\mathbbm{i}} \sigma_3 (-\Delta + \omega^{1}) t}f (0) \\& + \int_0^t
 e^{ {\mathbbm{i}} \sigma_3 (-\Delta + \omega^{1}) (t-t')}W(t,t')
\Big[ V f (t ') - \sum_{(\mu, \nu) \in \mathbf{M}}
z^\mu (t ') \overline{z}^\nu (t ')
 G_{\mu \nu}^0 + R_1(t ')+ R_2(t ') \Big]\,dt'.
\end{aligned}
 \end{equation*}
This implies, by standard arguments (cf. \cite[Sect. 11]{CM1}),
that $ W(0,t)e^{ {\mathbbm{i}} \sigma_3 (\Delta - \omega^{0}) t}f(t) \mathop{\longrightarrow}\limits\sb{t\to +\infty}f_+ $ in
$H^1(\R^3, \C^4)$. \qed

\textit{Completion of  the proof of Theorem~\ref{th:main}.}
Recall that expressing $u$ in terms of the coordinates in \eqref{eq:coordinate} we have
\begin{equation}\label{eq:u}
u(t)= e^{ -{\mathbbm{i}} \sigma_3\sum_{j=1}^{4}\tau_j'(t) \Diamond_j }
\big(\sqrt{1-| b'(t) |^2} + b '(t) \sigma_2 \mathbf{C}\big)
\big(\Phi_{ p' (t) } +P_{ p'(t)}r'(t)\big),
\end{equation}
where we denote by $ (p ',\tau ',b', r')$ the initial coordinates. Using the invariance $\Pi (u(t))= \Pi (u_0)$
we can express $ (p ', b')$ in terms of $r'$ obtaining the following:
\begin{equation}\label{eq:compl0}
\begin{aligned}&
 p'_j(t) = \Pi_j(u_0) - \Pi_j(r' (t)) + \resto^{1,2}_{\infty,\infty }\big(p^0_4,\Pi (r' (t)),r' (t)\big)
\quad\text{for $j=1,2,3,4$};
\\& b_R'(t) =({2p_4^0}) ^{-1} \Pi_{5}(r'(t))+ \resto^{2,0}_{\infty,\infty}\big( p^0_4,\Pi (r'(t))\big)
+\resto^{1,2}_{\infty,\infty} \big(p^0_4, \Pi (r'(t)), r'(t)\big);
\\&
 b_I' (t)= ({2p_4^0}) ^{-1}  \Pi_{6}(r'(t))+ \resto^{2,0}_{\infty,\infty} \big( p^0_4,\Pi (r'(t))\big)
+\resto^{1,2}_{\infty,\infty}\big(p^0_4,\Pi (r'(t)), r'(t)\big)
. \end{aligned} \end{equation}
 Furthermore we can express $r'$ in terms of the $(z,f)$ of the last coordinate system for $\ell = 2\mathbf{N}+1$
in Proposition~\ref{th:main}:
 \begin{equation}\label{eq:compl01}
\begin{aligned}  r'(t)= e^{{\mathbbm{i}} \sum  _{j=1}^{4}\sigma_3\resto^{0,2}_{k,m}
\big(p^0_4,\Pi (f(t)), z(t), f(t)\big)  \Diamond _j }& T(e^{ \sum  _{i=1}^{3}\resto^{0,2}_{k,m}(p^0_4, \Pi (f), z, f)
 {\mathbbm{i}} \sigma  _i}) \\&
\left(f(t)+ \mathbf{{S}}^{0,1}_{k,m}\big(p^0_4,\Pi (f(t)), z(t), f(t)\big)\right).
\end{aligned} \end{equation}
While the changes of coordinates in Lemma~\ref{lem:darboux} and in
 the normal forms   in Section~\ref{sec:Normal form}   involve loss of regularity of $f$,
in order to be differentiable so that the pullback of the symplectic forms makes sense,
nonetheless these maps are also continuous changes of coordinates inside in $H^1(\R^3, \C^2)$, see Lemma~\ref{lem:ODE} for $l=0$.
Notice that \eqref{Eq:NLS} leaves $\Sigma _k(\R^3, \C^2)$ invariant for any $k \in \N$ and that,  similarly, the system
leaves  $\C^{\mathbf{n}}\times (X_c\cap \Sigma _k(\R^3, \C^2))$ invariant.

\noindent By the well-posedness of \eqref{Eq:NLS} in $H^1(\R^3, \C^2)$ and of \eqref{eq:SystK1} in $\C^{\mathbf{n}}\times X_c$,
  a continuous change of coordinates  \eqref{eq:u}--\eqref{eq:compl01} maps solutions
 of \eqref{eq:SystK1} in $\C^{\mathbf{n}}\times X_c$ into solutions in $H^1(\R^3, \C^2)$ of \eqref{Eq:NLS},
 capturing the solutions of \eqref{Eq:NLS} in the statement of Theorem~\ref{th:main}. See also \cite[Sect. 8]{CM2}.

 \noindent By Lemma~\ref{lem:scattf} it is easy to conclude that $\resto^{0,2}_{k,m}\mathop{\longrightarrow}\limits\sb{t\to +\infty} 0 $ in $\R^7$ and $\mathbf{{S}}^{0,1}_{k,m}\mathop{\longrightarrow}\limits\sb{t\to +\infty} 0 $ in $ \Sigma_k(\R^3, \C^4)$ for the terms in \eqref{eq:compl01},
and that $\resto^{1,2}_{k,m}\mathop{\longrightarrow}\limits\sb{t\to +\infty} 0 $ for the terms in \eqref{eq:compl0}. Then for $1\le j\le 4$
we have
\begin{equation*}\begin{aligned}&
 \lim_{t\to +\infty} \Pi_j(r' (t)) = \lim_{t\to +\infty} \Pi_j(f (t))
 =\lim_{t\to +\infty}\Pi_j
\big(W(t,0)e^{\mathbbm{i} \sigma_3(-\Delta + \omega^{1}) t}f_+\big)
 = \Pi_j(f_+)
 \end{aligned}\end{equation*}
since $\Pi_j\big(W(t,0)e^{ {\mathbbm{i}} \sigma_3 (-\Delta + \omega^{1}) t}f_+\big)
= \Pi_j(f_+)$. Hence, since $p$ is characterized by the first four variables
(cf. \eqref{def p}), this defines $p_+$ in \eqref{eq:scattering}.

\noindent We consider a function $ g\in C^1(\R_+, \mathbf{G})$ such that
\begin{equation*}
 e^{ -{\mathbbm{i}} \sigma_3\sum_{j=1}^{4}\tau_j'(t) \Diamond_j }
\big(\sqrt{1-| b'(t) |^2} + b '(t) \sigma_2 \mathbf{C}\big) =T(g(t)).
\end{equation*}
By \eqref{eq:compl01} we have
\begin{equation}\label{expans1}\begin{aligned}&
 T(g(t))P_{ p'(t)}r'(t) = T(g(t)) e^{{\mathbbm{i}} \sigma_3\sum _{j=1}^4\resto^{0,2}_{k,m}   \Diamond _j }
 T(e^{\sum _{i=1}^3\resto^{0,2}_{k,m} {\mathbbm{i}}\sigma _i}) f + o_{\Sigma_{k}}(1),
 \end{aligned}\end{equation}
where $o_{\Sigma_{k}}(1) \mathop{\rightarrow}\limits\sb{t\to +\infty} 0 $ in $\Sigma_{k}(\R^3, \C^2)$.
We claim the following, with the proof in Appendix \ref{sec:flowspf}.
\begin{claim}\label{claim:represent}
\begin{equation}\label{eq:represent1}
    T(g(t)) e^{{\mathbbm{i}} \sigma_3\sum _{j=1}^4\resto^{0,2}_{k,m}   \Diamond _j }
 T(e^{\sum _{i=1}^3\resto^{0,2}_{k,m} {\mathbbm{i}}\sigma _i})= \widetilde{W}(0,t)
\end{equation}
 with $\widetilde{W}(t,s)$ the fundamental solution, in the sense of
\eqref{eq:fu sol}, of a system  of the form
\begin{equation}\label{eq:represent2}
 \dot u = {\mathbbm{i}} \sigma_3 \widetilde{\mathbf{v}} \cdot \Diamond u,
\qquad
\text{where}
\quad
\widetilde{\mathbf{v}} \cdot \Diamond = \sum_{j=1,...,7}{\mathbbm{i}} \sigma_3 \widetilde{\mathbf{v}}_j(t) \Diamond_j.
\end{equation}
\end{claim}

  Substituting \eqref{eq:coordinate-2} and \eqref{eq:compl01} into \eqref{Eq:NLS},
we get  for a  $G_1  \in C^0 (H^1(\R^3, \C^2), L^1(\R^3, \C^4))$ \begin{equation}\label{eq:compl3}
\begin{aligned} & \dot f=-{\mathbbm{i}} \sigma_3 \Delta f+ {\mathbbm{i}} \sigma_3 \widetilde{\mathbf{v}} \cdot \Diamond f +G_1(u),
\end{aligned}
\end{equation}
 while from \eqref{eq:equation f1}
 we have for a  $G_2  \in C^0 (H^1(\R^3, \C^2), L^1(\R^3, \C^4))$
 \begin{equation} \label{eq:equation f1-1}\begin{aligned} & \dot f =
-{\mathbbm{i}} \sigma_3 \Delta f +{\mathbbm{i}} \sigma_3 \omega^1 f
+ {\mathbbm{i}} \sigma_3 {\mathbf{v}} \cdot \Diamond f+G_2(u).
\end{aligned}\end{equation}
 The fact that  $G_1 ,\,G_2 \in C^0 (H^1(\R^3, \C^2), L^1(\R^3, \C^4))$ is rather simple.
   For example $G_2(u)$  is given by the sum of the r.h.s. of \eqref{eq:equation f1} with  a linear term
   $V _{\omega ^1}f$ where $V _{\omega ^1}\in S (\R ^3, M(\C ^4) )$ is the matrix valued function  in \eqref{eq:potlin}.  It is elementary to show that $u\to f$ is in
 $C^0 (H^1(\R^3, \C^2), L^2(\R^3, \C^4))$.

 \noindent The rest of $G_2(u)$ comes from the r.h.s. of \eqref{eq:equation f1},  obtained applying
 $\widehat{\nabla} _f$ to the terms $\mathbf{R}    | _{j=1} ^{3} $  in the expansion \eqref{eq:ExpH11}. It is elementary
  that this too is in $   C^0 (H^1(\R^3, \C^2), L^1(\R^3, \C^4))$.

 \noindent
  By comparing the equation  for $f$  with $G_1$ and the equation for $f$ with $G_2$,
 it follows that we necessarily have
 $ \widetilde{\mathbf{v}} \cdot \Diamond = \omega^1 + {\mathbf{v}} \cdot \Diamond $, see
  \cite[Lemma 13.8]{Cu3}. Hence, returning to \eqref{expans1}, we have
 \begin{equation} \begin{aligned}&
 T(g(t))P_{ p'(t)}r'(t) = \widetilde{W}(0,t) W(t,0) e^{ {\mathbbm{i}} \sigma_3 \omega^{1} t} e^{ - {\mathbbm{i}} \sigma_3 \Delta t}f_+ + o_{H^1}(1),
 \end{aligned}\nonumber \end{equation} for $W(t,0)$ defined by \eqref{eq:fu sol}
and where
\begin{equation*}
\partial_t (\widetilde{W}(0,t) W(t,0)e^{ {\mathbbm{i}} \sigma_3 \omega^{1} t})
= \widetilde{W}(0,t) {\mathbbm{i}} \sigma_3\left (({\mathbf{v}} - \widetilde{{\mathbf{v}}}) \cdot \Diamond
+\omega ^1 \right )  W(t,0) =0.
\end{equation*}
We conclude that there exists $g_0\in \mathbf{G}$ such that for $h_+= T(g_0) f_+$
one has
 \begin{equation} \begin{aligned}&
 T(g(t))P_{ p'(t)}r'(t) = e^{ - {\mathbbm{i}} \sigma_3 \Delta t} h_+ + o_{H^1}(1).
 \end{aligned}\nonumber \end{equation}
This completes the proof of \eqref{eq:scattering}.

Finally, we emphasize that the proof is predicated on the values  $\Pi _j (u_0)=p^0_j$ for $j\le 6$, with the coordinate
changes and the manifold $\mathscr{M}_{ 1}^{6}(p^0 )$ dependent on $p^0$.
However, since the symbols   $\mathcal{R}^{i,j}_{k, m}$ and $\mathbf{S}^{i,j}_{k, m}$
appearing in the coordinate changes depend continuously on $p^0$, the estimates
are uniform in $p^0$, as long as this is close enough to $p^1$. This completes the
proof  of Theorem~\ref{theorem-1.1}.\qed

\appendix
\section{Appendix.}\label{sec:flowspf}

Lemma \ref{lem:ODE}   is obtained  expressing $r$ in terms of $(z,f)$  from the following lemma, where we omit the dependence on the constant parameter $\Pi _4$.

\begin{lemma} \label{lem:ODEapp} For
$n,M,M_0, {s}, {s}',k,l\in  \N_0$ with $1\le l\le M$ such that \eqref{eq:index1} is satisfied, for $a\in A  $   a parameter,  with $A$ an open subset in $\R ^d$,  and for $\widetilde{\varepsilon}_0>0$,
consider
\begin{equation} \label{eq:ODEapp}\begin{aligned} &
 \dot r (t)
	= {\mathbbm{i}} \sigma_3\sum_{j\le 7}\resto^{0,M_0+1}_{n,M}(t,a, \Pi (r), r) \Diamond_j r +
\mathbf{S}^{i,M_0 }_{n,M}(t,a, \Pi (r), r), \end{aligned} \end{equation}
  Let $k\in \Z\cap [0,n-(l+1) ]$ and set for
$s''\ge 1$ and $\varepsilon >0$ \begin{equation} \label{eq:domain0app}
\begin{aligned} \U_{\varepsilon,k}^{ {s}''}
:= &\{ r \in T_{\Phi
_{p^1} }^{\perp_\Omega}\mathcal{M} \cap \Sigma_{ {s}''} \
 :\; \| r \|_{\Sigma_{-k }} + |\Pi (r)| \le \varepsilon \}.
\end{aligned} \end{equation}  Let $a_0\in A$. Then, for $\varepsilon >0$ small enough, \eqref{eq:ODEapp} 	 defines a flow
$\mathfrak{F}_{t } $
\begin{equation} \label{eq:ODE1app}
\begin{aligned} & \mathfrak{F}_{t } (r)
= e^{{\mathbbm{i}} \sigma_3  \sum _{j=1}^4\resto^{0,M_0+1 }_{ n- l-1, l}(t,a, \Pi (r), r)
   \Diamond _j }  T(  e^{   \sum _{i=1}^3  \resto^{0,M_0+1 }_{ n- l-1, l}(t,a, \Pi (r), r)
 {\mathbbm{i}} \sigma _i  }  )   \( r+\mathbf{{S}}^{i,M_0 }_{ n- l-1,l}(t,a, \Pi (r), r)\),
\end{aligned} \end{equation} where for
  and
 for $ \varepsilon_1> \varepsilon_2> 0$ sufficiently small we have
 \begin{equation} \label{eq:reg1app}\begin{aligned} &\mathfrak{F}_{t }
 \in C^l((-4,4)\times D _{\R^d}(a_0, \varepsilon_2 )\times \U_{\varepsilon_2,k}^{ {s}'}, \U_{\varepsilon_1,k}^{ {s}}
)
. \end{aligned} \end{equation}

 %Finally, we have \begin{equation}\label{eq:ODE11}\begin{aligned} & \mathfrak{F}_{t } (e^{{\mathbbm{i}} \sigma_3 \tau \cdot\Diamond } u) = e^{{\mathbbm{i}} \sigma_3 \tau \cdot \Diamond }\mathfrak{F}_{t } (u).\end{aligned} \end{equation} 	
 \end{lemma}
{\it Proof (sketch) } While the statement is the same of \cite[Lemma 3.8]{Cu0} and \cite[Lemma 3]{MR3053771},   we have to deal with
 operators $ \Diamond _j$  for $j=5,6,7$  which  don't commute.

For $\xi \in \mathbf{su}(2)$  and $q\in \R^4$ we consider
 $S:=   e^{-\mathbbm{i} \sigma _3 \sum _{j=1}^{4} q_j\Diamond _j } T ( e^{-\xi})r $, for $T$ the representation
 in \eqref{eq:rep1}.
   It is elementary that  for some $F_j\in C ^\infty$ we have
\begin{equation}\label{eq:momenta}\begin{aligned} &
   \Pi _j(r)  = \Pi _j(S)   \text  { for $j=1,2,3,4$ and } \\& \Pi _j(r)  =
   \Pi _j(S)+F_j(\xi,\Pi _k(S)| _{k=5} ^7) \text  { for $j=5,6,7$   }
\end{aligned} \end{equation}
where $F_j(0,* )\equiv 0\equiv F_j( *,0 )$ for any $*$ and where for $j=5,6,7$
the above equality is obtained proceeding like in Lemma  \ref{lem:lambdas1}. Then expressing the coefficients
of \eqref{eq:ODE1app}  in terms of the new variables, we have new coefficients
\begin{equation} \nonumber \begin{aligned} &
 \mathfrak{D}(t, \xi, \varrho , S):=\\&
	  e^{-\mathbbm{i} \sigma _3 \sum _{j=1}^{4} q_j\Diamond _j }T( e^{-\xi})
\mathbf{S}^{i,M_0 }_{n,M} \left (  t,  \varrho _l | _{l=1}^{4}, \varrho _l | _{l=5} ^7+
 F_l(\xi ,\varrho  _k | _{k=5} ^7)  | _{l=5} ^7
  ,e^{ \mathbbm{i} \sigma _3 \sum _{j=1}^{4} q_j\Diamond _j } T( e^{ \xi} )
   S   \right )  , \\& \mathfrak{A}_j(t, \xi, \varrho , S):= \resto^{0,M_0+1}_{n,M}
   \left (t,\varrho _l | _{l=1}^{4}, \varrho _l | _{l=5} ^7+ F_l(\xi ,\varrho  _k | _{k=5} ^7)
    | _{l=5} ^7 ,e^{ \mathbbm{i} \sigma _3 \sum _{j=1}^{4} q_j\Diamond _j } T( e^{ \xi} )  S    \right  ) .
\end{aligned}   \end{equation}
Notice that for $0\le \ell \le M$ we have
\begin{equation} \nonumber \begin{aligned} &
   \mathfrak{D} (t, \xi, \varrho , S)= \mathbf{S}^{i,M_0 }_{n-\ell,\ell}(t, \xi, \varrho , S)
    \text{ and }  \mathfrak{A}_j(t, \xi, \varrho , S)=\resto^{0,M_0+1}_{n-\ell,\ell}(t, \xi, \varrho , S).
\end{aligned}   \end{equation}
Then  consider the following system which we explain below:
\begin{align} &
    \dot S
	= \mathfrak{D}(t, \xi, \varrho , S)     ,\nonumber
 \\& \dot q _j = \mathfrak{A}_j(t, \xi, \varrho , S)   \text{  for $j=1,2,3,4$, with $q_j(0)=0$,}\nonumber
 \\& \sum _{k=1}^{\infty}
  \frac{1}{k!} \left ( \text{ad}(\xi) \right )  ^{k-1} \dot \xi  =   \sum_{i=1}^{3}
  \mathfrak{A}_j(t, \xi, \varrho , S)   {\mathbbm{i}} \sigma _i
  \text{    with $\xi(0)=0$,} \label{eq:par} \\& \dot \varrho _j= \langle   S ,   \Diamond _j
  \mathfrak{D}(t , \xi, \varrho, S)
   \rangle   \text{  for $j=1,2,3,4$ and } \nonumber\\& \dot \varrho _j= \langle   S ,
   \Diamond _j\mathfrak{D}(t, \xi, \varrho , S) \rangle-\partial _{\xi} F_j(\xi ,\varrho _k| _{k=5} ^7)
   \dot \xi -\sum _{l=5,6,7}
   \partial _{\varrho _l} F_j(\xi ,\varrho _k| _{k=5,6,7}  )  \dot \varrho _l   \text{  for $j=5,6,7$.}\nonumber
\end{align}
We explain now the above equations.
The  2nd and 3rd line are  defined in order to simplify the
equation for $S$. Indeed, when we substitute $S$ in the equation of $r$ we get
 \begin{align*} & \partial _t r   =
\partial _t(e^{ \mathbbm{i} \sigma _3 \sum _{j=1}^{4} q_j\Diamond _j } T ( e^{\xi})S)=
e^{ \mathbbm{i} \sigma _3 \sum _{j=1}^{4}  q_j\Diamond _j } T ( e^{\xi}) \left (
\mathbbm{i} \sigma _3 \sum _{j=1}^{4} \dot q_j\Diamond _j  S +T ( e^{-\xi}) \partial _t (T ( e^{\xi}))S + \dot S \right )=\label{large}\\&
 e^{ \mathbbm{i} \sigma _3 \sum _{j=1}^{4} q_j\Diamond _j }  \left (
  T ( e^{\xi}){\mathbbm{i}} \sigma_3\sum_{j=1} ^{4}\mathcal{A }_j(t, \Pi (r), r) \Diamond_j S
  +{\mathbbm{i}} \sigma_3\sum_{j=5} ^{7}\mathcal{A }_j(t, \Pi (r), r) \Diamond_j T ( e^{\xi})S\right) +
\mathcal{D}(t, \Pi (r), r).\nonumber
\end{align*}
By the choice made in the 2nd line of \eqref{eq:par} the summations over $j=1,2,3,4$ cancel out.
We will show that also the summations over $j=5,6,7$  cancel out. By the Baker--Campbell--Hausdorff formula, see
  \cite[p.15]{rossmann},      we have
\begin{equation}   \label{eq:lieder} \begin{aligned} &  \partial _t e^{\xi}=\left ( \sum _{k=1}^{\infty}
  \frac{1}{k!} \left ( \text{ad}(\xi) \right )  ^{k-1} \dot \xi \right ) e ^{\xi}
 \text{ where } \text{ad}(\xi): \mathbf{su}(2)\to \mathbf{su}(2) \text{ is }\text{ad}(\xi)\vartheta:= [\xi,\vartheta].
\end{aligned}   \end{equation}
So, for $\uno _{\C^2}$ the   unit element in  $\mathbf{SU}(2)$, we have
\begin{equation}  \label{eq:lieder1} \begin{aligned} &
\partial _t (T ( e^{\xi})) = dT(\uno _{\C^2}) \left ( \sum _{k=1}^{\infty}
  \frac{1}{k!} \left ( \text{ad}(\xi) \right )  ^{k-1} \dot \xi \right ) T (e ^{\xi}).
\end{aligned}   \end{equation}
  On the other hand,  by \eqref{eq:charge2} and \eqref{diffT} we have
\begin{equation} \nonumber \begin{aligned} &
\sum_{j=5,6,7} \mathcal{A }_j  {\mathbbm{i}} \sigma_3 \Diamond_j T (e ^{\xi}) =
\sum_{i=1,2,3} \mathcal{A }_{i+4} dT(\uno _{\C^2})  ( \mathbbm{i}\sigma _{i} )   T (e ^{\xi}).
\end{aligned}   \end{equation}
So the 3rd equation in \eqref{eq:par} yields the cancelation of these terms.  Hence we conclude that the 1st equation in \eqref{eq:par} is true.

\noindent We also derive equations for
$\varrho _j$ by  differentiating $\partial
_t\Pi _j(S)$ and by substituting $\Pi _j(S)$ with $\varrho _j$.

 \noindent Solving the last equation  in \eqref{eq:par}   in terms of  $\dot \varrho _j| _{j=5}^{7}$ and replacing in the last equation
 $\dot \xi$  by means   of the 3rd equation, we obtain  for $1\le \ell \le M$
\begin{equation} \label{eq:par1} \begin{aligned} &
    \dot S
	= \mathbf{S}^{i,M_0 }_{n-\ell,\ell}(t,\xi, \varrho ,  S)     ,
 \\& \dot q _j = \resto^{0,M_0+1}_{n-\ell,\ell}(t,\xi, \varrho ,   S)   \text{  for $j=1,2,3,4$, with $q_j(0)=0$,}
 \\& \dot \xi =\resto^{0,M_0+1}_{n-\ell,\ell}(t,\xi,  \varrho ,  S)
  \text{    with $\xi(0)=0$,}\\& \dot \varrho _j= \resto^{0,M_0+1}_{n-\ell-1,\ell}(t,\xi, \varrho ,  S)
     \text{  for $j=1,...,7$.}
\end{aligned}   \end{equation}
Taking as initial conditions $(r,0,0, \Pi (r))$, by elementary arguments, see
  \cite[Lemma 3.8]{Cu0}, we get from \eqref{eq:par1}  a flow
\begin{equation} \label{eq:par1} \begin{aligned}
      S(t)
&	= r+\int _0^t\mathbf{S}^{i,M_0 }_{n-\ell-1,\ell}(t',   \Pi (r) ,  r)dt'   =
r+ \mathbf{S}^{i,M_0 }_{n-\ell-1,\ell}(t ,   \Pi (r) ,  r)  ,
 \\    q _j (t)&= \int _0^t\resto^{0,M_0+1}_{n-\ell-1,\ell}(t',   \Pi (r) ,  r)dt'  =\resto^{0,M_0+1}_{n-\ell-1,\ell}(t ,  \Pi (r) ,  r)  \text{  for $j=1,2,3,4$,}
 \\    \xi (t)&=\sum _{i=1}^{3} \int _0^t\resto^{0,M_0+1}_{n-\ell-1,\ell}(t',   \Pi (r) ,  r)dt' \mathbbm{i} \sigma _i=
 \sum _{i=1}^{3} \resto^{0,M_0+1}_{n-\ell-1,\ell}(t ,  \Pi (r) ,  r) \mathbbm{i} \sigma _i
 ,\\ \Pi _j(S(t))&=\Pi _j(r) + \int _0^t\resto^{0,M_0+1}_{n-\ell-1,\ell}(t',   \Pi (r) ,  r)dt'\\& =\Pi _j(r)
 +  \resto^{0,M_0+1}_{n-\ell-1,\ell}(t ,   \Pi (r) ,  r)
     \text{  for $j=1,...,7$.}
\end{aligned}   \end{equation}
In view of \eqref{eq:momenta}  we get also
\begin{equation} \label{eq:par2} \begin{aligned} &
        \Pi _j(r(t))=\Pi _j(r) + \resto^{0,M_0+1}_{n-\ell-1,\ell}(t ,   \Pi (r) ,  r)
     \text{  for $j=1,...,7$ .}
\end{aligned}   \end{equation}
This ends the proof of the parts of  Lemma \ref{lem:ODE}  which are not the same of \cite[Lemma 3.8]{Cu0}.

\qed

\begin{lemma} \label{lem:ODEbisapp} Consider two systems for $\ell =1,2$:
\begin{equation*} \begin{aligned} &
 \dot r (t)
	= {\mathbbm{i}} \sigma_3\sum_{j=1,...,7} \mathcal{A }_j^{(\ell)}(t, \Pi (r), r) \Diamond_j r +
\mathcal{D}^{(\ell)}(t, \Pi (r), r),
\end{aligned} \end{equation*}
with the hypotheses of Lemma~\ref{lem:ODEapp} satisfied, and suppose that
\begin{equation*} \begin{aligned} &
 \mathcal{D}^{(1)}(t,\Pi(r),r)-\mathcal{D}^{(2)}(t,\Pi(r),r)=\mathbf{S}^{0,M_0+1}_{n,M}(t,\Pi(r),r).
\end{aligned} \end{equation*}
Let $r \mapsto r^t _{(\ell)} $ with $\ell =1,2$ be the  flow for each of the two
systems. Then for $ {s}, {s}'$ as in Lemma \ref{lem:ODEapp}
\begin{equation*} \begin{aligned} &
 \| r^1_{(1)}- r^1_{(2)}\|_{\Sigma_{- {s}' }} \le C \| r\|_{\Sigma_{ - {s} }}^{M_0+1}.
	 \end{aligned} \end{equation*}
\end{lemma}

\proof The proof is elementary. We consider
\begin{equation*} \begin{aligned} \sum_{\ell =1,2}(-)^\ell  \frac{d}{dt} r_{(\ell)}^t &= \sum_{\ell =1,2} (-)^\ell {\mathbbm{i}} \sigma_3 \resto^{0,M_0+1}_{n,M}(t, \Pi (r_{(\ell)}^t), r_{(\ell)}^t)\cdot \Diamond r_{(\ell)}^t\\& +\underbrace{\sum_{\ell =1,2}(-)^\ell\mathcal{D}^{(\ell)}(t, \Pi (r_{(2)}^t), r_{(2)}^t)}_{S^{0,M_0+1 }_{n,M}(t, \Pi (r_{(2)}^t), r_{(2)}^t)}+ \sum_{\ell =1,2}(-)^\ell\mathcal{D}^{(1)}(t, \Pi (r_{(\ell)}^t), r_{(\ell)}^t).
	 \end{aligned} \end{equation*}
Then for $\mathbf{x}_\ell ^t :=(\Pi (r_{(\ell)}^t), r_{(\ell)}^t) $
\begin{equation*} \begin{aligned} & \| r_{(2)} ^t- r_{(1)} ^t \|_{\Sigma_{-\mathbf{s} ' }}
\le \sum_\ell \int_0^t \| r_{(\ell)}  ^{t '} \|_{\Sigma_{- {s} }}^{M_0+2}\,dt'
 + \int_0^t \| r_{(2)} ^{t '}\|_{\Sigma_{- {s} }}^{M_0+1}\,dt'
 \\&+ \int_0^t
\int_0^1 \| \partial_{\Pi (r)} \mathcal{D}^{(1)}(t', \mathbf{x}_1^{t '}
+ \tau (\mathbf{x}_2^{t '}-\mathbf{x}_1^{t '})) \|_{\Sigma_{- {s} }}
|\Pi (r_{(2)} ^{t '}) - \Pi (r_{(1)} ^{t '}) |\,dt'
 \\& + \int_0^t
\int_0^1 \| \partial_{r} \mathcal{D}^{(1)}(t', \mathbf{x}_1^{t '}
+ \tau (\mathbf{x}_2^{t '}-\mathbf{x}_1^{t '})) \|_{\Sigma_{- {s} } \to \Sigma_{- {s} }}
\| r_{(2)} ^{t '} - r_{(1)}^{t '} \|_{\Sigma_{- {s} }}\,dt'.
	 \end{aligned} \end{equation*}
 Since there is a fixed $C>0$ such that
 \begin{equation*} \begin{aligned} & \| r_{(\ell)} (t')\|_{\Sigma_{- {s} }}\le
  C \| r \|_{\Sigma_{- {s} }} \text{ from \eqref{eq:ODE1}},
\\& |\Pi (r_{(2)} ^{t '}) - \Pi (r_{(1)} ^{t '}) | \le \| r \|_{\Sigma_{- {s} }}^{M_0+1}
 \text{ from the previous one and \eqref{eq:ODE1}}
 \\& \| \partial_{r} \mathcal{D}^{(1)}(t, \Pi, \varrho, r)
 \|_{\Sigma_{-\mathbf{s} } \to \Sigma_{- {s} } } \le C \| r \|_{\Sigma_{- {s} }}^{M_0-1},
 \\& \| \partial_{\varrho} \mathcal{D}^{(1)}(t, \Pi, \varrho, r) \|_{\Sigma_{- {s} } }
  \le C \| r \|_{\Sigma_{- {s} }}^{M_0 } \quad ,
\end{aligned} \end{equation*}
where the last inequalities follow from \eqref{eq:opSymb}, for some fixed constant $C>0$  we obtain
\begin{equation*}
\| r_{(2)} ^{t  }- r_{(1)}^{t  } \|_{\Sigma_{- {s} ' }}
%% \lesssim
\le
C\left(
t\| r \|_{\Sigma_{- {s} }}^{M_0+1} + \| r \|_{\Sigma_{- {s} }}^{M_0-1} \int_0^t
\| r_{(2)} ^{t '} - r_{(1)} ^{t '} \|_{\Sigma_{- {s} }}\,dt
\right),
\end{equation*}
 for $t\in[0,1]$,
which by Gronwall's inequality yields \eqref{eq:ODE1bis}.
\qed

\bigskip

\bigskip

\textit{Proof of {Claim} \ref{claim:represent}.}
Let $\mathbf{g} =\R ^4 \times  \mathbf{su}(2)$ be the Lie Algebra of $\mathbf{G}$.
 We can assume that the inverse of the  l.h.s. of \eqref{eq:represent1}
is equal to $  e ^{{\mathbbm{i}} \sigma_3 \sum _{j=1}^{4} X_j (t)\Diamond _j }T(e ^{\xi (t)}) $ with
$X \in C^ 1(\R _+, \R ^4) $  and $\xi \in  C^ 1(\R _+, \mathbf{su}(2)) $.
Then, by  \eqref{eq:lieder1},
 for $u(t):=e ^{{\mathbbm{i}} \sigma_3 \sum _{j=1}^{4} X_j (t)\Diamond _j }T(e ^{\xi (t)})u_0 $
we have
\begin{equation}   \nonumber \begin{aligned} &
\dot u(t) = {\mathbbm{i}} \sigma_3 \sum _{j=1}^{4} \dot X_j (t)\Diamond _ju(t) +
dT(\uno _{\C^2}) \left ( \sum _{k=1}^{\infty}
  \frac{1}{k!} \left ( \text{ad}(\xi (t) ) \right )  ^{k-1} \dot \xi (t)\right ) u (t).
\end{aligned}   \end{equation}
 We set $\widetilde{\mathbf{v}} _j(t)= \dot X_j (t)$ for $j\le 4$  and, exploiting that
 ${\mathbbm{i}} \sigma _l | _{l=1}^{3}$ is a basis of $\mathbf{su}(2)$,
  we define $\widetilde{\mathbf{v}} _j (t)| _{j=5}^{7}$  by
 \begin{equation*}
    \sum _{l=1}^3 \widetilde{\mathbf{v}}_{l+3}(t) {\mathbbm{i}} \sigma _l =  \sum _{k=1}^{\infty}
  \frac{1}{k!} \left ( \text{ad}(\xi (t) ) \right )  ^{k-1} \dot \xi (t).
 \end{equation*}
 Then we conclude that \eqref{eq:represent2}  it true
 for this choice of $u(t)$ and of $\widetilde{\mathbf{v}} _j (t)| _{j=1}^{7}$.
Then $u(t) = \widetilde{{W}} (t,0)u_0$   and $\widetilde{{W}} (0,t)=\widetilde{{W}} ^{-1} (t,0) $  is s.t.
 equality \eqref{eq:represent1}  it true. This yields Claim \ref{claim:represent}.\qed

\section{Appendix.}\label{sec:pullb}

Lemma \ref{lem:ODE1} can be obtained from the following lemma, expressing $r$ in terms
of $(z,f)$ and omitting again the dependence of the symbols on $\Pi _4$, which has constant value.

\begin{lemma} \label{lem:ODE1app} Consider
 $\mathfrak{F} =\mathfrak{F}^1 \circ \cdots \circ \mathfrak{F}^L$
 with $ \mathfrak{F}^j= \mathfrak{F}^j_{t=1}$ transformations as in Lemma~\ref{lem:ODEapp}
 on the manifold
   $\mathscr{M}_{ 1}^{6}(p^0 )  $.
 Suppose that for any $ \mathfrak{F}^j$ the $M_0 $ in Lemma~\ref{lem:ODEapp}
 equals $m_j$,
 where $1= m_1\le...\le m_L$ with the constant $i$ in
Lemma~\ref{lem:ODE} (ii)   equal to $1$
when $m_j=1$.
Fix   $M, k$ with $n_1\gg k\ge \mathbf{N}_0$ ($n_1$ picked in Lemma \ref{lem:modulation}).
Then there is a $ n=n(L,M,k)$ such that
if the  assumptions of Lemma~\ref{lem:ODE}
apply to each of operators $ \mathfrak{F}^j$  for  $(M,n)$,
  there exist
$\underline{{{\psi}}}(\varrho)\in C^\infty$ with $\underline{{\psi}}(\varrho)=O(|\varrho|^2)$
 and a small $\varepsilon >0$
 such that in $\U^{s}_{\varepsilon,k}$  for $s\ge n-(M+1)$      we have the expansion
\begin{align} \label{eq:back1app} &
K \circ \mathfrak{F}
= \underline{{{\psi}}} (\Pi (r)  ) + {2}^{-1}\Omega (\mathcal{H}_{ p} P_{ p}r, P_{ p}r) + \resto^{1,2}_{k,M} + E_P(P_{ p}r)+\textbf{R}'',
\\& \textbf{R}'':= \sum_{d=2,3,4}
\langle B_{d } (\Pi (r) , r), (P_{ p}r)^{ d} \rangle
 +\int_{\mathbb{R}^3}B_5 (x, \Pi (r) , r, r(x)) (P_{ p}r)^{ 5}(x)\,dx,
\nonumber
\end{align}
with:
\begin{itemize}
\item $B_2(0,0)=0$;
\item
$B_{d}(\varrho, r) \in C^{M } (\U_{-k},
\Sigma_k (\mathbb{R}^3, B (
(\mathbb{R}^{4 })^{\otimes d},\mathbb{R}))) $, $2\le d \le 4$,
with $\U_{-k}\subset \R^7 \times (T_{\Phi_{p^1} }^{\perp_\Omega}\mathcal{M} \cap \Sigma_{-k})$
an open neighborhood of $(0,0)$;
\item for
$ \zeta \in \mathbb{R}^{4 }$ $(\varrho, r) \in \U_{-k}$ we have for $i\le M$
\begin{equation*} \label{eq:B5a}
\begin{aligned} &
\| \nabla_{ r,\varrho,\zeta}^iB_5(\varrho,r,\zeta)\|_{\Sigma_k(\mathbb{R}^3,
B ((\mathbb{R}^{4 })^{\otimes 5},\mathbb{R})} \le C_i.
\end{aligned} \end{equation*}
\end{itemize}
\end{lemma}
\proof The proof is in \cite{Cu0}, but  we sketch it.
First of all, by \eqref{eq:reg1app} we have, for $k\le  n-L(M+1) $,
\begin{equation}\label{eq:lossreg}
 \U_{\varepsilon_{L+1},k}^{n-(M+1)}
 \stackrel{\mathfrak{F}^L }{\longrightarrow} \U_{\varepsilon_L,k}^{n- 2(M+1)}...\stackrel{\mathfrak{F}^2 }{\longrightarrow} \U_{\varepsilon_2,k}^{n- L (M+1)}
 \stackrel{\mathfrak{F}^1 }{\longrightarrow} \U_{\varepsilon_1,k}^{n-(L+1)(M+1)}\subset \U_{\varepsilon_1,k}^{ k+3} \subset \U_{\varepsilon_1,k}^{\mathbf{N}_0} ,
\end{equation}
where each map is $C^M$ if we  pick $n_1\ge n = n (L,M,k) := k+3 +(L+1)(M+1)$ and then
we get $\mathfrak{F}\in C^M (\U_{\varepsilon_{L+1},k}^{n-(M+1)},  \U_{\varepsilon_1,k}^{ k+3}) $.

\noindent By \eqref{eq:ODE1app}, the $r$--th component of $\mathfrak{F}$ is of the form
\begin{equation} \label{eq:structf}\begin{aligned} & \mathfrak{F} (
 \varrho, r) = e^{{\mathbbm{i}} \sigma_3\sum _{j=1}^4\resto^{1,1 }_{ k+3,M}(\varrho, r)  \Diamond _j }
  T(  e^{   \sum _{i=1}^3  \resto^{1,1 }_{ k+3,M}(\varrho, r)
 {\mathbbm{i}} \sigma _i  }  )   (r+
\mathbf{{S}}^{1,1 }_{ k+3,M}
 (\varrho, r)).
\end{aligned} \end{equation}
Then by $[\Diamond _j,\Diamond _k]=0$ for all $k$ if $j\le 4$ we have
 \begin{equation*}
  \Pi _j(r) | _{j=1}^{4}\circ \mathfrak{F} =   \Pi _j(r+
\mathbf{{S}}^{1,1 }_{ k+3,M}
 \Pi (r), r)) | _{j=1}^{4} = \Pi (r)| _{j=1}^{4}
 +\mathcal{R}^{1, 2}_{k+2,M}(\Pi (r), r) .
 \end{equation*}
 From \eqref{eq:variables} we have \begin{equation*} \label{eq:ODE2}\begin{aligned} & p\circ \mathfrak{F} = p +\mathcal{R}^{1, 2}_{k+2,M} \text{ and so } \Phi_{ p}\circ \mathfrak{F} = \Phi_{p} + \mathbf{S}^{1, 2}_{k+2,M}.
\end{aligned} \end{equation*}
 Then we have
\begin{equation*} \begin{aligned}
& E(u\circ \mathfrak{F})=
E\Big(e^{ -{\mathbbm{i}} \sigma_3\sum_{j=1}^{4}\tau_j \Diamond_j } (\sqrt{1-| b |^2} + b \sigma_2 \mathbf{C}) (\Phi_{ p } +P_{ p}r)\circ \mathfrak{F}\Big)
\\& = E((\Phi_{ p } +P_{ p}r)\circ \mathfrak{F}) = E(\Phi_{ p} + \mathbf{S}^{1, 2}_{k+2,M} +P_{ p }(e^{{\mathbbm{i}} \sigma_3\resto^{1,1 }_{ k+2,M} \cdot \Diamond } (r+
\mathbf{{S}}^{1,1 }_{ k+2,M}))
\\&= E(\Phi_{ p } +P_{ p}r +\mathbf{S}^{1, 2}_{k+2,M} + P_{ p}\mathbf{S}^{1, 1}_{k+2,M}),
\end{aligned} \end{equation*}
where we use   the commutation (for the proof see \cite[Lemma 4.1]{Cu0})
 \begin{equation*}\begin{aligned}
&
 [P_{ p },e^{{\mathbbm{i}} \sigma_3\sum _{j=1}^4\resto^{1,1 }_{ k+3,M}  \Diamond _j }
  T(  e^{   \sum _{i=1}^3  \resto^{1,1 }_{ k+3,M} ) \mathbbm{i} \sigma _i  } ]r\\& =
   [e^{{\mathbbm{i}} \sigma_3\sum _{j=1}^4\resto^{1,1 }_{ k+3,M}  \Diamond _j }
  T(  e^{   \sum _{i=1}^3  \resto^{1,1 }_{ k+3,M} ) \mathbbm{i} \sigma _i  },\widehat{P}_{ p }]r=\mathbf{S}^{1, 2}_{k+2,M}.
  \end{aligned}
 \end{equation*}
 We get similarly for $1\le j \le 4$
 \begin{equation*} \begin{aligned} \Pi _j(u\circ \mathfrak{F}) | _{j=1}^{4} &= \Pi _j(\Phi_{ p } +P_{ p}r +\mathbf{S}^{1, 2}_{k+2,M} + P_{ p}\mathbf{S}^{1, 1}_{k+2,M}) | _{j=1}^{4}= \Pi_j (\Phi_{ p } +P_{ p}r)| _{j=1}^{4} + \mathcal{R}^{1,2}_{k,m}.
\end{aligned} \end{equation*}
Then
 \begin{equation}
\label{k-f-u}
 \begin{aligned} &
 {K}(\mathfrak{F}(u))= E (\Phi_{ p } +P_{ p}r +\mathbf{S}^{1, 2}_{k+2,M} + P_{ p}\mathbf{S}^{1, 1}_{k+2,M})- E\left (\Phi_{ {p}^0}\right) \\&
\qquad\qquad\qquad
-\sum _{j\le 4}(\lambda_j(p) + \resto^{1,2}_{k+2,m}) \left (\Pi_j (\Phi_{ p } +P_{ p}r) + \mathcal{R}^{1,2}_{k,m} -\Pi_j\left (\Phi_{ {p}^0}\right) \right).
\end{aligned}
\end{equation}
Like in \cite[Lemma 4.3]{Cu0}, we set
\begin{equation*}
\Psi = \Phi_{ p } +\mathbf{S}^{1, 2}_{k+2,M} + P_{ p}\mathbf{S}^{1, 1}_{k+2,M};
\end{equation*}
we need to analyze
$E(\Psi +P_{ p}r)$ which we break into
(cf. \eqref{eq:energyfunctional})
\[
E(\Psi +P_{ p}r)
=
E_P(\Psi +P_{ p}r)
+
E_K(\Psi +P_{ p}r).
\]
It is also shown in \cite[Lemma 4.3]{Cu0} that
\begin{equation} \label{eq:back2}
\begin{aligned} &
E_P(\Psi +P_{ p}r)= E_P(\Psi) + E_P(P_{ p}r)
+ \text{terms that can be incorporated into $ \textbf{R}''  $}
\\&
\qquad\qquad\qquad
+\sum_{ j=0,1} \int_{\R^3}dx \int_{[0,1]^2}
 \frac{t^j}{j!} (\partial_t^{j+1})\at{t=0} \partial_s [B (|s \Psi+tP_{ p}r |^2) ]\,dt\,ds.
\end{aligned}
\end{equation}
The second line of \eqref{eq:back2} equals
\begin{equation} \label{eq:back3}
\begin{aligned}
&
\int_{\R^3}dx \int_{[0,1]^2}\,dt\,ds
\sum_{ j=0,1}\frac{t^j}{j!}(\partial_t^{j+1})\at{t=0} \partial_s \ \Big \{ \ \ B (|s \Phi_{ p}+tP_{ p}r |^2)
+
\\
&
\qquad\qquad
+ \int_0^1 d\tau \partial_\tau [B (|s (\Phi_{ p}
+\tau (\mathbf{S}^{1, 2}_{k+2,M} + P_{ p}\mathbf{S}^{1, 1}_{k+2,M}) +tP_{ p}r |^2) ] \ \ \Big \}.
\end{aligned}
\end{equation}
The contribution from the last line of \eqref{eq:back3} can be incorporated into $\textbf{R}''+\resto^{1,2}_{k,m}$.
Notice that from the
$j=0$ term in the first line of \eqref{eq:back3} we get
\begin{equation} \label{eq:back31-1}
2\int_{\R^3}dx \int_{0}^1\,ds
 \partial_s [B' (|s \Phi_p |^2) s \Phi_{ p}\cdot P_{ p}r ]
=2\int_{\R^3}dx B' (| \Phi_{ p} |^2) \Phi_{ p}\cdot P_{ p}r
=\langle \nabla E_P(\Phi_{ p}), P_{ p} r\rangle.
\end{equation}
The
$j=1$ term in the first line of \eqref{eq:back3} is
${2^{-1}}\langle\nabla^2 E_P(\Phi_{p})P_{p}r,P_{p}r\rangle$;
thus,
\begin{eqnarray}\label{e-p}
E_P(\Psi +P_{ p}r)
=
E_P(\Psi) + E_P(P_{ p}r)
+
\langle \nabla E_P(\Phi_{ p}), P_{ p} r\rangle
+
{2}^{-1}\langle\nabla^2 E_P(\Phi_{p})P_{p}r,P_{p}r\rangle
+\textbf{R}''+\resto^{1,2}_{k,m}.
\end{eqnarray}
Then,
\begin{eqnarray} \label{eq:back31}
E_K(\Psi +P_{p}r)
= E_K(\Psi) -\langle \Delta \Phi_{ p}, P_{ p}r \rangle
+\overbrace{\langle -\Delta (\mathbf{S}^{1, 2}_{k+2,M}
+P_{ p}\mathbf{S}^{1, 1}_{k+2,M}), P_{ p}r \rangle}^{\resto^{1,2}_{k,m}}
+E_K(P_{ p}r).
\end{eqnarray}
Using \eqref{eq:energyfunctional}, \eqref{eq:Omega}, \eqref{eq:eqphi} and the fact that
${\mathbbm{i}} \sigma_3\lambda (p)\cdot \Diamond \Phi_{ p }\in T_{ \Phi
_{ p} }\mathcal{M}$, see \eqref{eq:Tmanifold},  we have
\begin{equation*}\begin{aligned}
 \langle -\Delta \Phi_{ p}+\nabla E_P(\Phi_{ p}), P_{ p} r\rangle &= \langle \nabla E (\Phi_{ p }), P_{ p} r\rangle = -\Omega ({\mathbbm{i}} \sigma_3 \nabla E (\Phi_{ p }), P_{ p} r)
\\&
= -\Omega ({\mathbbm{i}} \sigma_3 \lambda (p) \cdot \Diamond \Phi_{ p }, P_{ p} r) =0.
\end{aligned}
\end{equation*}
Adding \eqref{e-p} and \eqref{eq:back31}
and using
the cancellation of the sum of the
second term in the right-hand side of \eqref{eq:back31}
with the term \eqref{eq:back31-1}
which follows from the above relation,
we arrive at
\begin{eqnarray}\label{e-p-k}
E(\Psi +P_{ p}r)
=
E(\Psi) + E(P_{ p}r)
+
{2}^{-1}\langle\nabla^2 E_P(\Phi_{p})P_{p}r,P_{p}r\rangle
+\textbf{R}''+\resto^{1,2}_{k,m},
\end{eqnarray}
where we used \eqref{eq:energyfunctional}.
From \eqref{eq:eqphi},
\begin{equation} \label{lastform}\begin{aligned}
E (\Psi)
&= E (\Phi_{ p }) +\overbrace{\langle \nabla E(\Phi_{ p }), P_{ p }\mathbf{S}^{1,1}_{k+2,M}\rangle}^{0} +
 \overbrace{\langle \nabla E(\Phi_{ p }), \mathbf{S}^{1,2 }_{k+2,M} \rangle }^{ \resto^{1,2 }_{k+2,M}}+ \resto^{1,2 }_{k,M}
%%\\&
= E (\Phi_{ p }) + \resto^{1,2 }_{k,M},
\end{aligned}
\end{equation}
where the $\resto^{1,2 }_{k,M}$ in the right-hand side
is absorbed by $\resto^{1,2 }_{k,M}$ in \eqref{eq:back1app}.

We have
\begin{equation}\label{previousformula} \begin{aligned} -\lambda (p)\cdot \Pi (\Phi_{ p } +P_{ p}r) &= -\lambda (p)\cdot \Pi (\Phi_{ p })-\lambda (p)\cdot \Pi (P_{ p}r)
- \langle \lambda (p)\cdot \Diamond \Phi_{ p },P_{ p}r \rangle
\\&
= -\lambda (p)\cdot \Pi (\Phi_{ p })-\lambda (p)\cdot \Pi (P_{ p}r),
\end{aligned}
\end{equation}
where we used $\langle \lambda (p)\cdot \Diamond \Phi_{ p },P_{ p}r \rangle = \Omega (-{\mathbbm{i}} \sigma_3\lambda (p)\cdot \Diamond \Phi_{ p },P_{ p}r)=0$.

Substituting \eqref{e-p-k} (where we apply \eqref{lastform}) and \eqref{previousformula}
into \eqref{k-f-u}, we have:
\begin{equation} \nonumber \begin{aligned}   {K}(\mathfrak{F}(u))
&=
E(\Phi_p) + E(P_{ p}r)
+
{2}^{-1}\langle\nabla^2 E_P(\Phi_{p})P_{p}r,P_{p}r\rangle
-E(\Phi_{p^0})\\&  -\lambda (p)\cdot \Pi (\Phi_{ p })-\lambda (p)\cdot \Pi (P_{ p}r)
+
\lambda(p)\cdot\Pi(\Phi_{ {p}^0})
+\textbf{R}''+\resto^{1,2}_{k,m}.   .\end{aligned}
\end{equation}
By \eqref{def-d-p},
$d(p)  =E(\Phi_{ p }) - \lambda (p) \cdot \Pi (\Phi_{ p })$.
Then we have
\begin{align} \label{eq:back4}
E(\Phi_{ p }) - E(\Phi_{ {p}^0 })- \lambda (p)\cdot (\Pi (\Phi_{ p }) - \Pi (\Phi_{ {p}^0 }))
&=
d(p)- d({p}^0)- (\lambda ({p}^0)-\lambda (p))\cdot {p}^{0 }
\nonumber
\\
&=
K (\Phi_{ p })= O((\Pi _j(r)| _{j=1}^{4})^2)+\resto^{2,2}_{\infty, \infty},
\end{align}
where $O((\Pi  _j(r)| _{j=1}^{4})^2)$ is $\underline{{{\psi}}} (\Pi  _j(r) )$  in
\eqref{eq:back1app} and $\resto^{2,2}_{\infty, \infty}$  is absorbed inside $\resto^{1,2}_{k,M}$.
Thus,
\[
{K}(\mathfrak{F}(u))
=
 \underline{{{\psi}}} (\Pi   (r) )
+E(P_{ p}r)
+
{2}^{-1}\langle\nabla^2 E_P(\Phi_{p})P_{p}r,P_{p}r\rangle
-\lambda (p)\cdot \Pi (P_{ p}r)
+\textbf{R}''+\resto^{1,2}_{k,m}.
\]
Breaking $E(P_{ p}r)=E_P(P_{ p}r)+E_K(P_{ p}r)$
and using the relation
\begin{equation*}
\begin{aligned}&
 {2}^{-1} \langle \nabla^2 E_P(\Phi_{ p })P_{ p } r, P_{ p } r\rangle + E_K(P_{ p }r)-\lambda (p) \cdot
 \Pi (P_{ p } r) \\
&\qquad
= {2}^{-1}\langle (\nabla^2 E (\Phi_{ p }) -\lambda (p) \cdot \Diamond)P_{ p } r, P_{ p } r\rangle = {2}^{-1}\Omega (\mathcal{H}_{ p} P_{ p}r, P_{ p}r),
\end{aligned}
\end{equation*}
we arrive at the conclusion of the lemma.
\qed

The following   is an elementary consequence of Lemma~\ref{lem:ODE1app}
and is proved in \cite[Lemma 4.4]{Cu0}.
\begin{lemma}
 \label{lem:back11app} Under the hypotheses and notation of Lemma~\ref{lem:ODE1},
 for $\textbf{R}'$ like $\textbf{R}''$, for $ {{{\psi}}} \in C^\infty (\R^4, \R)$
 with $ {{\psi}} (\varrho) =O(| \varrho |^2)$, we have
\begin{align} \label{eq:back11app} &
K \circ \mathfrak{F}= {{\psi}} (\Pi  _j(r)| _{j=1,...,4}) +  {2}^{-1}\Omega (\mathcal{H}_{ {p}^1} r, r)
 + \resto^{1,2}_{k,m} 	 +E_P(r)+\textbf{R}',
\\& \textbf{R} ':= \sum_{d=2,3,4}
\langle B_{d } (\Pi (r), r), r^{ d} \rangle
 +\int_{\mathbb{R}^3}
B_5 (x, \Pi (r), r, r(x)) r)^{ 5}(x)\,dx, \nonumber
\end{align}
 the $B_d$ for $2\le d\le 5$ with similar properties of the functions in Lemma~\ref{lem:back}.

 \end{lemma} \proof
The proof, for whose details we refer to \cite{Cu0}, is obtained by writing
\begin{equation*}
 P_{ p} r=r+(P_{ p} -P_{ {p}^1}) r=r +\mathbf{S}^{1,1}_{\infty,\infty}
\end{equation*}
and substituting $P_{ p} r=r +\mathbf{S}^{1,1}_{\infty,\infty}$ inside
\eqref{eq:back1app}. That from $E_P(P_{ p}r)+\textbf{R}''$ in
\eqref{eq:back1app} we   obtain a term which is contained in
$\resto^{1,2}_{k,m} +E_P(r)+\textbf{R}'$ in \eqref{eq:back11app} is elementary and is discussed in \cite{Cu0}.
 We have
\begin{equation} \label{eq:back012} \begin{aligned}& {2}^{-1}\Omega (\mathcal{H}_{ p} P_{ p}r, P_{ p}r) = {2}^{-1} \langle -\Delta P_{ p}r, P_{ p}r \rangle - \lambda (p) \cdot \Pi (P_{ p}r) + {2}^{-1} \langle \nabla^2 E_P (\Phi_{ p}) P_{ p}r, P_{ p}r \rangle.\end{aligned}
\end{equation}
Then
 \begin{equation} \label{eq:back13}
\begin{aligned} \langle -\Delta P_{ p}r, P_{ p}r \rangle
&= \langle -\Delta r, r \rangle +\resto^{1,2}_{k,m}, \quad
\Pi (P_{ p}r) = \Pi (r) +\resto^{1,2}_{k,m},
\\[1ex]
\langle \nabla^2 E_P (\Phi_{ p})P_{ p}r, P_{ p}r \rangle & =
\langle \nabla^2 E_P (\Phi_{ {p}^1}) r, r \rangle + \resto^{1,2}_{k,m}
+ \langle (\nabla^2 E_P (\Phi_{ p})-\nabla^2 E_P (\Phi_{ {p}^1})) r, r \rangle,
\\[1ex]
\lambda (p) &=\lambda ({p}^1)+ \resto^{1,0}_{\infty, \infty}(\Pi  _j(r)| _{j=1}^{4}) + \resto^{1,2}_{k,m},
\end{aligned} \nonumber
\end{equation}
where for the last line we considered \eqref{eq:variables} which implies
\begin{equation*}
\begin{aligned} & p=\Pi -\Pi (r) + \resto^{1,2}_{\infty,\infty }
\end{aligned} \end{equation*}
and where $\resto^{1,0}_{\infty, \infty}(\Pi (r))$ is smooth in the argument and is $O(|\Pi (r)|)$.

\noindent Then we conclude that the right hand side of \eqref{eq:back012} is
\begin{equation} \label{eq:back120} \begin{aligned} &
 \overbrace{2^{-1} \langle (-\Delta - \lambda ({p}^1) \cdot \Diamond +\nabla^2 E_P (\Phi_{ {p}^1})) r, r \rangle }^{2^{-1}\Omega (\mathcal{H}_{ {p}^1} r, r) }+\resto^{2,0}_{\infty, \infty}(\Pi  _j(r)| _{j=1}^{4}) +\resto^{1,2}_{k,m} \\&+ 2^{-1}\langle (\nabla^2 E_P (\Phi_{ {p} })-\nabla^2 E_P (\Phi_{ {p}^1})) r, r \rangle,
\end{aligned}
\end{equation}
where the last term can be absorbed in the $d=2$ term of $\textbf{R}'$ by \eqref{eq:variables}.
Setting $\psi (\varrho) =\underline{\psi }(\varrho) +\resto^{2,0}_{\infty,\infty} (\varrho)$ with the $\resto^{2,0}_{\infty,\infty}$
in \eqref{eq:back120}, we get the desired result.
\qed

\section*{Acknowledgments}
  S.C., and a visit of A.C. in Trieste,
was   funded by    grants  FRA 2015 and    FRA 2018 from the University of Trieste.

\bibliographystyle{amsplain}

\end{document}